\documentclass[journal,9pt,onecolumn]{IEEEtran}
\usepackage[noadjust]{cite}
\usepackage[pdftex]{graphicx}
\usepackage{caption}
\usepackage{subcaption}
\usepackage{amsmath,amsfonts,amssymb}
\usepackage{dsfont}
\usepackage{algorithmic}
\usepackage{array}
\usepackage{stfloats}
\usepackage{bigints}
\usepackage{url}
\usepackage{multirow}
\usepackage{tikz}
\newtheorem{proposition}{\textbf{Proposition}}
\newtheorem{fact}{\textbf{$\blacktriangleright$ Fact}}
\newtheorem{example}{\textbf{$\blacktriangleright$ Example}}

\begin{document}
\title{Gaussian distributions on Riemannian symmetric spaces\,:
statistical learning with structured covariance matrices}

\author{Salem~Said, Hatem~Hajri, Lionel Bombrun, Baba~C.~Vemuri}

\maketitle

\abstract
The Riemannian geometry of covariance matrices has been essential to several successful applications, in computer vision, biomedical signal and image processing, and radar data processing. For these applications, an important ongoing challenge is to develop Riemannian-geometric tools which are adapted to structured covariance matrices. The present paper proposes to meet this challenge by introducing a new class of probability distributions, \textit{Gaussian distributions of structured covariance matrices}. These are Riemannian analogs of Gaussian distributions, which only sample from covariance matrices having a preassigned structure, such as complex, Toeplitz, or block-Toeplitz. The usefulness of these distributions stems from three features\,: (1) they are completely tractable, analytically or numerically, when dealing with large covariance matrices, (2) they provide a statistical foundation to the concept of structured Riemannian barycentre (\textit{i.e.} Fr\'echet or geometric mean), (3) they lead to efficient statistical learning algorithms, which realise, among others, density estimation and classification of structured covariance matrices. The paper starts from the observation that several spaces of structured covariance matrices, considered from a geometric point of view, are Riemannian symmetric spaces. Accordingly, it develops an original theory of Gaussian distributions on Riemannian symmetric spaces, of their statistical inference, and of their relationship to the concept of Riemannian barycentre. Then, it uses this original theory to give a detailed description of Gaussian distributions of three kinds of structured covariance matrices, complex, Toeplitz, and block-Toeplitz. Finally, it describes algorithms for density estimation and classification of structured covariance matrices, based on Gaussian distribution mixture models.  
\section{Introduction} \label{sec:intro}
Data with values in the space of covariance matrices have become unavoidable in many applications (here, the term \textit{covariance matrix} is taken to mean \textit{positive definite matrix}). For example, these data arise as ``tensors'' in mechanical engineering and medical imaging~\cite{moakher2,MR,pennec2}, or as ``covariance descriptors'' in signal or image processing and computer vision~\cite{congedo,arnaudon,yang}\!\!~\cite{image1,image2,tuzel,caseiro}. Most often, the first step in any attempt to analyse these data consists in computing some kind of central tendancy, or average\,: the Riemannian barycentre~\cite{moakher2}\cite{pennec2}\cite{jeurissurvey} (\textit{i.e.} Fr\'echet or geometric mean), the Riemannian median~\cite{arnaudon}\cite{yang}, the log-det or Stein centre~\cite{stein}\cite{anoop}, \textit{etc}. Therefore, an important objective is to provide a rigorous statistical framework for the use of these various kinds of averages.  In the case of the Riemannian barycentre, this objective was recently realised, through the introduction of Gaussian distributions on the space of covariance matrices~\cite{vemuri}\cite{said}. 

Gaussian distributions on the space of covariance matrices were designed to extend the useful properties of normal distributions, from the usual context of Euclidean space, to that of the Riemannian space of covariance matrices equipped with its affine-invariant metric described in~\cite{pennec2}. Accordingly, they were called ``generalised normal distributions'' in \cite{vemuri} and ``Riemannian Gaussian distributions'' in~\cite{said}. The most attractive feature of these Gaussians distributions is that they verify the following property, which may also be taken as their definition,
\begin{equation} \label{eq:intro1}
\mbox{defining property of Gaussian distributions\,: \textit{maximum likelihood is equivalent to Riemannian barycentre}}
\end{equation}
This is an extension of the defining property of normal distributions on Euclidean space, dating all the way back to the work of Gauss\!~\cite{ghistory}. Indeed, in a Euclidean space, the concept of Riemannian barycentre reduces to that of arithmetic average. By virtue of Property (\ref{eq:intro1}), efficient algorithms were developed, which use Riemannian barycentres in carrying out prediction and regression~\cite{vemuri}\cite{regress}, as well as classification and density estimation~\cite{rosu,eusipco,gretsi}, for data in the space of covariance matrices. 

The present paper addresses the problem of extending the definition of Gaussian distributions, from the space of covariance matrices to certain spaces of structured covariance matrices. The motivation for posing this problem is the following. The structure of a covariance matrix is a highly valuable piece of information. Therefore, it is important to develop algorithms which deal specifically with data in the space of covariance matrices, when these are known to possess a certain structure. To this end, the paper introduces a new class of probability distributions, \textit{Gaussian distributions of structured covariance matrices}. These are Gaussian distributions which only sample from covariance matrices having a certain preassigned structure. Using these new distributions, 
it becomes possible to develop new versions of the algorithms from~\cite{vemuri}\cite{regress}\cite{rosu,eusipco,gretsi}, which deal specifically with structured covariance matrices.

In the following, Gaussian distributions of structured covariance matrices will be given a complete description for three structures\,: complex, Toeplitz, and block-Toeplitz. This is achieved using the following observation, which was made in~\cite{vemuri}\cite{said}. The fact that it is possible to use Property (\ref{eq:intro1}), in order to define Gaussian distributions on the space of covariance matrices, is a result of the fact that this space, when equipped with its affine-invariant metric, is a Riemannian symmetric space of non-positive curvature. Starting from this observation, the following new contributions are developed below, in Sections \ref{sec:gaussian}, \ref{sec:matrices} and \ref{sec:mixture}. These employ the mathematical background provided in Section \ref{sec:geo} and Appendix \ref{app:roots}. 

\indent --- \textbf{In Section \ref{sec:gaussian}\,:} Gaussian distributions on Riemannian symmetric spaces of non-positive curvature are introduced. If $\mathcal{M}$ is a Riemannian symmetric space of non-positive curvature, then a Gaussian distribution $G(\bar{x},\sigma)$ is defined on $\mathcal{M}$ by its probability density function, with respect to the Riemannian volume element of $\mathcal{M}$,
\begin{equation} \label{eq:intro2}
\mbox{p.d.f. of a Gaussian distribution\,:}\hspace{0.5cm}  p(x|\,\bar{x},\sigma) =  \frac{1}{Z(\sigma)} \times \exp \left[-\,\frac{d^{\,\scriptscriptstyle 2}(x,\bar{x})}{2\sigma^{\,\scriptscriptstyle 2}}\right] \hspace{0.25cm} \bar{x} \in \mathcal{M}\,,\, \sigma > 0
\end{equation}
where $d : \mathcal{M} \times \mathcal{M} \rightarrow \mathbb{R}_+$ denotes Riemannian distance in $ \mathcal{M}$. In all generality, as long as $\mathcal{M}$ is a Riemannian symmetric space of non-positive curvature, Proposition \ref{prop:z} provides the exact expression of the normalising factor $Z(\sigma)$, and Proposition \ref{prop:sample} provides an algorithm for sampling from the Gaussian distribution $G(\bar{x},\sigma)$. Also, Propositions \ref{prop:mle} and \ref{prop:asymp} characterise the maximum likelihood estimates of the parameters $\bar{x}$ and $\sigma$. In particular, Proposition \ref{prop:mle} shows that definition (\ref{eq:intro2}) implies Property (\ref{eq:intro1}). 

\indent --- \textbf{In Section \ref{sec:matrices}\,:} as special cases of the definition of Section \ref{sec:gaussian}, Gaussian distributions of structured covariance matrices are described for three spaces of structured covariance matrices\,: complex, Toeplitz, and block-Toeplitz. Each one of these three spaces of structured covariance matrices becomes a Riemannian symmetric space of non-positive curvature, when equipped with a so-called Hessian metric\!~\cite{barbaresco,barbaresco1,jeuris1,jeuris2}, 
\begin{equation} \label{eq:intro3}
  \mbox{ Riemannian metric } = \mbox{ Hessian of entropy} 
\end{equation}
where the entropy function, for a covariance matrix $x$, is $H(x) = -\log \det(x)$, as in the case of a multivariate normal model~\cite{cover}. For each one of the three structures, the metric (\ref{eq:intro3}) is given explicitly, the normalising factor $Z(\sigma)$ is computed analytically or numerically, and an algorithm for sampling from a Gaussian distribution with known parameters is provided, (these algorithms are described in Appendix \ref{app:sample}).   

\indent --- \textbf{In Section \ref{sec:mixture}\,:} efficient statistical learning algorithms are derived, for classification and density estimation of data in a Riemannian symmetric space of non-positive curvature, $\mathcal{M}$. These algorithms are based on the idea that any probability density on $\mathcal{M}$, at least if sufficiently regular, can be approximated to any desired precision by a finite mixture of Gaussian densities, 
\begin{equation} \label{eq:intro4}
\mbox{mixture of Gaussian densities\,:}\hspace{0.5cm}  p(x) = \sum^{\scriptscriptstyle K}_{\scriptscriptstyle \kappa \,= 1} \omega_{\scriptscriptstyle \kappa} \times p(x | \,\bar{x}_{\scriptscriptstyle \kappa},\sigma_{\scriptscriptstyle \kappa})  
\end{equation}
where $\omega_{\scriptscriptstyle 1}, \ldots, \omega_{\scriptscriptstyle K} > 0$ satisfy  $\omega_{\scriptscriptstyle 1} + \ldots + \omega_{\scriptscriptstyle K} = 1$ and where each density $p(x | \,\bar{x}_{\scriptscriptstyle \kappa},\sigma_{\scriptscriptstyle \kappa})$ is given by (\ref{eq:intro2}). A BIC criterion for selecting the mixture order $K$, and an EM (expectation-maximisation) algorithm for estimating mixture parameters $\lbrace (\omega_{\scriptscriptstyle \kappa},\bar{x}_{\scriptscriptstyle \kappa},\sigma_{\scriptscriptstyle \kappa})\rbrace$ are provided. Applied in the setting of Section \ref{sec:matrices}, these lead to new versions of the algorithms in~\cite{said}\cite{rosu,eusipco,gretsi}, which deal specifically with structured covariance matrices, and which perform, with a significant reduction in computational complexity, the same tasks as the recent non-parametric methods of~\cite{chevallier3,chevallier1,chevallier2}. 

Table 1. below summarises some of the new contributions of Sections \ref{sec:gaussian} and \ref{sec:matrices}. The table is intended to help readers navigate directly to desired results. \\[0.03cm]
\indent The main geometric ingredient in the defining Property (\ref{eq:intro1}) of Gaussian distributions is the Riemannian barycentre. Recall that the Riemannian barycentre of a set of points $x_{\scriptscriptstyle 1},\ldots, x_{\scriptscriptstyle N}$ in a Riemannian manifold $\mathcal{M}$ is any global minimiser $\hat{x}_{\scriptscriptstyle N}$ of the variance function
\begin{equation} \label{eq:intro5}
\mbox{variance function\,:} \hspace{0.5cm} E_{\scriptscriptstyle N}(x) = \, \frac{1}{N} \, \sum^{\scriptscriptstyle N}_{\scriptscriptstyle n=1} 
d^{\scriptscriptstyle \,2}(x,x_{\scriptscriptstyle n})
\end{equation}
where $\hat{x}_{\scriptscriptstyle N}$ exists and is unique whenever $\mathcal{M}$ is a Riemannian space of non-positive curvature~\cite{afsari}\cite{moakher1}. The Riemannian barycentre of covariance matrices arises when $\mathcal{M} \, = \, \mathcal{P}_{\scriptscriptstyle n\,}$, the space of $n \times n$ covariance matrices, equipped with its affine-invariant metric. Its existence and uniqueness follow from the fact that $\mathcal{P}_{\scriptscriptstyle n}$ is a Riemannian symmetric space of non-positive curvature~\cite{vemuri}\cite{said}. 

Due to its many applications, the Riemannian barycentre of covariance matrices is the subject of intense ongoing research~\cite{jeurissurvey}\cite{congedo1}\cite{congedo2}. Recently, attention has been directed to the problem of computing the Riemannian barycentre of structured covariance matrices~\cite{jeuris1}\cite{jeuris2}\cite{bini}. This problem arises because if $Y_{\scriptscriptstyle 1},\ldots,Y_{\scriptscriptstyle N} \in \mathcal{P}_{\scriptscriptstyle n}$ are covariance matrices which have some additional structure, then their Riemannian barycentre $\hat{Y}_{\scriptscriptstyle N}$ need not have this same structure~\cite{bini}. In order to resolve this issue, two approaches have been proposed, as follows.

Note that a space of structured covariance matrices is the intersection of $\mathcal{P}_{\scriptscriptstyle n}$ with some real vector space $V$ of $n \times n$ matrices. This is here denoted $V\mathcal{P}_{\scriptscriptstyle n} = V \, \cap \, \mathcal{P}_{\scriptscriptstyle n\,}$. The vector space $V$ can be that of matrices which have complex, Toeplitz, or block-Toeplitz structure, \textit{etc}. The first approach to computing the Riemannian barycentre $\hat{Y}_{\scriptscriptstyle N}$ of $Y_{\scriptscriptstyle 1},\ldots,Y_{\scriptscriptstyle N} \in V\mathcal{P}_{\scriptscriptstyle n}$ was proposed in~\cite{bini}. It 
defines $\hat{Y}_{\scriptscriptstyle N}$ as a minimiser of the variance function (\ref{eq:intro5}) over $V\mathcal{P}_{\scriptscriptstyle n\,}$, rather than over the whole of $\mathcal{P}_{\scriptscriptstyle n\,}$. This is an extrinsic approach, since it considers $V\mathcal{P}_{\scriptscriptstyle n}$ as a subspace of $\mathcal{P}_{\scriptscriptstyle n\,}$. By contrast, the second approach to computing $\hat{Y}_{\scriptscriptstyle N}$ is an intrinsic approach, proposed in~\cite{barbaresco,barbaresco1,jeuris1,jeuris2}. It considers $V\mathcal{P}_{\scriptscriptstyle n}$ as a ``stand-alone'' Riemannian manifold $\mathcal{M} = V\mathcal{P}_{\scriptscriptstyle n}$ equipped with the Hessian metric (\ref{eq:intro3}). Then, $\hat{Y}_{\scriptscriptstyle N}$ is given by the general definition of the Riemannian barycentre, applied to this Riemannian manifold $\mathcal{M}$.  

For the present paper, the intrinsic approach to the definition of the Riemannian barycentre of structured covariance matrices is preferred. In this approach, a space $V\mathcal{P}_{\scriptscriptstyle n}$ of structured covariance matrices is equipped with the Hessian metric (\ref{eq:intro3}). If $V\mathcal{P}_{\scriptscriptstyle n}$ is the space of complex, Toeplitz, or block-Toeplitz covariance matrices, then it becomes a Riemannian symmetric space of non-positive curvature, when equipped with this metric\!~\cite{barbaresco,barbaresco1,jeuris1,jeuris2}. In turn, this allows for the definition of Gaussian distributions, based on Section \ref{sec:gaussian}. In general, note that $V\mathcal{P}_{\scriptscriptstyle n}$ is an open convex cone in the real vector space $V$. If $V\mathcal{P}_{\scriptscriptstyle n}$ is moreover a homogeneous and self-dual cone, then it becomes a Riemannian symmetric space of non-positive curvature, when equipped with a Hessian metric of the form (\ref{eq:intro3}), which is given in~\cite{barbaresco}\cite{faraut}. Thus, the definition of Gaussian distributions, based on  Section \ref{sec:gaussian}, extends to any space $V\mathcal{P}_{\scriptscriptstyle n}$ of structured covariance matrices, whenever $V\mathcal{P}_{\scriptscriptstyle n}$ is a homogeneous and self-dual cone.

Throughout the following, the normalising factor $Z(\sigma)$ appearing in (\ref{eq:intro2}) will play a major role. Proposition \ref{prop:z} establishes its two main properties\,: (i) although (\ref{eq:intro2}) contains two parameters $\bar{x}$ and $\sigma$, this normalising factor does not depend on $\bar{x}$ but only on $\sigma$, (ii) when expressed in terms of the ``natural parameter'' $\eta = -1/2\sigma^{\scriptscriptstyle \,2}$, this normalising factor is a strictly log-convex function. Proposition \ref{prop:mle} shows these two properties are crucial ingredients in the definition of Gaussian distributions. Precisely, if $x_{\scriptscriptstyle 1},\ldots, x_{\scriptscriptstyle N}$ are independent samples from a Gaussian distributions $G(\bar{x},\sigma)$ given by (\ref{eq:intro2}), then the log-likelihood function for the parameters $\bar{x}$ and $\sigma$ is given by
\begin{equation} \label{eq:intro6}
\mbox{log-likelihood function of a Gaussian distribution\,:} \hspace{0.25cm} L(\bar{x},\sigma) \, = \, -N\, \log Z(\sigma) \, - \, \frac{1}{\mathstrut 2\sigma^{\scriptscriptstyle 2}} \, \sum^{\scriptscriptstyle N}_{\scriptscriptstyle n=1}
d^{\scriptscriptstyle \,2}(\bar{x},x_{\scriptscriptstyle n})
\end{equation}
Then, (i) guarantees that the defining Property  (\ref{eq:intro1}) is verified. Indeed, since the first term in (\ref{eq:intro6}) does not depend on $\bar{x}$, maximisation of $L(\bar{x},\sigma)$ with respect to $\bar{x}$ is equivalent to minimisation of the sum of squared distances appearing in its second term. However, this is the same as the variance function (\ref{eq:intro5}). On the other hand, (ii) guarantees that maximisation of $L(\bar{x},\sigma)$ with respect to $\sigma$ is equivalent to a one-dimensional convex optimisation problem, (expressed as a Legendre transform in equation (\ref{eq:legendre1})). Note the exact expression of $Z(\sigma)$ is given by formulae (\ref{eq:zint}) and (\ref{eq:zprod}) in the Proof of Proposition \ref{prop:z}. In~\cite{paolo}, a Monte Carlo integration technique was designed which uses these formulae to generate the graph of $Z(\sigma)$. In Section \ref{sec:matrices}, this technique is shown to effectively compute $Z(\sigma)$ when the space $\mathcal{M}$ is a space of large structured covariance matrices. 

In conclusion, it is expected that Gaussian distributions on Riemannian symmetric spaces of non-positive curvature will be valuable tools in designing statistical learning algorithms which deal with structured covariance matrices. Such algorithms would combine efficiency with reduced computational complexity, since the main requirement for maximum likelihood estimation of Gaussian distributions is computation of Riemannian barycentres, a task for which there exists an increasing number of high-performance routines. As discussed in Section \ref{sec:mixture}, this is the essential advantage of the definition of Gaussian distributions adopted in the present paper, over other recent definitions, also used in designing statistical learning algorithms, such as those considered in~\cite{chevallier1,chevallier2,chevallier3}. 
\vspace{0.14cm}


\begin{center}
\begin{tabular}{lcccc} 
\multicolumn{5}{c}{\underline{Table 1.\,: Gaussian distributions on symmetric spaces and spaces of structured covariance matrices}}\\[0.3cm]
 & \;\;Notation\;\; & \;\;Riemannian metric\;\; & \phantom{xx} \;\;Normalising factor $Z(\sigma)$\;\;  & \;\;Sampling algorithm\;\; \\[0.15cm]
\multirow{2}{*}{\parbox{26mm}{Section \ref{sec:gaussian}\,:\\ symmetric space $\mathcal{M}$}}\phantom{xx} &  & & 
\multirow{2}{*}{\parbox{26mm}{$\phantom{xxx}$(\ref{eq:zint}) and (\ref{eq:zprod}) \\ $\phantom{xx}$in Proposition \ref{prop:z}}}
& \multirow{2}{*}{\parbox{26mm}{$\phantom{mm}$Algorithm 1.\\ based on Proposition \ref{prop:sample}}} \\
& & & 
 & 
\\[0.15cm]
\multirow{2}{*}{\parbox{26mm}{Paragraph \ref{subsec:2t2}\,: $n \times n$\\complex covariance \\matrices }} &  & & 
\multirow{2}{*}{\parbox{26mm}{$\phantom{mm}$Monte Carlo\\ integration from (\ref{eq:hnzint})}}
& \multirow{2}{*}{\parbox{26mm}{$\phantom{mm}$Algorithm 2.\\$\phantom{m}$ in Appendix \ref{app:sample}}} \\
& $\mathcal{M} \,=\, \mathcal{H}_{\scriptscriptstyle n}$ & Hessian metric (\ref{eq:hnmetric}) &   \\
& & & \\[0.15cm]
\multirow{2}{*}{\parbox{26mm}{Paragraph \ref{subsec:tn}\,: $n \times n$\\Toeplitz covariance \\matrices }} &  &  & &
\multirow{2}{*}{\parbox{26mm}{$\phantom{mm}$Algorithm 3.\\$\phantom{m}$ in Appendix \ref{app:sample}}} \\
& $\mathcal{M} \,=\, \mathcal{T}_{\scriptscriptstyle n}$ & Hessian metric (\ref{eq:tnmetric}) & Analytic expression (\ref{eq:tnz0})  \\
& & & \\[0.15cm]
\multirow{2}{*}{\parbox{26mm}{Paragraph \ref{subsec:btn}\,: $n \times n$ block-Toeplitz \\ covariance matrices  \\ with $N\times N$ blocks}} &  &  & \multirow{2}{*}{\parbox{26mm}{$\phantom{mm}$Monte Carlo\\ $\phantom{xx}$integration from \\ $\phantom{xx}\,$(\ref{eq:zdn}) and (\ref{eq:btnz0})}} &
\multirow{2}{*}{\parbox{26mm}{$\phantom{mm}$Algorithm 4.\\$\phantom{m}$ in Appendix \ref{app:sample}}} \\
& $\mathcal{M} \,=\, \mathcal{T}^{\scriptscriptstyle N}_{\scriptscriptstyle n}$ & Hessian metric (\ref{eq:btnmetric}) &  \\
& & &
\end{tabular}
\end{center}
\section{Mathematical Background on Riemannian symmetric spaces} \label{sec:geo}
This section gives five facts about Riemannian symmetric spaces, Facts \ref{fact:1}--\ref{fact:5} below, which will be essential for the following.
These are here stated in their general form, but three concrete examples are presented in Section \ref{sec:matrices}.

A Riemannian symmetric space is a Riemannian manifold $\mathcal{M}$ such that, for each point $x \in \mathcal{M}$, there exists an isometry transformation $I_x\,: \mathcal{M} \rightarrow \mathcal{M}$ whose effect is to reverse geodesic curves passing through $x$~\cite{helgason} (Chapter IV, Page $170$). Any Riemannian symmetric space is a Riemannian homogeneous space, although the converse is not true. In particular, Riemannian symmetric spaces share in two important properties of Riemannian homogeneous spaces, given in Facts \ref{fact:1} and \ref{fact:2} below\,: invariance of distance, and invariance of integrals. 

To say that $\mathcal{M}$ is a Riemannian homogeneous space means a Lie group $G$ acts transitively and isometrically on $\mathcal{M}$~\cite{helgason} (Chapter II, Page $113$). Each element $g \in G$ defines a geometric transformation of $\mathcal{M}$ which maps each point $x \in \mathcal{M}$ to its image $g\cdot x$. These transformations satisfy the group action property,
\begin{equation} \label{eq:action}
  (\,g_{\scriptscriptstyle 1}g_{\scriptscriptstyle 2}) \cdot x \, =  g_{\scriptscriptstyle 1}\cdot (\, g_{\scriptscriptstyle 2}\cdot x)
\end{equation}
Then, transitivity of group action means that for any two points $x$ and $y$ of $\mathcal{M}$, there exists $g \in G$ such that $g\cdot x = y$. 

When $G$ acts on $\mathcal{M}$ transitively, one says that $\mathcal{M}$ is a homogeneous space. Homogeneous spaces admit a general description, as the quotient of $G$ by a compact Lie subgroup $H$. Precisely, choose some point $o \in \mathcal{M}$ as origin, and let $H$ be the subgroup of all $h \in G$ such that $h\cdot o = o$~\cite{helgason} (Chapter II, Page $113$). Then, $\mathcal{M}$ is identified with the quotient $G/H$.

When $G$ acts on $\mathcal{M}$ isometrically, the transformations $x \mapsto g\cdot x$ preserve the Riemannian metric and Riemannian distance of $\mathcal{M}$. This is expressed in Fact \ref{fact:1} below, which uses the following notation~\cite{berger}. The Riemannian metric of $\mathcal{M}$ is denoted $ds^{\scriptscriptstyle 2}$. For $x \in \mathcal{M}$, this defines a quadratic form\,: $ds^{\scriptscriptstyle 2}_{\scriptscriptstyle x}(dx) = \mbox{squared length of a displacement } dx \mbox{ away from }x$. The Riemannian distance corresponding to the metric $ds^{\scriptscriptstyle 2}$ is denoted $d(x,y)$ for $x,y \in \mathcal{M}$. Fact \ref{fact:1} is not a theorem, but really part of the definition of a Riemannian homogeneous space~\cite{helgason}.
\vspace{0.1cm}
\begin{fact}[Invariant distance] \label{fact:1}
 For $g \in G$ and $x,y \in \mathcal{M}$,
\begin{subequations} \label{eq:invmetric}
\begin{equation} \label{eq:invmetric1}
   \mbox{preservation of metric\,: }  \hspace{1.5cm} z = g\cdot x \, \Rightarrow \, ds^{\scriptscriptstyle 2}_{\scriptscriptstyle z}(dz) = ds^{\scriptscriptstyle 2}_{\scriptscriptstyle x}(dx) \hspace{2.1cm}
\end{equation}
\begin{equation} \label{eq:invmetric2}
\hspace{-0.35cm} \mbox{preservation of distance\,:} \hspace{0.8cm} \hspace{0.5cm} \;d(\,g\cdot x,g\cdot y) = d(x,y)  \hspace{3cm}
\end{equation}
\end{subequations}
\end{fact}

\indent Turning to Fact \ref{fact:2}, recall that the Riemannian metric $ds^{\scriptscriptstyle 2}$ defines a Riemannian volume element $dv$~\cite{berger}. Fact \ref{fact:1} states that the transformations $x \mapsto g\cdot x$ preserve the metric $ds^{\scriptscriptstyle 2}$. Therefore, these transformations must also preserve the volume element $dv$. This leads to the invariance of integrals with respect to $dv$, as follows~\cite{helgason} (Chapter X, Page $361$). 
\begin{fact}[Invariant integrals] \label{fact:2}
 For any integrable function $f : \mathcal{M} \rightarrow \mathbb{R}$ and $g \in G$,
 \begin{equation} \label{eq:invintegral}
   \mbox{preservation of volume\,: } \hspace{0.5cm} \int_{\scriptscriptstyle \mathcal{M}}\, f(g\cdot x) \, dv(x) = \int_{\scriptscriptstyle \mathcal{M}}\, f(x) \, dv(x) \hspace{0.7cm}
 \end{equation}
\end{fact}
\indent Facts \ref{fact:1} and \ref{fact:2} above hold for any Riemannian homogeneous space. The following fact \ref{fact:3} requires that $\mathcal{M}$ should have non-positive curvature. This guarantees the existence and uniqueness of Riemannian barycentres in $\mathcal{M}$~\cite{afsari}.

Let $\pi$ be a probability distribution on $\mathcal{M}$. The Riemannian barycentre of $\pi$ is a point $\hat{x}_{\scriptscriptstyle \pi}$ defined as follows. Consider $E_{\scriptscriptstyle \pi} : \mathcal{M} \rightarrow \mathbb{R}_+\,$,
\begin{equation} \label{eq:variance}
 \mbox{variance function\,:} \hspace{0.5cm} E_{\scriptscriptstyle \pi}(x) = \int_{\scriptscriptstyle \mathcal{M}}\, d^{\scriptscriptstyle\, 2}(x,z) \, d\pi(z)
\end{equation}
The Riemannian barycentre $\hat{x}_{\scriptscriptstyle \pi}$ is then a global minimiser of $E_{\scriptscriptstyle \pi}$. Fact \ref{fact:3} states that $\hat{x}_{\scriptscriptstyle \pi}$ exists and is unique~\cite{afsari} (Theorem 2.1., Page $659$).
\begin{fact}[Riemannian barycentre] \label{fact:3}
  Assume that $\mathcal{M}$ has non-positive curvature. For any probability distribution $\pi$ on $\mathcal{M}$, the Riemannian barycentre $\hat{x}_{\scriptscriptstyle \pi}$ exists and is unique. Moreover, $\hat{x}_{\scriptscriptstyle \pi}$ is the unique stationary point of the variance function $E_{\scriptscriptstyle \pi}$.
\end{fact}
\indent Fact \ref{fact:2} above states that integrals with respect to the Riemannian volume element $dv$ are invariant. However, it does not show how these integrals can be computed. This is done by the following Fact \ref{fact:4}, which requires that $\mathcal{M}$ should be a Riemannian symmetric space of non-compact type. This means that the group $G$ is a semisimple Lie group of non-compact type~\cite{helgason} (Chapter V). It should be noted that 
this property implies that $\mathcal{M}$ has (strictly) negative curvature, so that Fact \ref{fact:3} holds true.

To state Fact \ref{fact:4}, assume that $G$ is a semisimple Lie group of non-compact type. Recall that $H$ is a compact Lie subgroup of $G$. Let $\mathfrak{g}$ and $\mathfrak{h}$ denote the Lie algebras of $G$ and $H$ respectively, and let $\mathfrak{g}$ have the following Iwasawa decomposition~\cite{helgason} (see details in Chapter VI, Page $219$),
\begin{equation} \label{eq:iwasawa}
  \mbox{Iwasawa decomposition\,:} \hspace{0.5cm} \mathfrak{g} \,=\, \mathfrak{h} \, + \, \mathfrak{a} \, + \, \mathfrak{n}
\end{equation}
where $\mathfrak{a}$ is an Abelian subalgebra of of $\mathfrak{g}$, and $\mathfrak{n}$ is a nilpotent subalgebra of $\mathfrak{g}$. Each $x \in \mathcal{M}$ can be written in the following form\!~\cite{helgason} (Chapter X),
\begin{equation} \label{eq:polar}
 \mbox{polar coordinates\,:} \hspace{0.5cm} x = \exp\left(\mathrm{Ad}(h)\, a\right) \cdot o \hspace{1cm}  a \in \mathfrak{a}\, , \, h \in H
\end{equation}
where $\exp : \mathfrak{g} \rightarrow G$ is the Lie group exponential, and $\mathrm{Ad}$ is the adjoint representation. Here, it will be said that $(a,h)$ are polar coordinates of $x$, which is then written $x = x(a,h)$ (this is an abuse of notation, since $(a,h)$ are not unique, but it will provide useful intuition).

The following integral formula holds~\cite{helgason} (Chapter X, Page $382$). It makes use of the concept of roots. These are linear mappings $\lambda : \mathfrak{a} \rightarrow \mathbb{R}$, which characterise the Lie algebra $\mathfrak{g}$. For a minimal, introduction to this concept, see Appendix \ref{app:roots}.
\begin{fact}[Integration in polar coordinates] \label{fact:4}
 For any integrable function $f : \mathcal{M} \rightarrow \mathbb{R}$,
 \begin{equation} \label{eq:integralp}
     \int_{\scriptscriptstyle \mathcal{M}}\, f( x) \, dv(x) = \mathrm{C} \times \int_{\scriptscriptstyle H} \int_{\scriptscriptstyle \mathfrak{a}} \, f(a,h) \, D(a) \, da \, dh
 \end{equation}
 where $\mathrm{C}$ is a constant which does not depend on the function $f$, $dh$ is the normalised Haar measure on $H$ and $da$ is the Lebesgue measure on $\mathfrak{a}$. The function $D : \mathfrak{a} \rightarrow \mathbb{R}_+$ is given by 
 \begin{equation} \label{eq:D}
  D(a) = \prod_{\scriptscriptstyle \lambda > 0} \, \sinh^{\scriptscriptstyle m_{\lambda}} (|\lambda(a)|)
 \end{equation}
 where the product is over positive roots $\lambda : \mathfrak{a} \rightarrow \mathbb{R}$, and $m_{\scriptscriptstyle \lambda}$ is the dimension of the root space corresponding to $\lambda$. 
 \end{fact}
\indent Fact \ref{fact:5} provides a geometric interpretation of the polar coordinates $(a,h)$  introduced in (\ref{eq:polar}). Just like usual polar coordinates, these give a direct description of geodesics and distance, with respect to the origin $o$~\cite{helgason} (Chapter IV). 
\begin{fact}[Geometric interpretation of polar coordinates] \label{fact:5} The two following properties hold true. \\[0.1cm]
(i) Geodesic curves passing through the origin $o \in \mathcal{M}$ are exactly the curves $x : \mathbb{R} \rightarrow \mathcal{M}$ given by $x(t) = x(t \,a, h)$ where $a$ and $h$ are constant. \\[0.1cm]
(ii) The Riemannian distance between the origin $o \in \mathcal{M}$ and $x(a,h)$ is given by
\begin{equation} \label{eq:killing}
\mbox{distance to origin\,:}\hspace{0.5cm}  d^{\scriptscriptstyle \,2}(o,x) \, = \, B(a,a)
\end{equation}
where $B : \mathfrak{g} \times \mathfrak{g} \rightarrow \mathbb{R}$ is the Killing form, an $\mathrm{Ad}$-invariant symmetric bilinear form on $\mathfrak{g}$, which is positive definite on $\mathfrak{a}$~\cite{helgason} (Chapter II). 
\end{fact}
\section{Gaussian distributions on Riemannian symmetric spaces} \label{sec:gaussian}
This section introduces Gaussian distributions on Riemannian symmetric spaces, and studies maximum likelihood estimation of their parameters. These distributions establish a connection between the statistical concept of maximum likelihood estimation, and the geometric concept of Riemannian barycentre.

Here, $\mathcal{M}$ is a Riemannian symmetric space of non-positive curvature. Recall $d(x,y)$ denotes Riemanniain distance between $x,y \in \mathcal{M}$, and $dv$ denotes the Riemannian volume element of $\mathcal{M}$. A Gaussian distribution is one that has the following probability
density with respect to $dv$,
\begin{equation} \label{eq:gaussianpdf}
  p(x|\,\bar{x},\sigma) =  \frac{1}{Z(\sigma)} \times \exp \left[-\,\frac{d^{\,\scriptscriptstyle 2}(x,\bar{x})}{2\sigma^{\,\scriptscriptstyle 2}}\right]
\end{equation}
where $\bar{x} \in \mathcal{M}$ and $\sigma > 0$. This distribution will be denoted $G(\bar{x},\sigma)$. 

The main feature of the probability density (\ref{eq:gaussianpdf}) is that the normalising factor $Z(\sigma)$ does not depend on the parameter $\bar{x}$. This will be proved as part of Proposition \ref{prop:z} below, and is the key ingredient in the relationship between Gaussian distributions and the concept of Riemannian barycentre. In~\cite{vemuri}\cite{said}, this relationship was developed in the special case where $\mathcal{M} = \mathcal{P}_{\scriptscriptstyle n\,}$, the symmetric space of real $n \times n$ covariance matrices.  The present paper considers the general case where $\mathcal{M}$ is any symmetric space of non-positive curvature. This allows for Gaussian distributions of the form (\ref{eq:gaussianpdf}) to be extended to spaces of structured covariance matrices, such as complex, Toeplitz, and block-Toeplitz covariance matrices, as done in Section \ref{sec:matrices} below. 

In the following, Paragraph \ref{subsec:gaussiandef} gives Proposition \ref{prop:z} which characterises the normalising factor $Z(\sigma)$, and Proposition \ref{prop:sample} which explains how to sample from a Gaussian distribution with given parameters. Paragraph \ref{subsec:inference} is concerned with maximum likelihood estimation of the parameters of a Gaussian distribution. In particular, it contains Proposition \ref{prop:mle} which shows that the maximum likelihood estimate $\hat{x}_{\scriptscriptstyle N}$ of the parameter $\bar{x}$ is precisely the Riemannian barycentre of available observations $x_{\scriptscriptstyle 1},\ldots, x_{\scriptscriptstyle N}$. Finally, Proposition \ref{prop:asymp} can be used to obtain the consistency and asymptotic normality of $\hat{x}_{\scriptscriptstyle N}$.

\subsection{Definition\,: the role of invariance} \label{subsec:gaussiandef}
This paragraph shows how the invariance properties of the Riemannian symmetric space $\mathcal{M}$, described in Section \ref{sec:geo}, can be used to characterise Gaussian distributions on $\mathcal{M}$. Precisely, for a Gaussian distribution $G(\bar{x},\sigma)$, it is possible to obtain a) an exact expression of the normalising factor $Z(\sigma)$ appearing in (\ref{eq:gaussianpdf}), and b) a generally applicable method for sampling from this distribution $G(\bar{x},\sigma)$. 

Proposition \ref{prop:z} is concerned with the normalising factor $Z(\sigma)$. A fully general expression of $Z(\sigma)$ is given by (\ref{eq:zint}) and (\ref{eq:zprod}) in the proof of this proposition. To state Proposition \ref{prop:z}, consider the following notation,
\begin{equation} \label{eq:profile}
f(x|\,\bar{x},\sigma) \, = \, \exp \left[-\,\frac{d^{\,\scriptscriptstyle 2}(x,\bar{x})}{2\sigma^{\,\scriptscriptstyle 2}}\right]
\hspace{1cm} Z(\bar{x},\sigma) \, = \, \int_{\scriptscriptstyle \mathcal{M}}\, f(x|\,\bar{x},\sigma)\, dv(x)
\end{equation}
In other words, $f(x|\,\bar{x},\sigma)$ is the ``profile'' of the probability density $p(x|\,\bar{x},\sigma)$.
 \vspace{0.1cm}
\begin{proposition}[Normalising factor] \label{prop:z}
The two following properties hold true. \\[0.1cm]
(i) $Z(\bar{x},\sigma) = Z(o\,,\sigma)$ for any $\bar{x} \in \mathcal{M}$, where $o$ is some point chosen as origin of $\mathcal{M}$. \\[0.1cm]
(ii)  $Z(\sigma) = Z(o\,,\sigma)$ is a strictly log-convex function of the parameter $\eta = -1/2\sigma^{\scriptscriptstyle \,2}$.
\end{proposition}
 \vspace{0.1cm}
\textbf{Proof\,: } the proof of (i) follows from Facts \ref{fact:1} and \ref{fact:2}. Since $G$ acts transitively on $\mathcal{M}$, any $\bar{x} \in \mathcal{M}$ is of the form $x = g\cdot o$ for some $g \in G$. Applying (\ref{eq:invmetric2}) from Fact \ref{fact:2} to (\ref{eq:profile}), it then follows that $f(x|\,\bar{x},\sigma) \, = \, f(g^{\scriptscriptstyle -1}\cdot x|\,o\,,\sigma)$. Accordingly,
$$
Z(\bar{x},\sigma) \, = \, \int_{\scriptscriptstyle \mathcal{M}}\, f(x|\,\bar{x},\sigma)\, dv(x) \, = \, 
\int_{\scriptscriptstyle \mathcal{M}}\, f(g^{\scriptscriptstyle -1}\cdot x|\,o\,,\sigma) \, dv(x) 
$$
Now, by Fact \ref{fact:2}, (and using $g^{\scriptscriptstyle -1}$ instead of $g$ in (\ref{eq:invintegral})),
$$
\int_{\scriptscriptstyle \mathcal{M}}\, f(g^{\scriptscriptstyle -1}\cdot x|\,o\,,\sigma) \, dv(x) \,= \,
\int_{\scriptscriptstyle \mathcal{M}}\, f( x|\,o\,,\sigma) \, dv(x) \, = \, Z(o\,,\sigma)
$$
so that (i) is indeed true. The proof of (ii) is based on Facts \ref{fact:4} and \ref{fact:5}. Assume first that $\mathcal{M}$ is a Riemannian symmetric space of non-compact type. By (\ref{eq:integralp}) and (\ref{eq:killing}),
\begin{equation} \label{eq:zgeneral}
Z(\sigma) \,=\, \int_{\scriptscriptstyle \mathcal{M}}\, f( x|\,o\,,\sigma) \, dv(x) \, = \,\mathrm{C} \times \int_{\scriptscriptstyle H} \int_{\scriptscriptstyle \mathfrak{a}} \,  
\exp \left[-\,\frac{B(a,a)}{2\sigma^{\,\scriptscriptstyle 2}}\right]\, D(a) \, da \, dh
\end{equation}
But the function under the integral does not depend on $h$. Since $dh$ is the normalised Haar measure on $H$,
\begin{equation} \label{eq:zint}
Z(\sigma) \,=\, \mathrm{C} \times \int_{\scriptscriptstyle \mathfrak{a}} \,  
\exp \left[-\,\frac{B(a,a)}{2\sigma^{\,\scriptscriptstyle 2}}\right]\, D(a) \, da
\end{equation}
In general, a Riemannian symmetric space $\mathcal{M}$ of non-positive curvature decomposes into a product Riemannian symmetric space~\cite{helgason} (Chapter V, Page $208$),
\begin{equation} \label{eq:ssdecomposition}
\mathcal{M} =  \mathcal{M}_{\scriptscriptstyle 1} \times \ldots \times \mathcal{M}_{\scriptscriptstyle r}
\end{equation}
where each $\mathcal{M}_{\scriptscriptstyle p}$ is either a Euclidean space or a Riemannian symmetric space of non-compact type. Now, each $x \in \mathcal{M}$ is an $r$-tuple $x = (x_{\scriptscriptstyle 1},\ldots,x_{\scriptscriptstyle r})$ where $x_{\scriptscriptstyle p} \in \mathcal{M}_{\scriptscriptstyle p}$. Moreover, for $x,y \in \mathcal{M}$, 
\begin{equation} \label{eq:proddistance}
  d^{\scriptscriptstyle\, 2}(x,y) = \sum^{\scriptscriptstyle r}_{\scriptscriptstyle p=1} \, d^{\scriptscriptstyle\, 2}_{\scriptscriptstyle p}(x_{\scriptscriptstyle p}\,,\,y_{\scriptscriptstyle p})
\end{equation}
where $d_{\scriptscriptstyle p}(x_{\scriptscriptstyle p}\,,\,y_{\scriptscriptstyle p})$ is the Riemannian distance in $\mathcal{M}_{\scriptscriptstyle p}$ between
$x_{\scriptscriptstyle p}$ and $y_{\scriptscriptstyle p}$ (where $y = (y_{\scriptscriptstyle 1},\ldots,y_{\scriptscriptstyle r})$). Replacing (\ref{eq:proddistance}) in (\ref{eq:profile}), it follows that,
\begin{equation} \label{eq:zprod}
  Z(\sigma) \, = \,   Z_{\scriptscriptstyle 1}(\sigma) \times \ldots \times Z_{\scriptscriptstyle r}(\sigma)
\end{equation}
where $Z_{\scriptscriptstyle p}(\sigma)$ is the normalising factor for the space $\mathcal{M}_{\scriptscriptstyle p}$. Recall that a product of strictly log-convex functions is strictly log-convex. If $\mathcal{M}_{\scriptscriptstyle p}$ is a Euclidean space of dimension $q$, then $Z_{\scriptscriptstyle p}(\sigma) = (2\pi\sigma^{\scriptscriptstyle 2})^{\scriptscriptstyle -q/2}$, which is a strictly log-convex function of $\eta$. If $\mathcal{M}_{\scriptscriptstyle p}$ is a Riemannian symmetric space of non-compact type, $Z_{\scriptscriptstyle p}(\sigma)$ is given by (\ref{eq:zint}). Let $\rho = B(a,a)$ and $\mu(d\rho)$ the image of $D(a)da$ under the mapping $a \mapsto \rho$. It follows that
\begin{equation} \label{eq:mgf}
  Z(\sigma) = C \times \int^{\scriptscriptstyle \infty}_{\scriptscriptstyle 0} \, e^{\eta\, \rho}\,\mu(d\rho)
\end{equation}
As a function of $\eta\,$, this is the moment-generating function of the positive measure $\mu(d\rho)$, and therefore a strictly log-convex function~\cite{bndrf}, (Chapter 7, Page $103$). Thus, $Z(\sigma)$ is a product (\ref{eq:zprod}) of strictly log-convex functions of $\eta$, and therefore itself strictly log-convex. This proves (ii). \hfill$\blacksquare$

The following Proposition \ref{prop:sample} shows that polar coordinates (\ref{eq:polar}) can be used to effectively sample from a Gaussian distribution $G(\bar{x},\sigma)$. This proposition is concerned with the special case where $\mathcal{M}$ is a Riemannian symmetric space of non-compact type. However, the general case where $\mathcal{M}$ is a Riemannian symmetric space of non-positive curvature can then be obtained using decomposition (\ref{eq:ssdecomposition}). This is explained in the remark following the proposition.

In the statement of Proposition \ref{prop:sample}, the following usual notation is used. If $X$ is a random variable and $P$ is a probability distibution, $X\sim P$ means that $X$ follows the probability distribution $P$.
\vspace{0.1cm}
\begin{proposition}[Gaussian distribution via polar coordinates] \label{prop:sample}
Let $\mathcal{M}$ be a Riemannian symmetric space of non-compact type. Let $\mathcal{M} = G/H$ where $\mathfrak{g}$ has Iwasawa decomposition (\ref{eq:iwasawa}). Let $h$ and $a$ be independent random variables in $H$ and $\mathfrak{a}$, respectively. Assume $h$ is uniformly distributed on $H$, and $a$ has probability density, with respect to the Lebesgue measure $da$ on $\mathfrak{a}$,
\begin{equation} \label{eq:apdf}
   p(a) \, \propto \, \exp \left[-\,\frac{B(a,a)}{2\sigma^{\,\scriptscriptstyle 2}}\right]\, D(a) 
\end{equation}
where $\propto$ indicates proportionality. The following hold true.\\[0.1cm]
(i) If $x = x(a,h)$ as in (\ref{eq:polar}), then $x \sim G(o,\sigma)$. \\[0.1cm]
(ii) If $x \sim G(o,\sigma)$ and $\bar{x} = g\cdot o\,$, then $g\cdot x \sim G(\bar{x},\sigma)$. 
\end{proposition}
 \vspace{0.1cm}
\textbf{Proof\,: } the proof of (ii) is a direct application of Fact \ref{fact:2}, (compare to~\cite{said}, proof of Proposition 5). Only the proof of (i) will be detailed. This follows using Fact \ref{fact:4}.  Recall that a uniformly distributed random variable $h$ in $H$ has constant probability density, identically equal to $1$, with respect to the normalised Haar measure $dh$. Let $h$ and $a$ be as in the statement of the proposition. Since $h$ and $a$ are independent, (\ref{eq:apdf}) implies that for any function $\varphi : H \times \mathfrak{a} \rightarrow \mathbb{R}$, the expectation of $\varphi(h,a)$ is equal to
$$
\int_{\scriptscriptstyle H} \int_{\scriptscriptstyle \mathfrak{a}} \,  \varphi(h,a) \, p(a) \, da \, dh \, \propto \,
\int_{\scriptscriptstyle H} \int_{\scriptscriptstyle \mathfrak{a}} \,  \varphi(h,a) \, \exp \left[-\,\frac{B(a,a)}{2\sigma^{\,\scriptscriptstyle 2}}\right]\, D(a) \, da \, dh 
$$
Then, if $x = x(a,h)$ as in (\ref{eq:polar}) and $\tilde{\varphi}: \mathcal{M} \rightarrow \mathbb{R}$, the expectation of $\tilde{\varphi}(x)$ is proportional to,
$$
\int_{\scriptscriptstyle H} \int_{\scriptscriptstyle \mathfrak{a}} \,  \tilde{\varphi}(x(h,a)) \, \exp \left[-\,\frac{B(a,a)}{2\sigma^{\,\scriptscriptstyle 2}}\right]\, D(a) \, da \, dh \, \propto \, \int_{\scriptscriptstyle \mathcal{M}} \, \tilde{\varphi}(x) \, p(x|\,o,\sigma) \, dv(x) 
$$
as follows from (\ref{eq:integralp}), (\ref{eq:killing}) and (\ref{eq:gaussianpdf}). Since the function $\tilde{\varphi}$ is arbitrary, this shows that the probability density of $x$ with respect to $dv$ is equal to $p(x|\,o,\sigma)$. In other words, $x \sim G(o,\sigma)$. \hfill$\blacksquare$ \\[0.15cm]
\indent $\blacktriangleright$ \textbf{Remark 1} : Proposition \ref{prop:sample} only deals with the special case where $\mathcal{M}$ is a Riemannian symmetric space of non-compact type. The general case where $\mathcal{M}$ is a Riemannian symmetric space of non-positive curvature is obtained as follows. Let $\mathcal{M}$ have decomposition (\ref{eq:ssdecomposition}) and $x = (x_{\scriptscriptstyle 1},\ldots,x_{\scriptscriptstyle r})$ be a random variable in $\mathcal{M}$ with distribution $G(\bar{x},\sigma)$. If $\bar{x} = (\bar{x}_{\scriptscriptstyle 1},\ldots,\bar{x}_{\scriptscriptstyle r})$, then it follows from  (\ref{eq:gaussianpdf}) and (\ref{eq:proddistance}) that
\begin{equation} \label{eq:prodpdf}
  p(x|\,\bar{x},\sigma) = p(x_{\scriptscriptstyle 1}|\,\bar{x}_{\scriptscriptstyle 1},\sigma) \, \times \, \ldots \, \times \,
 p(x_{\scriptscriptstyle r}|\,\bar{x}_{\scriptscriptstyle r},\sigma)
\end{equation}
Therefore, $x_{\scriptscriptstyle 1},\ldots, x_{\scriptscriptstyle r}$ are independent and each $x_{\scriptscriptstyle p}$ has distribution $G(\bar{x}_{\scriptscriptstyle p\,},\sigma)$. Conversely, (\ref{eq:prodpdf}) implies that if $x_{\scriptscriptstyle 1},\ldots, x_{\scriptscriptstyle r}$ are independent and each $x_{\scriptscriptstyle p}$ has distribution $G(\bar{x}_{\scriptscriptstyle p\,},\sigma)$, then $x$ has distribution $G(\bar{x},\sigma)$. \hfill $\blacksquare$ \\[0.1cm]
\indent Proposition \ref{prop:sample}, along with the subsequent remark, provides a generally applicable method for sampling from a Gaussian distribution $G(\bar{x},\sigma)$ on a Riemanniain symmetric space of non-positive curvature $\mathcal{M}$. This proceeds as follows, using decomposition (\ref{eq:ssdecomposition}). 
\vspace{0.3cm}
\hrule
\vspace{0.3cm}

\noindent \texttt{l1}.\; \texttt{for} $p = 1,\ldots, r$ \texttt{do}\\[0.1cm]
\texttt{l2}.\; \; \; \texttt{if} $\mathcal{M}_{\scriptscriptstyle p}$ is a Euclidean space\\[0.1cm]
\texttt{l3}.\; \; \; \; \; sample $x_{\scriptscriptstyle p}$ from a multivariate normal distribution $\mathcal{N}(\bar{x}_{\scriptscriptstyle p},\sigma)$ \\[0.1cm] 
\texttt{l4}.\;  \; \; \texttt{end if} \\[0.1cm]
\texttt{l5}.\; \; \; \texttt{if} $\mathcal{M}_{\scriptscriptstyle p} = G/H$ is a Riemannian symmetric space of non-compact type \\[0.1cm]
\texttt{l6}.\; \; \; \; \; sample $h$ from a uniform distribution on $H$ \\[0.1cm]
\texttt{l7}.\; \; \; \; \; sample $a$ from the multivariate density (\ref{eq:apdf}) \\[0.1cm]
\texttt{l8}.\; \; \; \; \; put $x_{\scriptscriptstyle p} = x(a,h)$ as in (\ref{eq:polar})  \\[0.1cm]
\texttt{l9}.\; \; \; \; \; reset $x_{\scriptscriptstyle p}$ to  $g\cdot x_{\scriptscriptstyle p}$ where $\bar{x} = g\cdot o$ as in (ii) of Proposition \ref{prop:sample}.  \\[0.1cm]
\texttt{l10}.\; \; \;  \texttt{end if} \hspace{8cm} {\small \sc Algorithm 1.\,: sampling from a Gaussian}\\[0.1cm]
\texttt{l11}.\! \texttt{end for} \hspace{8.5cm}\,\,\,\,\,\; {\small \sc distribution on a symmetric space $\mathcal{M}$} \\[0.1cm]
\texttt{l12}.\! put $x = (x_{\scriptscriptstyle 1},\ldots,x_{\scriptscriptstyle r})$  
\vspace{0.2cm}
\hrule
\vspace{0.3cm}

\indent  The above steps are completely general, and their concrete implementation depends on the specific space $\mathcal{M}$ at hand. When $\mathcal{M} = \mathcal{P}_{\scriptscriptstyle n\,}$, the symmetric space of real $n \times n$ covariance matrices, this was described in~\cite{said} (Page $8$). Three new implementations, for spaces of structured covariance matrices, will be given in Appendix \ref{app:sample}.
\subsection{Maximum likelihood estimation} \label{subsec:inference}
The current paragraph studies maximum likelihood estimation of the parameters $\bar{x}$ and $\sigma$ of a Gaussian distribution $G(\bar{x},\sigma)$, on a Riemannian symmetric space of non-positive curvature $\mathcal{M}$. Proposition \ref{prop:mle} characterises the resulting maximum likelihood estimates, and Proposition \ref{prop:asymp} yields their asymptotic behavior. Let $\hat{x}_{\scriptscriptstyle N}$ and $\hat{\sigma}_{\scriptscriptstyle N}$ be maximum likelihood estimates, based on independent samples $x_{\scriptscriptstyle 1},\ldots, x_{\scriptscriptstyle N}$ from $G(\bar{x},\sigma)$. Proposition \ref{prop:mle} shows that $\hat{x}_{\scriptscriptstyle N}$ and $\hat{\sigma}_{\scriptscriptstyle N}$ provide a geometric description of these samples. Indeed, the proposition states that $\hat{x}_{\scriptscriptstyle N}$ is the Riemannian barycentre of $x_{\scriptscriptstyle 1},\ldots, x_{\scriptscriptstyle N}$, and that $\hat{\sigma}_{\scriptscriptstyle N}$ measures the dispersion (that is, mean squared distance) of $x_{\scriptscriptstyle 1},\ldots, x_{\scriptscriptstyle N}$ away from this Riemannain barycentre. 

For Proposition \ref{prop:mle}, recall the Riemannian barycentre of $x_{\scriptscriptstyle 1}, \ldots, x_{\scriptscriptstyle N}$ is a global minimiser of the function $E_{\scriptscriptstyle N} : \mathcal{M} \rightarrow \mathbb{R}$, defined as follows\!~\cite{afsari},
\begin{equation} \label{eq:empiricalvar}
\mbox{empirical variance function\,:} \hspace{0.5cm} E_{\scriptscriptstyle N}(x) = \, \frac{1}{N} \, \sum^{\scriptscriptstyle N}_{\scriptscriptstyle n=1} 
d^{\scriptscriptstyle \,2}(x,x_{\scriptscriptstyle n})
\end{equation}
This is called an empirical variance function, since it is a special case of the variance function (\ref{eq:variance}), corresponding to the empirical distribution $\pi = N^{\scriptscriptstyle -1}\,{\scriptstyle \sum^{\scriptscriptstyle N}_{\scriptscriptstyle n=1} \delta_{\scriptscriptstyle x_n}\,}$, where $\delta_{\scriptscriptstyle x}$ denotes the Dirac distribution concentrated at $x \in \mathcal{M}$. Accordingly, existence and uniqueness of the Riemannian barycentre of $x_{\scriptscriptstyle 1}, \ldots, x_{\scriptscriptstyle N}$ follow from Fact \ref{fact:3}. 
\begin{proposition}[MLE and Riemannian barycentre] \label{prop:mle}
 Let $x_{\scriptscriptstyle 1},\ldots, x_{\scriptscriptstyle N}$ be independent samples from $G(\bar{x},\sigma)$.  The maximum likelihood estimates $\hat{x}_{\scriptscriptstyle N}$ and $\hat{\sigma}_{\scriptscriptstyle N}$, based on these samples, exist and are unique. Precisely, \\[0.1cm]
(i) $\hat{x}_{\scriptscriptstyle N}$ is the Riemannian barycentre of $x_{\scriptscriptstyle 1},\ldots, x_{\scriptscriptstyle N}$. \\[0.1cm]
(ii) There exists a strictly increasing bijective function $\Phi : \mathbb{R}_{+} \rightarrow \mathbb{R}_+$, such that
\begin{equation} \label{eq:phi}
 \hat{\sigma}_{\scriptscriptstyle N} \, = \, \Phi \left( {\scriptstyle \frac{1}{\mathstrut N}} \, {\scriptstyle \sum^{\scriptscriptstyle N}_{\scriptscriptstyle n=1}} 
d^{\scriptscriptstyle \,2}(\hat{x}_{\scriptscriptstyle N},x_{\scriptscriptstyle n}) \right)
\end{equation}
\end{proposition}
 \vspace{0.1cm}
\textbf{Proof\,: } from (\ref{eq:gaussianpdf}), the log-likelihood function for the parameters $\bar{x}$ and $\sigma$ is given by
$$
L(\bar{x},\sigma) \, = \, -N\, \log Z(\sigma) \, - \, \frac{1}{\mathstrut 2\sigma^{\scriptscriptstyle 2}} \, \sum^{\scriptscriptstyle N}_{\scriptscriptstyle n=1}
d^{\scriptscriptstyle \,2}(\bar{x},x_{\scriptscriptstyle n})
$$
By definition, $\hat{x}_{\scriptscriptstyle N}$ and $\hat{\sigma}_{\scriptscriptstyle N}$ are found by maximising this function. The sum in the second term on the right hand side does not depend on $\sigma$. Therefore, $\hat{x}_{\scriptscriptstyle N}$ can be found separately, by minimising this sum over the values of $\bar{x}$ (minimisation instead of maximisation is due to the minus sign ahead of the sum). However, this is the same as minimising the empirical variance function (\ref{eq:empiricalvar}). Now, the unique global minimiser of the empirical variance function is the Riemannian barycentre of $x_{\scriptscriptstyle 1},\ldots, x_{\scriptscriptstyle N}$, whose existence and uniqueness follow from Fact \ref{fact:3} as noted before the proposition. This proves (i). Now, assume $\hat{x}_{\scriptscriptstyle N}$ has been found. Then, to find $\hat{\sigma}_{\scriptscriptstyle N}$, it is convenient to maximise $L(\hat{x}_{\scriptscriptstyle N}\,,\sigma)$ over the values of the parameter $\eta = -1/2\sigma^{\scriptscriptstyle \,2\,}$, which was already used in Proposition \ref{prop:z}. Indeed, note that
$$
\frac{1}{N} \, L(\hat{x}_{\scriptscriptstyle N}\,,\sigma) \, = \, \eta\, \hat{\rho} - \psi(\eta)
$$
where $\hat{\rho} = N^{\scriptscriptstyle -1}\,{\scriptstyle \sum^{\scriptscriptstyle N}_{\scriptscriptstyle n=1}}\, d^{\scriptscriptstyle \,2}(\hat{x}_{\scriptscriptstyle N},x_{\scriptscriptstyle n})$ and $\psi(\eta) = \log Z(\sigma)$. Thus, the maximum likelihood estimate $\hat{\eta}_{\scriptscriptstyle N}$ is given by
\begin{equation} \label{eq:legendre1}
   \hat{\eta}_{\scriptscriptstyle N} \, = \, \mathrm{argmax}_{\scriptscriptstyle \eta} \, \left \lbrace \eta\, \hat{\rho} - \psi(\eta) \right \rbrace
\end{equation}
Recall that (ii) of Proposition \ref{prop:z} states that $\psi$ is a strictly convex function, so its derivative $\psi^{\scriptscriptstyle \prime}$ is strictly increasing. But, from (\ref{eq:legendre1}), $\hat{\eta}_{\scriptscriptstyle N}$ is given by $\psi^{\scriptscriptstyle \prime}(\hat{\eta}_{\scriptscriptstyle N}) = \hat{\rho}$. It follows that $\hat{\eta}_{\scriptscriptstyle N}$ exists and is unique, whenever $\hat{\rho}$ is in the range of $\psi^{\scriptscriptstyle \prime}$. It is possible to show, using (\ref{eq:mgf}), that this range is equal to $\mathbb{R}_+$, which implies existence and uniqueness of $\hat{\eta}_{\scriptscriptstyle N}$ in general. Now, to complete the proof of (ii), it is enough to perform the change of parameter from $\hat{\eta}_{\scriptscriptstyle N}$ back to $\hat{\sigma}_{\scriptscriptstyle N}$.\hfill $\blacksquare$ \\[0.15cm]
\indent Proposition \ref{prop:mle} shows how the maximum likelihood estimates $\hat{x}_{\scriptscriptstyle N}$ and $\hat{\sigma}_{\scriptscriptstyle N}$ can be computed, based on the samples $x_{\scriptscriptstyle 1}, \ldots, x_{\scriptscriptstyle N}$. Since $\hat{x}_{\scriptscriptstyle N}$ is the Riemannian barycentre of these samples, $\hat{x}_{\scriptscriptstyle N}$ can be computed using a variety of existing algorithms for computation of Riemannian barycentres. These include gradient descent\!~\cite{pennec1}\cite{lenglet}, Newton method~\cite{groisser}\cite{ferreira}, recursive estimators\!~\cite{sturm}\cite{aistats}, and stochastic methods\!~\cite{bonnabel}\cite{miclo}. Some of these algorithms are specifically adapted to Riemannian symmetric spaces of non-positive curvature (that is, to the context of the present paper)\!~\cite{aistats,sturm,ferreira}. Least computationally expensive are the recursive estimators~\cite{sturm}\cite{aistats}. On the other hand, the most computationally expensive, involving exact calculation of the Hessian of squared Riemannian distance, is the Newton method of~\cite{ferreira}. Computation of $\hat{\sigma}_{\scriptscriptstyle N}$ amounts to solving the equation $\psi^{\scriptscriptstyle \prime}(\hat{\eta}_{\scriptscriptstyle N}) = \hat{\rho}$. This is a scalar non-linear equation, which can be solved using standard methods, such as the Newton-Raphson algorithm.   
\\[0.15cm]
\indent $\blacktriangleright$ \textbf{Remark 2} : the entropy of a Gaussian distribution $G(\bar{x},\sigma)$ is directly related to the function $\psi$ introduced in the proof of Proposition \ref{prop:z}. By definition, the entropy of $G(\bar{x},\sigma)$ is~\cite{cover},
$$
h(\bar{x},\sigma) = \int_{\scriptscriptstyle \mathcal{M}} \, \log p(x|\,\bar{x},\sigma) \times p(x|\,\bar{x},\sigma) \, dv(x)
$$
Using Fact \ref{fact:2} as in the proof of (i) in Proposition \ref{prop:z}, it follows $h(\bar{x},\sigma) = h(o,\sigma)$, so the entropy does not depend on $\bar{x}$. Evaluataion of $h(o,\sigma)$ gives, in the notation of (\ref{eq:mgf}),
$$
h(o,\sigma) = C \times \int^{\scriptscriptstyle \infty}_{\scriptscriptstyle 0} \, \left( {\eta\, \rho}-\psi(\eta)\right) \times e^{{\eta\, \rho}-\psi(\eta)}\,\mu(d\rho)
$$
Following~\cite{bndrf} (Chapter 9), let $\psi^{\scriptscriptstyle *}$ be the Legendre transform of $\psi$ and $\bar{\rho} = \psi^{\scriptscriptstyle \prime}(\eta)$. Then, the entropy $h(o,\sigma)$ is equal to
\begin{equation} \label{eq:entropy}
  \mbox{entropy of a Gaussian distribution\,:} \hspace{0.5cm} \psi^{\scriptscriptstyle *}(\bar{\rho}) \, = \, \eta \, \bar{\rho} - \psi(\eta)
\end{equation}
Essentially, this is due to the fact that, when $\bar{x}$ is fixed, $G(\bar{x},\sigma)$ takes on the form of an exponential distribution with natural parameter $\eta$ and sufficient statistic $\rho$. \hfill $\blacksquare$ \\[0.1cm]
\indent The following Proposition \ref{prop:asymp} clarifies the significance of the parameters $\bar{x}$ and $\sigma$ of a Gaussian distribution $G(\bar{x},\sigma)$. This Proposition is an asymptotic counterpart of Proposition \ref{prop:mle}. In particular, it can be used to obtain the consistency and asymptotic normality of $\hat{x}_{\scriptscriptstyle N}$. This is explained in the remark after the proposition. For the statement of Proposition \ref{prop:asymp}, recall the definition (\ref{eq:variance}) of the variance function $E_{\scriptscriptstyle \pi} : \mathcal{M} \rightarrow \mathbb{R}_+$ of a probability distribution $\pi$ on $\mathcal{M}$. If $\pi = G(\bar{x},\sigma)$, denote $E_{\scriptscriptstyle \pi}(x) = E(x|\,\bar{x},\sigma)$.
\vspace{0.1cm}
\begin{proposition}[significance of the parameters] \label{prop:asymp}
  The parameters $\bar{x}$ and $\sigma$ of a Gaussian distribution $G(\bar{x},\sigma)$ are given by, 
 \begin{subequations} \label{eq:asymp}
   \begin{equation} \label{eq:asymp1}
  \bar{x} \mbox{ is the Riemannian barycentre of } G(\bar{x},\sigma)\mbox{\,:} \hspace{0.5cm}   \bar{x} = \mathrm{argmin}_{\scriptscriptstyle x \in \mathcal{M}}\; E(x|\,\bar{x},\sigma)
   \end{equation}
   \begin{equation} \label{eq:asymp2}
      \sigma \, = \, \Phi \left(\, {\textstyle \bigintsss_{\scriptscriptstyle \mathcal{M}}} \;  
d^{\scriptscriptstyle \,2}(\bar{x},z)\, p(z|\,\bar{x},\sigma)\, dv(z) \right)
   \end{equation} 
 \end{subequations}
 where $\Phi : \mathbb{R}_{+} \rightarrow \mathbb{R}_+$ is the strictly increasing function introduced in Proposition \ref{prop:mle}.
\end{proposition}
 \vspace{0.1cm}
\textbf{Proof\,: } the proof can be carried out by direct generalisation of the one in~\cite{said} (proof of Proposition 9). An alternative, shorter proof of (\ref{eq:asymp1}), based on the definition of a Riemannian symmetric space, can be obtained by generalising~\cite{vemuri} (proof of Theorem 2.3.). \hfill$\blacksquare$\\[0.15cm]
\indent $\blacktriangleright$ \textbf{Remark 3} : Proposition \ref{prop:asymp} yields the consistency and asymptotic normality of the maximum likelihood estimate $\hat{x}_{\scriptscriptstyle N}$. Precisely,
\begin{subequations} \label{eq:consnorm}
\begin{equation} \label{eq:consnorm1}
 \mbox{consistency of } \hat{x}_{\scriptscriptstyle N} \mbox{ \,:}\hspace{0.5cm}  \lim_{\scriptscriptstyle N}\; \hat{x}_{\scriptscriptstyle N} \, = \, \bar{x}
\hspace{1.1cm}
\end{equation}
\begin{equation} \label{eq:consnorm2}
 \mbox{asymptotic normality of } \hat{x}_{\scriptscriptstyle N} \mbox{ \,:} \hspace{0.5cm} \sqrt{N}\;\, \mathrm{Log}_{\scriptscriptstyle \bar{x}}  (\hat{x}_{\scriptscriptstyle N}) \Rightarrow N_{\scriptscriptstyle d}(0\,,C^{\scriptscriptstyle -1}) 
\end{equation}
\end{subequations}
where $\mathrm{Log}$ denotes the Riemannian logarithm mapping\!~\cite{berger}, and $N_{\scriptscriptstyle d}$ a normal distribtuion in $d$ dimensions, ($d$ being the dimension of $\mathcal{M}$). Here,
\begin{equation} \label{eq:C}
  C \, = \, {\textstyle \bigintsss_{\scriptscriptstyle \mathcal{M}}} \;  
\left(\sigma^{\scriptscriptstyle -2}\,\mathrm{Log}_{\scriptscriptstyle \bar{x}}(z)\right) \otimes \left(\sigma^{\scriptscriptstyle -2}\,\mathrm{Log}_{\scriptscriptstyle \bar{x}}(z)\right)\, p(z|\,\bar{x},\sigma)\, dv(z)
\end{equation}
where $\otimes$ denotes exterior product of tangent vectors to $\mathcal{M}$. Indeed, (\ref{eq:consnorm1}) follows from Propositions \ref{prop:mle} and \ref{prop:asymp}. According to~\cite{bhatta1} (Theorem $2.3.$, Page $8$), the Riemannian barycentre $\hat{x}_{\scriptscriptstyle N}$ of independent samples $x_{\scriptscriptstyle 1},\ldots,x_{\scriptscriptstyle N}$ from $G(\bar{x},\sigma)$ converges to the Riemannien barycentre of $G(\bar{x},\sigma)$. But this is exactly $\bar{x}$. As for (\ref{eq:consnorm2}), it can be proved by direct generalisation of the proof in~\cite{said} (Proposition 11). \hfill$\blacksquare$

\section{Application to spaces of structured covariance matrices} \label{sec:matrices}
The previous Section \ref{sec:gaussian} introduced the notion of a Gaussian distribution on a Riemannian symmetric space of non-positive curvature, $\mathcal{M}$. In the present section, the development of Section \ref{sec:gaussian} is applied to the special case where $\mathcal{M}$ is a space of structured covariance matrices. Three variants are considered\,: $\mathcal{M} \, = \,$ space of complex covariance matrices, $\mathcal{M} \, = \,$ space of Toeplitz covariance matrices, and $\mathcal{M} \, = \,$ space of block-Toeplitz covariance matrices. In~\cite{vemuri}\cite{said}, Gaussian distributions were studied for $\mathcal{M} \, = \,$ space of real covariance matrices. Here, they are extended to spaces of covariance matrices with additional structure.

It seems that\!~\cite{barbaresco}\cite{barbaresco1}, in general, a space of structured covariance matrices becomes a Riemannian symmetric space of non-positive curvature, when it is equipped with a so-called Hessian metric. Roughly, a Hessian metric is a Riemannian metric arising from the entropy function,
\begin{equation} \label{eq:hessian1}
  \mbox{ Riemannian metric } = \mbox{ Hessian of entropy} 
\end{equation}
where, if $\mathcal{M}$ is a space of structured covariance matrices, the entropy of $x \in \mathcal{M}$ is minus its log-determinant, $H(x) = -\log \det(x)$~\cite{cover} (Chapter 8, Page $254$), as in the usual formula for the entropy of a multivariate normal distribution, with covariance $x$. 

For each of the three variants mentioned above, the following development will be carried out. First, (\ref{eq:hessian1}) is used to obtain a Riemannian metric $ds^{\scriptscriptstyle 2}$ on $\mathcal{M}$. It turns out that, as intended, this metric makes $\mathcal{M}$ into a Riemannian symmetric space of non-positive curvature. Second, Gaussian distributions on $\mathcal{M}$ are characterised, as in Section \ref{sec:gaussian}. Precisely, the graphs of the functions $\log Z : \mathbb{R}_+ \rightarrow \mathbb{R}$ and $\Phi : \mathbb{R}_{+} \rightarrow \mathbb{R}_+$ are given based on Propositions \ref{prop:z} and \ref{prop:mle}. Further, the reader is referred to Appendix \ref{app:sample}, for algorithms which provide samples from a Gaussian distribution on $\mathcal{M}$. 

\subsection{Complex covariance matrices} \label{subsec:2t2}
Here, $\mathcal{M}  = \mathcal{H}_{\scriptscriptstyle n\,}$, the space of $n \times n$ complex covariance matrices. Precisely, each element $Y \in \mathcal{H}_{\scriptscriptstyle n}$ is an $n \times n$ Hermitian positive definite matrix. This space $\mathcal{H}_{\scriptscriptstyle n}$ is equipped with a Hessian metric as follows. \\[0.1cm]
\indent \textbf{-- Hessian metric\,:} application of the general prescription (\ref{eq:hessian1}) leads to a Riemannian metric on $\mathcal{H}_{\scriptscriptstyle n}$ which is identical to the well-known affine-invariant metric, used in~\cite{bhatia}\cite{pennec2}. Indeed, the Hessian of the entropy function $H(Y) = -\log \det(Y)$ has the following expression,
\begin{subequations}
$$
\mbox{Hessian of entropy\,:} \hspace{0.5cm} D^{\scriptscriptstyle 2}H(Y)[v,w] \, = \, \mathrm{tr} \left[\, Y^{\scriptscriptstyle -1}v\,Y^{\scriptscriptstyle -1}w\,\right] \hspace{0.5cm} Y \in \mathcal{H}_{\scriptscriptstyle n}
$$
for any complex matrices $v$ and $w$, where $\mathrm{tr}$ denotes the trace. This expression can be found using the matrix differentiation formulae, as in~\cite{absil} (Appendix A, Page $196$),
\begin{equation} \label{eq:absil}
D \log \det(Y)[w] = \mathrm{tr}\left[\, Y^{\scriptscriptstyle -1}w\,\right] \hspace{1cm} DY^{\scriptscriptstyle -1}[v] = -\,Y^{\scriptscriptstyle -1}v\,Y^{\scriptscriptstyle -1}
\end{equation}
Replacing in (\ref{eq:hessian1}), gives the affine-invariant metric of $\mathcal{H}_{\scriptscriptstyle n}$,
\begin{equation} \label{eq:hnmetric}
  ds^{\scriptscriptstyle 2}_{\scriptscriptstyle Y}(dY) \, =  D^{\scriptscriptstyle 2}H(Y)[dY,dY] 
\, = \, \mathrm{tr}\left[\, Y^{\scriptscriptstyle -1}dY\,Y^{\scriptscriptstyle -1}dY\,\right] \hspace{0.5cm} 
Y \in \mathcal{H}_{\scriptscriptstyle n}
\end{equation}
The Riemannian distance corresponding to the metric (\ref{eq:hnmetric}) has the following expression~\cite{bhatia}\cite{pennec2},
\begin{equation} \label{eq:hndistance}
 d^{\scriptscriptstyle \,2}(X,Y) \, = \, \mathrm{tr}\left[\log\left(X^{\scriptscriptstyle -1/2}YX^{\scriptscriptstyle -1/2}\right)\right]^{\scriptscriptstyle 2} \hspace{0.5cm} X,Y \in \mathcal{H}_{\scriptscriptstyle n}
\end{equation}
where all matrix functions, (matrix power and logarithm), are Hermitian matrix functions~\cite{higham}. Moreover, the Riemannian volume element $dv$ is given by, (this can be checked using~\cite{nagar}, Lemma 2.1),
\begin{equation} \label{eq:hnvolume}
   dv(Y) = \det(Y)^{\scriptscriptstyle -n} \, \prod_{\scriptscriptstyle i\leq j} \mathrm{Re}\,dY_{\scriptscriptstyle ij} \,
   \prod_{\scriptscriptstyle i< j} \mathrm{Im}\,dY_{\scriptscriptstyle ij}
\end{equation}
where subscripts denote matrix elements, and $\mathrm{Re}$ and $\mathrm{Im}$ denote real and imaginary parts. \\[0.1cm]
\indent \textbf{-- Symmetric space\,:} equipped with the Riemannian metric (\ref{eq:hnmetric}),  $\mathcal{H}_{\scriptscriptstyle n}$ is a Riemannian symmetric space of non-positive curvature~\cite{helgason}. Accordingly, it is possible to apply Facts \ref{fact:1}--\ref{fact:5}. To do so, note that $G = GL(n,\mathbb{C})$, the Lie group of $n \times n$ invertible complex matrices, acts on $\mathcal{H}_{\scriptscriptstyle n}$ transitively and isometrically by congruence transformations~\cite{bhatia},
\begin{equation} \label{eq:congruence}
   g\cdot Y \, = \, g\,Yg^{\scriptscriptstyle H} \hspace{0.5cm} g \in G \, , \, Y \in \mathcal{H}_{\scriptscriptstyle n}
\end{equation}
where $^{\scriptscriptstyle H}$ denotes the conjugate transpose. Then, Facts \ref{fact:1} and \ref{fact:2} state that definitions
(\ref{eq:hnmetric})--(\ref{eq:hnvolume}) verify the general identities (\ref{eq:invmetric}) and (\ref{eq:invintegral}). Furthermore, Fact \ref{fact:3} states the existence and uniqueness of Riemannian barycentres in $\mathcal{H}_{\scriptscriptstyle n}$. 

The application of Facts \ref{fact:4} and \ref{fact:5} to $\mathcal{H}_{\scriptscriptstyle n}$ leads to the two following formulae, which are proved in Example \ref{ex:sln} of Appendix \ref{app:roots},
\begin{equation} \label{eq:hnintegralp}
  \int_{\scriptscriptstyle \mathcal{H}_{\scriptscriptstyle n}} \, f(Y)\, dv(Y) \, = \, \mathrm{C} \times \int_{\scriptscriptstyle U(n)} \int_{\scriptscriptstyle \mathbb{R}^n} \, f(r,U) \, \prod_{i < j} \sinh^{\scriptscriptstyle 2}(|r_{\scriptscriptstyle i} - r_{\scriptscriptstyle j}|/2) \, dr \, dU
\end{equation} 
\begin{equation} \label{eq:hnkilling}
 d^{\scriptscriptstyle \,2}(I,Y) \, = \, \Vert r \Vert^{\scriptscriptstyle 2} \, = \, \sum^{\scriptscriptstyle n}_{\scriptscriptstyle j=1} \, r^{\scriptscriptstyle 2}_{\scriptscriptstyle j} \hspace{0.75cm} (I = n \times n \mbox{ identity matrix})
\end{equation} 
For $\mathcal{H}_{\scriptscriptstyle n}$, formulae (\ref{eq:hnintegralp}) and (\ref{eq:hnkilling}) play the same role as formulae (\ref{eq:integralp}) and (\ref{eq:killing}), for a general symmetric space $\mathcal{M}$. Here $Y = U e^{\scriptscriptstyle r} U^{\scriptscriptstyle H}$ is the spectral decomposition of $Y$, with $U \in U(n)$, the group of $n \times n$ unitary matrices, and $e^{\scriptscriptstyle r} = \mathrm{diag}(e^{\scriptscriptstyle r_1},\ldots, e^{\scriptscriptstyle r_n})$ for $r = (r_{\scriptscriptstyle 1},\ldots,r_{\scriptscriptstyle n}) \in \mathbb{R}^{\scriptscriptstyle n}$. Moreover, $dr$ denotes the Lebesgue measure on $\mathbb{R}^{\scriptscriptstyle n}$ and $dU$ denotes the normalised Haar measure on $U(n)$. \\[0.1cm]
\indent \textbf{-- Gaussian distributions\,:} the functions $\log Z : \mathbb{R}_+ \rightarrow \mathbb{R}$ and $\Phi : \mathbb{R}_+ \rightarrow \mathbb{R}_+$, for Gaussian distributions on $\mathcal{H}_{\scriptscriptstyle n}$, can be obtained using formulae (\ref{eq:hnintegralp}) and (\ref{eq:hnkilling}). These yield the following expression for $Z(\sigma)$, similar to (\ref{eq:zgeneral}) in the proof of Proposition \ref{prop:z},
$$
 Z(\sigma) \,=\, \int_{\scriptscriptstyle \mathcal{H}_n}\, \exp \left[-\,\frac{d^{\,\scriptscriptstyle 2}(I,Y)}{2\sigma^{\,\scriptscriptstyle 2}}\right]
  \, dv(Y) \, = \,\mathrm{C} \times \int_{\scriptscriptstyle U(n)} \int_{\scriptscriptstyle \mathbb{R}^n} \,  
\exp \left[-\,\frac{\Vert r \Vert^{\scriptscriptstyle 2}}{2\sigma^{\,\scriptscriptstyle 2}}\right]\, \prod_{i < j} \sinh^{\scriptscriptstyle 2}(|r_{\scriptscriptstyle i} - r_{\scriptscriptstyle j}|/2) \, dr \, dU
$$
Since the function under the double integral does not depend on $U$, this integral simplifies to a formula similar to formula (\ref{eq:zint}) in the proof of Proposition \ref{prop:z},
\begin{equation} \label{eq:hnzint}
 Z(\sigma) \,=\, \mathrm{C} \times \int_{\scriptscriptstyle \mathbb{R}^n} \,  
\exp \left[-\,\frac{\Vert r \Vert^{\scriptscriptstyle 2}}{2\sigma^{\,\scriptscriptstyle 2}}\right]\, \prod_{i < j} \sinh^{\scriptscriptstyle 2}(|r_{\scriptscriptstyle i} - r_{\scriptscriptstyle j}|/2) \, dr 
\end{equation}
For any moderate value of $n$, (in practice, up to $n \approx 40$ has been considered), this formula can be numerically evaluated on a desktop computer, using a specifically designed Monte Carlo technique~\cite{paolo}, leading to the graph of $Z(\sigma)$. Then, obtaining the graphs of $\log Z : \mathbb{R}_+ \rightarrow \mathbb{R}$ and $\Phi : \mathbb{R}_+ \rightarrow \mathbb{R}_+$ is straightforward. In particular, it follows from (\ref{eq:legendre1}) in the proof of Proposition \ref{prop:mle} that the function $\Phi$ is obtained by solving, for each $\rho \in \mathbb{R}_{+}$, the equation $\rho = \psi^{\scriptscriptstyle \prime}(\eta)$. Denoting $\hat{\eta}$ the unique solution of this equation, it follows that $\Phi(\rho) = \hat{\sigma}$ where $\hat{\eta} = -1/2\hat{\sigma}^{\scriptscriptstyle 2}$. 

Figure \ref{fig:hn} shows the graphs of $\log Z$ and $\Phi$ when $n = 20$. The constant factor $\mathrm{C}$ appearing in (\ref{eq:hnzint}), which depends on $n$ through an involved expression (compare to~\cite{said}, Proposition 4), is ignored in these graphs. This has no influence on the graph of $\Phi$.
\begin{figure} 
\centering
\begin{subfigure}{.5\textwidth}
  \centering
  \includegraphics[width=1\linewidth]{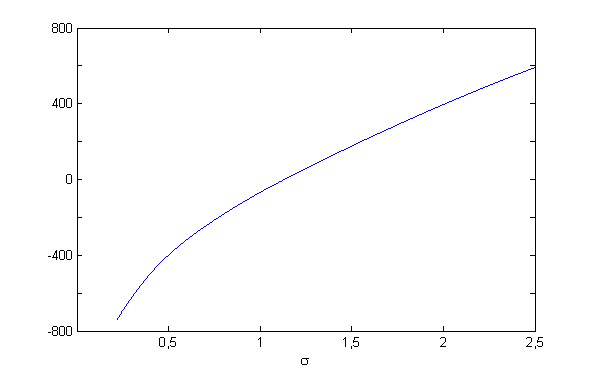}
  \caption{Graph of $\log Z : \mathbb{R}_+ \rightarrow \mathbb{R}$}
  \label{fig:hn1}
\end{subfigure}%
\begin{subfigure}{.5\textwidth}
  \centering
  \includegraphics[width=1\linewidth]{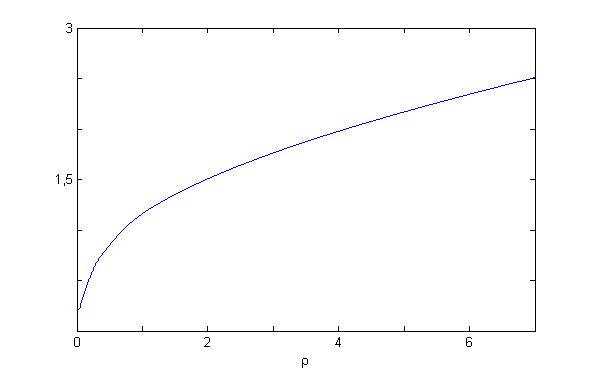}
  \caption{Graph of $\Phi:\mathbb{R}_+ \rightarrow \mathbb{R}_+$}
  \label{fig:hn2}
\end{subfigure}
\caption{$\log Z$ and $\Phi$ for $20 \times 20$ complex covariance matrices}
\label{fig:hn}
\end{figure}
\\[0.15cm]
\indent $\blacktriangleright$ \textbf{Remark 4} : the space $\mathcal{H}_{\scriptscriptstyle 2}$ of $2 \times 2$ complex covariance matrices plays a central role in polarimetric synthetic aperture radar (PolSAR) image analysis~\cite{polsar1}\cite{polsar2}. Each pixel in a PolSAR image is given by a covariance matrix $Y \in \mathcal{H}_{\scriptscriptstyle 2}$, which is the coherency matrix of a polarised lightwave. For the special case of the space $\mathcal{H}_{\scriptscriptstyle 2}$, the analytic expression of the function $Z(\sigma)$ can be found by integrating  (\ref{eq:hnzint}),
\begin{equation} \label{eq:h2z}
  Z(\sigma) \, = \, \mathrm{C} \, \times \, \pi \sigma^{\scriptscriptstyle 2} \, \left( e^{\scriptscriptstyle \sigma^2} - 1 \right)
\end{equation}
Formula (\ref{eq:h2z}) gives the same behavior as the function represented in Figure \ref{fig:hn1}. Indeed, it is clear from this formula that $\log Z(\sigma) \rightarrow -\infty$ when $\sigma \rightarrow 0$ and that $\log Z(\sigma) \rightarrow +\infty$ when $\sigma \rightarrow +\infty$. \hfill $\blacksquare$ 
\end{subequations}
\subsection{Toeplitz covariance matrices} \label{subsec:tn}
Here, $\mathcal{M} = \mathcal{T}_{\scriptscriptstyle n\,}$, the space of $n \times n$ Toeplitz covariance matrices. Precisely, each element $T \in \mathcal{T}_{\scriptscriptstyle n}$ is an $n \times n$ Hermitian positive definite matrix, which also has Toeplitz structure. Toeplitz covariance matrices are of fundamental importance in signal processing, since they are the autocovariance matrices of wide-sense stationary random signals~\cite{therrien}. A Hessian metric on $\mathcal{T}_{\scriptscriptstyle n}$ was introduced in~\cite{barbaresco}\cite{barbaresco1}. This metric is expressed in terms of so-called reflection or partial autocorrelation coefficients, as given by (\ref{eq:tnmetric}) below. 

It should be noted that $\mathcal{T}_{\scriptscriptstyle n}$ is a subspace of the space $\mathcal{H}_{\scriptscriptstyle n}$ of complex covariance matrices, studied in \ref{subsec:2t2}. As suggested in~\cite{bini}, this means $\mathcal{T}_{\scriptscriptstyle n}$ can be equipped with a Riemannian metric which is simply the restriction of the Riemannian metric (\ref{eq:hnmetric}) from $\mathcal{H}_{\scriptscriptstyle n\,}$. Unlike the Hessian metric of~\cite{barbaresco}\cite{barbaresco1}, this metric does not make $\mathcal{T}_{\scriptscriptstyle n}$ into a Riemannian symmetric space. Therefore, it cannot be used in the present paper. \\[0.15cm]
\indent $\blacktriangleright$ \textbf{Remark 5} : in signal processing, reflection coefficients arise naturally in linear prediction problems~\cite{therrien} (Chapter 8). Equivalently, they also arise from the study of orthogonal polynomials on the unit circle~\cite{simon}.  Here, it is suitable to recall their definition from this analytical point of view~\cite{simon} (see Chapter 1, Page $55$). Recall each element $T \in \mathcal{T}_{\scriptscriptstyle n}$ has Toeplitz structure. That is, there exist $r_{\scriptscriptstyle 0} > 0$ and $r_{\scriptscriptstyle 1},\ldots,r_{\scriptscriptstyle n-1} \in \mathbb{C}$ such that $T_{\scriptscriptstyle ij} = T^{\scriptscriptstyle *}_{\scriptscriptstyle ji} =  r_{\scriptscriptstyle i-j\,}$, (where $^*$ denotes the complex conjugate). Since $T$ is positive definite, it can be used to define a scalar product on polynomials of degree $\leq n-1$, which depend on a complex variable $z$ on the unit circle, $|z|^{\scriptscriptstyle 2} = zz^{\scriptscriptstyle *} = 1$. This scalar product is given by
\begin{subequations} \label{eq:opuc}
\begin{equation} \label{eq:opuc1}
(z^{\scriptscriptstyle i},z^{\scriptscriptstyle j}) = (z^{\scriptscriptstyle j},z^{\scriptscriptstyle i})^{\scriptscriptstyle *} = r_{\scriptscriptstyle i-j} 
\end{equation}
Let $\phi_{\scriptscriptstyle 0}, \ldots, \phi_{\scriptscriptstyle n-1}$ be orthogonal monic polynomials with respect to this scalar product (monic means the leading degree coefficient is equal to $1$). Reflection coefficients are the complex numbers $\alpha_{\scriptscriptstyle 1},\ldots, \alpha_{\scriptscriptstyle n-1}$ defined by the Szeg\"{o}-Levinson recursion,  
\begin{equation} \label{eq:szego}
  \alpha_{\scriptscriptstyle j} \, = \, -(1,z\phi_{\scriptscriptstyle j-1})/\delta_{\scriptscriptstyle j-1} \hspace{0.5cm} \phi_{\scriptscriptstyle j} \ = \, z\phi_{\scriptscriptstyle j-1} + \alpha_{\scriptscriptstyle j}\tilde{\phi}_{\scriptscriptstyle j-1} \hspace{0.5cm} \delta_{\scriptscriptstyle j} \, = \, \delta_{\scriptscriptstyle j-1}(1 - |\alpha_{\scriptscriptstyle j}|^{\scriptscriptstyle 2})
\end{equation}
where $\phi_{\scriptscriptstyle 0} = 1$ and $\delta_{\scriptscriptstyle 0} = r_{\scriptscriptstyle 0\,}$, and $\tilde{\phi}$ is the polynomial obtained from $\phi$ by conjugating and reversing the order of its coefficients. The significance of the numbers $\delta_{\scriptscriptstyle j}$ is that $\delta_{\scriptscriptstyle j} = (\phi_{\scriptscriptstyle j},\phi_{\scriptscriptstyle j})$, so $\delta_{\scriptscriptstyle j}$ must be strictly positive. From the third formula in (\ref{eq:szego}), this implies that each reflection coefficient $\alpha_{\scriptscriptstyle j}$ verifies $|\alpha_{\scriptscriptstyle j}| < 1$. Reflection coefficients allow a useful expression of the determinant of the matrix $T$. Indeed, since the polynomials $\phi_{\scriptscriptstyle j}$ are orthogonal with respect to the scalar product defined by $T$,
\begin{equation} \label{eq:toepdet}
 \det(T) = \prod^{\scriptscriptstyle n-1}_{j=0} (\phi_{\scriptscriptstyle j},\phi_{\scriptscriptstyle j}) = r_{\scriptscriptstyle 0} \prod^{\scriptscriptstyle n-1}_{j=1} \delta_{\scriptscriptstyle j} = r^{\scriptscriptstyle n}_{\scriptscriptstyle 0} \prod^{\scriptscriptstyle n-1}_{j=1} \,(1 - |\alpha_{\scriptscriptstyle j}|^{\scriptscriptstyle 2}\,)^{\scriptscriptstyle\, n-j}
\end{equation}
\end{subequations}
This will be useful in obtaining the Hessian metric (\ref{eq:tnmetric}) below. 
\hfill $\blacksquare$ \\[0.1cm]
\indent \textbf{-- Hessian metric\,:} consider now the Hessian metric on $\mathcal{T}_{\scriptscriptstyle n}$, introduced in~\cite{barbaresco}\cite{barbaresco1}.  This is based on the following identification,
\begin{equation} \label{eq:tnssdecomposition}
 \mathcal{T}_{\scriptscriptstyle n} \,=\, \mathbb{R}_{+} \,\times\, \mathbb{D} \,\times\, \ldots \,\times \mathbb{D} \hspace{0.5cm} (n \mbox{ factors}) 
\end{equation}
where $\mathbb{D}$ denotes the set of complex $\alpha$ which verify $|\alpha| < 1$. Precisely~\cite{schou}\cite{ramsey}, $\mathcal{T}_{\scriptscriptstyle n}$ is diffeomorphic to the product space on the right hand side, and the required diffeomorphism is given in the notation of  (\ref{eq:opuc}) by $T \mapsto ( r_{\scriptscriptstyle 0}, \alpha_{\scriptscriptstyle 1},\ldots,\alpha_{\scriptscriptstyle n-1})$. The advantage of parameterising $T$ by means of $( r_{\scriptscriptstyle 0}, \alpha_{\scriptscriptstyle 1},\ldots,\alpha_{\scriptscriptstyle n-1})$ instead of $(r_{\scriptscriptstyle 0},r_{\scriptscriptstyle 1},\ldots, r_{\scriptscriptstyle n-1})$ is that the latter are subject to a complicated set of non-linear constraints, expressing the fact that $T$ is positive definite, whereas the former are subject to no constraints, except $r_{\scriptscriptstyle 0} \in \mathbb{R}_+$ and $\alpha_{\scriptscriptstyle j} \in \mathbb{D}$. In the following, $r_{\scriptscriptstyle 0}$ is denoted $r$, since there is no risk of confusion.

Using (\ref{eq:tnssdecomposition}), the Hessian metric on $\mathcal{T}_{\scriptscriptstyle n}$ can be found directly. Note first the entropy function $H(T) = -\log \det(T)$ follows from (\ref{eq:toepdet}),
\begin{subequations}
$$
 \mbox{entropy function\,:} \hspace{0.5cm}  H(T) \,=\,  -n\log r - \, \sum^{\scriptscriptstyle n-1}_{j=1} (n-j) \, \log(1-|\alpha_{\scriptscriptstyle j}|^{\scriptscriptstyle 2})
$$
Note the following elementary calculations,
$$
\frac{\partial^{\scriptscriptstyle 2} H}{\partial r^{\scriptscriptstyle 2}} = n\, \frac{1}{r^{\scriptscriptstyle 2}} \hspace{0.5cm}
\frac{\partial^{\scriptscriptstyle 2} H}{\partial \alpha_{\scriptscriptstyle j} \, \partial \alpha^{\scriptscriptstyle *}_{\scriptscriptstyle j}} =
\, (n-j)\,\left[1 - |\alpha_{\scriptscriptstyle j}|^{\scriptscriptstyle 2}\,\right]^{\scriptscriptstyle\,-2}
$$
These immediately yield the Riemannian metric~\cite{barbaresco}\cite{barbaresco1},
\begin{eqnarray} \label{eq:tnmetric}
  ds^{\scriptscriptstyle 2}_{\scriptscriptstyle T}(dr,d\alpha) \,=\, n\, \left(\frac{dr}{r}\right)^{\scriptscriptstyle 2} \,+\, 
   \sum^{\scriptscriptstyle n-1}_{j=1} (n-j)\, ds^{\scriptscriptstyle 2}_{\scriptscriptstyle\alpha_ j}(d\alpha_{\scriptscriptstyle j}) \\[0.1cm]
   \nonumber \mbox{where } \;\; ds^{\scriptscriptstyle 2}_{\scriptscriptstyle\alpha}(d\alpha) \,=\, \left.|d\alpha|^{\scriptscriptstyle 2}\middle/\left[1 - |\alpha|^{\scriptscriptstyle 2}\,\right]^{\scriptscriptstyle\,2} \right. \;\;\; \alpha \in \mathbb{D} \;\;\,\;\,\,
\end{eqnarray}
Here, $ds^{\scriptscriptstyle 2}_{\scriptscriptstyle\alpha}$ is recognised as the Poincar\'e metric on $\mathbb{D}$~\cite{terras1} (Chapter III). Accordingly, the Riemannian distance and Riemannian volume element on $\mathcal{T}_{\scriptscriptstyle n}$ have the following expressions,
\begin{eqnarray} \label{eq:tndistance}
  d^{\scriptscriptstyle 2}(T^{\scriptscriptstyle (1)},T^{\scriptscriptstyle (2)}) \,=\, n\, \left| \log\left(r^{\scriptscriptstyle (2)}\,\right) - \log\left(r^{\scriptscriptstyle (1)}\,\right)\,\right|^{\scriptscriptstyle 2} \,+\, 
   \sum^{\scriptscriptstyle n-1}_{j=1} \,(n-j)\, d^{\scriptscriptstyle\, 2}_{\scriptscriptstyle \mathbb{D}}(\alpha^{\scriptscriptstyle (1)}_{\scriptscriptstyle j},\alpha^{\scriptscriptstyle (2)}_{\scriptscriptstyle j}) \\[0.1cm]
   \nonumber \mbox{where } \;\; d_{\scriptscriptstyle \mathbb{D}}(\alpha,\beta) \,=\, \,\,\mathrm{atanh}\left|\frac{\alpha - \beta}{1 - \alpha^{\scriptscriptstyle *}\beta }\right|\;\;\; \alpha,\beta \in \mathbb{D} \hspace{3.1cm}\,\,\,\,\,
\end{eqnarray}
\begin{eqnarray} \label{eq:tnvolume}
  dv(T) \,=\, {\scriptstyle \sqrt{n}\prod^{n-1}_{j=1} (n-j)} \, \times \, \frac{dr}{r} \, \prod^{\scriptscriptstyle n-1}_{j=1} \, dv(\alpha_{\scriptscriptstyle j})  \hspace{1.4cm} \\[0.1cm]
   \nonumber \mbox{where } \;\; dv(\alpha) \,=\, \frac{\mathrm{Re}\, d\alpha \, \mathrm{Im} \, d\alpha}{\strut\left[1 - |\alpha|^{\scriptscriptstyle 2}\,\right]^{\scriptscriptstyle\,2}}\;\;\; \alpha \in \mathbb{D} \;\;\,\,\, \hspace{1.7cm}
\end{eqnarray}
Here, $d^{\scriptscriptstyle\, 2}_{\scriptscriptstyle \mathbb{D}}(\alpha,\beta)$ and $dv(\alpha)$ denote the Riemannian distance and Riemannian volume element on $\mathbb{D}$, corresponding to the Poincar\'e metric $ds^{\scriptscriptstyle 2}_{\scriptscriptstyle\alpha}$ appearing in (\ref{eq:tnmetric})~\cite{terras1}. \\[0.1cm]
\indent \textbf{-- Symmetric space\,:} to see that the Riemannian metric (\ref{eq:tnmetric}) makes $\mathcal{T}_{\scriptscriptstyle n}$ into a Riemannian symmetric space of non-positive curvature, it is enough to notice that decomposition (\ref{eq:tnssdecomposition}) is of the form (\ref{eq:ssdecomposition}), as in the proof of Proposition \ref{prop:z} in \ref{subsec:gaussiandef}. Indeed, in this decomposition, $\mathbb{R}_+$ appears as a one-dimensional Euclidean space, by identification with $\mathbb{R}$ through the logarithm function, while each copy of $\mathbb{D}$ appears as a scaled version of the Poincar\'e disc, with the scale factor $(n-j)^{\scriptscriptstyle 1/2\,}$, and therefore as a Riemannian symmetric space of non-compact type~\cite{terras1} (Chapter III), (the Poincar\'e disc is the set $\mathbb{D}$ with the Poincar\'e metric $ds^{\scriptscriptstyle 2}_{\scriptscriptstyle\alpha\,}$, and is a space of constant negative curvature).

To apply Facts \ref{fact:1}--\ref{fact:5} to $\mathcal{T}_{\scriptscriptstyle n\,}$, consider the group $G = \mathbb{R}_+ \, \times \, SU(1,1) \, \times \, \ldots \, \times SU(1,1)$, where $SU(1,1)$ is the group of $2 \times 2$ complex matrices $u$ such that
\begin{equation} \label{eq:gsieg}
  u \, = \, \left(\begin{array}{cc} a & b \\ c & d \end{array}\right): \hspace{0.35cm}
  u^{\scriptscriptstyle T} \left(\begin{array}{cc} 0 & -1 \\ 1 & 0 \end{array}\right) u \;= \left(\begin{array}{cc} 0 & -1 \\ 1 & 0 \end{array}\right) \hspace{0.25cm}\, , \hspace{0.25cm}  
  u^{\scriptscriptstyle H} \left(\begin{array}{cc} -1 & 0 \\ 0 & 1 \end{array}\right) u \;= \left(\begin{array}{cc} -1 & 0 \\ 0 & 1 \end{array}\right)
\end{equation}
with $a,b,c,d \in \mathbb{C}$ and $^{\scriptscriptstyle T}$ denotes the transpose. This group acts transitively and isometrically on $\mathcal{T}_{\scriptscriptstyle n}$ in the following way,
\begin{eqnarray} \label{eq:tngroup}
  g \cdot T \,=\,  (s\,r,u_{\scriptscriptstyle 1}\cdot \alpha_{\scriptscriptstyle 1}, \ldots, u_{\scriptscriptstyle n-1}\cdot \alpha_{\scriptscriptstyle n-1})
  \hspace{0.25cm},\hspace{0.25cm} g = (s,u_{\scriptscriptstyle 1},\ldots, u_{\scriptscriptstyle n-1})
  \\[0.1cm]
   \nonumber \mbox{where } \;\; u \cdot \alpha \,=\, \frac{a\,\alpha + b }{c \,\alpha + d} \hspace{0.2cm} u \in SU(1,1) \, , \, \alpha \in \mathbb{D}\;\;\,\,\hspace{2cm}\,
\end{eqnarray}
Then, Facts \ref{fact:1} and \ref{fact:2} state that definitions (\ref{eq:tnmetric})--(\ref{eq:tnvolume}) verify the general identities (\ref{eq:invmetric}) and (\ref{eq:invintegral}). Furthermore, Fact \ref{fact:3} states the existence and uniqueness of Riemannian barycentres in $\mathcal{T}_{\scriptscriptstyle n}$. 

The application of Facts \ref{fact:4} and \ref{fact:5} to $\mathcal{T}_{\scriptscriptstyle n}$ leads to the two following formulae, which are proved in Example \ref{ex:sp1} of Appendix \ref{app:roots},
\begin{eqnarray} \label{eq:tnintegralp}
  \int_{\scriptscriptstyle \mathcal{T}_{\scriptscriptstyle n}} \, f(T)\, dv(T) \, = \, \mathrm{C} \times \int^{\scriptscriptstyle +\infty}_{\scriptscriptstyle -\infty} \,
\int_{\scriptscriptstyle \mathbb{D}} \ldots
\int_{\scriptscriptstyle \mathbb{D}}
\, f(t,\alpha_{\scriptscriptstyle 1},\ldots,\alpha_{\scriptscriptstyle n-1})  \, dv(\alpha_{\scriptscriptstyle 1}) \ldots dv(\alpha_{\scriptscriptstyle n-1}) \,  dt \\[0.1cm]
\nonumber \int_{\scriptscriptstyle \mathbb{D}} \, f(\alpha)\, dv(\alpha) \, = \, \mathrm{C} \times \int^{\scriptscriptstyle 2\pi}_{\scriptscriptstyle 0} \int^{\scriptscriptstyle +\infty}_{\scriptscriptstyle -\infty} \, f(\rho,\theta) \, \sinh(|\rho|)  d\rho \, d\theta \hspace{3.3cm} 
\end{eqnarray}
\begin{equation} \label{eq:tnkilling}
 d^{\scriptscriptstyle \,2}(I,T) \, = \, n \,t^{\scriptscriptstyle 2} + \sum^{\scriptscriptstyle n-1}_{\scriptscriptstyle n=1} (n-j)\rho^{\scriptscriptstyle \,2}_{\scriptscriptstyle j} \hspace{0.75cm} (I = n \times n \mbox{ identity matrix})\, \hspace{2.75cm}
\end{equation} 
For $\mathcal{T}_{\scriptscriptstyle n}$, formulae (\ref{eq:tnintegralp}) and (\ref{eq:tnkilling}) play the same role as formulae (\ref{eq:integralp}) and (\ref{eq:killing}), for a general symmetric space $\mathcal{M}$. Here $r = e^{\scriptscriptstyle t}$ and $\alpha_{\scriptscriptstyle j} = \mathrm{tanh}(\rho_{\scriptscriptstyle j})\,e^{\scriptscriptstyle \mathrm{i}\theta_{\scriptscriptstyle j}}$ for $j = 1,\ldots, n-1$, where $\mathrm{i} = \sqrt{-1}$. \\[0.1cm]
\indent \textbf{-- Gaussian distributions\,:} using formulae (\ref{eq:tnintegralp}) and (\ref{eq:tnkilling}), it is possible to find the exact analytic expression of the function $\log Z : \mathbb{R}_+ \rightarrow \mathbb{R}$ for Gaussian distributions on $\mathcal{T}_{\scriptscriptstyle n}$. From this expression, which is given in (\ref{eq:tnz}) below, the function $\Phi : \mathbb{R}_+ \rightarrow \mathbb{R}_+$ is found in the same way discussed in \ref{subsec:2t2}.

Unlike Gaussian distributions on $\mathcal{H}_{\scriptscriptstyle n}$, for which the function $\log Z$ cannot be found analytically, but only using the Monte Carlo technique of~\cite{paolo}, Gaussian distributions on $\mathcal{T}_{\scriptscriptstyle n}$ admit an exact analytic expression of $\log Z$ for any value of $n$. This means there is no computational limitation on the size (given by $n$) of Toeplitz covariance matrices which can be modelled using Gaussian distributions. 

Consider the calculation of $Z(\sigma)$ using (\ref{eq:tnintegralp}) and (\ref{eq:tnkilling}). Recall from (\ref{eq:zgeneral}) in the proof of Proposition 
\ref{prop:z}, 
$$
 Z(\sigma) \,=\, \int_{\scriptscriptstyle \mathcal{T}_n}\, \exp \left[-\,\frac{d^{\,\scriptscriptstyle 2}(I,T)}{2\sigma^{\,\scriptscriptstyle 2}}\right]
  \, dv(T) 
$$
Upon replacing (\ref{eq:tnintegralp}) and (\ref{eq:tnkilling}), this becomes,
$$
 Z(\sigma) \,=\, \mathrm{C} \times \int^{\scriptscriptstyle +\infty}_{\scriptscriptstyle -\infty} \,
\exp\left[ -\frac{n \,t^{\scriptscriptstyle 2}}{2\sigma^{\scriptscriptstyle 2}}\right] \, dt \, \times \, \prod^{\scriptscriptstyle n-1}_{\scriptscriptstyle j=1} \,\,
 \int^{\scriptscriptstyle 2\pi}_{\scriptscriptstyle 0} \int^{\scriptscriptstyle +\infty}_{\scriptscriptstyle -\infty} \, \exp\left[ -\frac{(n-j) \,\rho^{\scriptscriptstyle 2}_{\scriptscriptstyle j}}{2\sigma^{\scriptscriptstyle 2}}\right] 
 \, \sinh(|\rho_{\scriptscriptstyle j}|)  d\rho_{\scriptscriptstyle j} \, d\theta_{\scriptscriptstyle j}
$$
by computing each integral in the product, it follows that
\begin{eqnarray} \label{eq:tnz0}
 Z(\sigma) \,=\, \mathrm{C} \times \sqrt{2\pi/n} \; \sigma \, \times \, \prod^{\scriptscriptstyle n-1}_{\scriptscriptstyle j=1} \, Z_{\scriptscriptstyle \mathbb{D}}\left(\sigma\middle/\sqrt{n-j}\,\right) \\[0.1cm]
\nonumber  Z_{\scriptscriptstyle \mathbb{D}}(\sigma) \, = \, \left( 2\pi \right)^{\scriptscriptstyle 2/3} \, \sigma \, e^{\scriptscriptstyle \frac{\sigma^2}{2}} \times
\mathrm{erf}\left(\frac{\sigma}{\sqrt{2}}\right) \hspace{1.35cm}
\end{eqnarray}
where $\mathrm{erf}$ denotes the error function~\cite{aba}. Finally, this yields the analytic expression of $\log Z$,
\begin{equation} \label{eq:tnz}
 \log Z(\sigma) \,=\, \mathrm{C} \, + \, \log\left( \sigma\middle/\sqrt{\scriptstyle n}\,\right) \, + \,  \sum^{\scriptscriptstyle n-1}_{\scriptscriptstyle j=1} \, \log Z_{\scriptscriptstyle \mathbb{D}}
 \left( \sigma\middle/\sqrt{\scriptstyle n-j}\,\right)
\end{equation}
Using (\ref{eq:tnz0}) and (\ref{eq:tnz}), Figure \ref{fig:tn} shows the graphs of $\log Z$ and $\Phi$ when $n = 20$. These graphs ignore the constant $\mathrm{C}$ appearing in (\ref{eq:tnz}).
\begin{figure} 
\centering
\begin{subfigure}{.5\textwidth}
  \centering
  \includegraphics[width=1\linewidth]{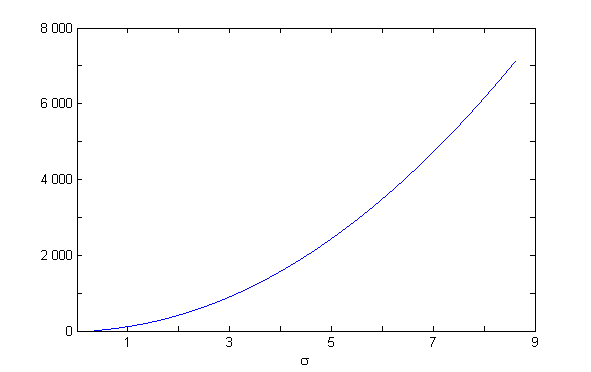}
  \caption{Graph of $\log Z : \mathbb{R}_+ \rightarrow \mathbb{R}$}
  \label{fig:tn1}
\end{subfigure}%
\begin{subfigure}{.5\textwidth}
  \centering
  \includegraphics[width=1\linewidth]{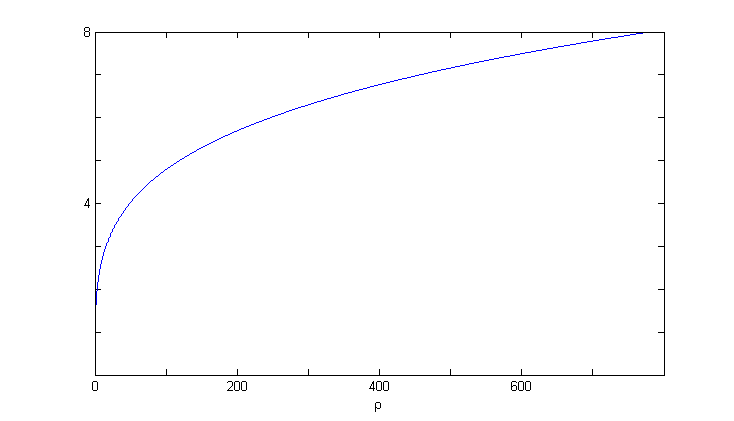}
  \caption{Graph of $\Phi:\mathbb{R}_+ \rightarrow \mathbb{R}_+$}
  \label{fig:tn2}
\end{subfigure}
\caption{$\log Z$ and $\Phi$ for $20 \times 20$ Toeplitz covariance matrices}
\label{fig:tn}
\end{figure}
 \end{subequations}
\subsection{Block-Toeplitz covariance matrices} \label{subsec:btn}
Here, $\mathcal{M} = \mathcal{T}^{\scriptscriptstyle N}_{\scriptscriptstyle n\,}$, the space of $n \times n$ block-Toeplitz covariance matrices which have $N \times N$ Toeplitz blocks. Precisely, each element $T \in \mathcal{T}^{\scriptscriptstyle N}_{\scriptscriptstyle n\,}$ is an $nN \times nN$ Hermitian positive definite matrix, which has a block structure of the following form
$$
T \, = \left( \begin{array}{llllll} {\scriptstyle T_{\scriptscriptstyle 0}} & \hspace{-0.0cm} {\scriptstyle T^{\scriptscriptstyle *}_{\scriptscriptstyle 1}} & &  &  \dots &   {\scriptstyle T^{\scriptscriptstyle *}_{\scriptscriptstyle n-1}} \\[0.1cm] {\scriptstyle T_{\scriptscriptstyle 1}} & \hspace{-0.0cm} {\scriptstyle T_{\scriptscriptstyle 0}} & \hspace{-0.0cm} {\scriptstyle T^{\scriptscriptstyle *}_{\scriptscriptstyle 1}} &  & \dots & {\scriptstyle T^{\scriptscriptstyle *}_{\scriptscriptstyle n-2}} \\[0.1cm]& \hspace{-0.0cm} {\scriptstyle T_{\scriptscriptstyle 1}} & \hspace{-0.0cm} {\scriptstyle T_{\scriptscriptstyle 0}} & & & \\
\vdots & \vdots& &  &  & \vdots \\[0.1cm]
{\scriptstyle T_{\scriptscriptstyle \tiny n-1}} & \hspace{-0.0cm} {\scriptstyle T_{\scriptscriptstyle \tiny n-2}} &  &  & \dots  & {\scriptstyle T_{\scriptscriptstyle 0}}  \end{array}\! \right) \hspace{1cm} T_{\scriptscriptstyle 0},\ldots, T_{\scriptscriptstyle n-1} \in \mathcal{T}_{\scriptscriptstyle N}
$$
Block-Toeplitz covariance matrices which have Toeplitz blocks play a fundamental role in signal and image processing, where they arise in multi-channel and two-dimensional linear prediction and filtering problems~\cite{tbbt1}\cite{tbbt2}. A Hessian metric on  $\mathcal{T}^{\scriptscriptstyle N}_{\scriptscriptstyle n\,}$ was considered in~\cite{jeuris1}, see (\ref{eq:btnmetric}) below. This generalises the Hessian metric (\ref{eq:tnmetric}) on the space $\mathcal{T}_{\scriptscriptstyle n}$ of Toeplitz covariance matrices, considered in the previous Paragraph \ref{subsec:tn}. Precisely, in the same way as the metric (\ref{eq:tnmetric}) was expressed in terms of scalar reflection coefficients $\alpha_{\scriptscriptstyle 1},\ldots,\alpha_{\scriptscriptstyle n-1}$ which belong to the complex unit disc $\mathbb{D}$, the metric (\ref{eq:btnmetric}) will be expressed in terms of matrix reflection coefficients $\Omega_{\scriptscriptstyle 1},\ldots,\Omega_{\scriptscriptstyle n-1}$ which belong to the matrix unit disc $\mathbb{D}_{\scriptscriptstyle N}$\,: the set of $N \times N$ symmetric complex matrices $\Omega$ which verify $I - \Omega^{\scriptscriptstyle *}\Omega \succ 0$, where $I$ denotes the $N \times N$ identity matrix and $\succ$ the Loewner order, ($X \succ Y$ means that $X - Y$ is positive definite)~\cite{bhatia}.
\\[0.1cm]
\indent $\blacktriangleright$ \textbf{Remark 6} :  in signal and image processing, matrix reflection coefficients arise naturally in multi-channel and two-dimensional linear prediction problems~\cite{tbbt1}\cite{tbbt2}. Here, they are introduced through the formalism of orthogonal matrix polynomials on the unit circle~\cite{damanik}\cite{holger} (compare to Remark 5, Paragraph \ref{subsec:tn}). A matrix polynomial $P$ on the unit circle is an expression $P(z) = P_{\scriptscriptstyle 0} +\ldots + P_{\scriptscriptstyle m}\, z^{\scriptscriptstyle m}$, where $z$ is a complex variable on the unit circle, and $P_{\scriptscriptstyle 0}, \ldots, P_{\scriptscriptstyle m}$ are complex $N \times N$ matrices (here, $m  \leq n-1$). Following~\cite{damanik} (Section 3, Page $43$), consider for a fixed $T \in \mathcal{T}^{\scriptscriptstyle N}_{\scriptscriptstyle n}$ the so-called left and right matrix scalar products of matrix polynomials. These are given by, 
\begin{subequations} \label{eq:mopuc}
\begin{equation} \label{eq:mopuc1}
  \langle P\, z^{\scriptscriptstyle i} \,,\, Qz^{\scriptscriptstyle j} \, \rangle_{\scriptscriptstyle L} \, = \, P\, T_{\scriptscriptstyle i-j} \,Q^{\scriptscriptstyle H} \hspace{0.5cm} \langle P\, z^{\scriptscriptstyle i} \,,\, Qz^{\scriptscriptstyle j} \, \rangle_{\scriptscriptstyle R} \, = \, P^{\scriptscriptstyle H}\, T_{\scriptscriptstyle i-j} \,Q
\end{equation}
for any complex $N \times N$ matrices $P$ and $Q\,$, where $T_{\scriptscriptstyle i-j} = T^{\scriptscriptstyle *}_{\scriptscriptstyle j-i\,}$. Let $\varphi^{\scriptscriptstyle L}_{\scriptscriptstyle 0},\ldots, \varphi^{\scriptscriptstyle L}_{\scriptscriptstyle n-1}$ and $\varphi^{\scriptscriptstyle R}_{\scriptscriptstyle 0},\ldots, \varphi^{\scriptscriptstyle R}_{\scriptscriptstyle n-1}$ be orthonormal matrix polynomials with respect to these matrix scalar products. Matrix reflection coefficients $\Omega_{\scriptscriptstyle 1},\ldots,\Omega_{\scriptscriptstyle n-1}$ are given by the generalised Szeg\"{o}-Levinson recursion,
\begin{equation} \label{eq:mopuc2}
 z \varphi^{\scriptscriptstyle L}_{\scriptscriptstyle j - 1} - \rho_{\scriptscriptstyle j} \varphi^{\scriptscriptstyle L}_{\scriptscriptstyle j} = \widetilde{\varphi^{\scriptscriptstyle R}_{\scriptscriptstyle j - 1}\Omega_{\scriptscriptstyle j}} \hspace{0.5cm}
 z \varphi^{\scriptscriptstyle R}_{\scriptscriptstyle j - 1} - \varphi^{\scriptscriptstyle R}_{\scriptscriptstyle j}\rho^{\scriptscriptstyle *}_{\scriptscriptstyle j} = \widetilde{\Omega_{\scriptscriptstyle j} \varphi^{\scriptscriptstyle L}_{\scriptscriptstyle j - 1}}
\end{equation}
where $\rho_{\scriptscriptstyle j} = (I - \Omega^{\scriptscriptstyle *}_{\scriptscriptstyle j}\,\Omega^{\phantom *}_{\scriptscriptstyle j})^{\scriptscriptstyle 1/2}$, and where $\tilde{P}$, for a matrix polynomial $P$, denotes the polynomial obtained by taking the conjugate transposes of the coefficients of $P$ and reversing their order. It follows from~\cite{holger} (Definition 2.3., Page $1615$), that $\Omega_{\scriptscriptstyle 1},\ldots,\Omega_{\scriptscriptstyle n-1} \in \mathbb{D}_{\scriptscriptstyle N}$. Moreover,  $\Omega_{\scriptscriptstyle 1},\ldots,\Omega_{\scriptscriptstyle n-1}$ allow a useful expression of the determinant of the matrix $T$~\cite{jeuris1} (Page $10$), 
\begin{equation} \label{eq:mopuc3}
  \det(T) \, =  \,\det(T_{\scriptscriptstyle 0})^{\scriptscriptstyle n} \, \prod^{\scriptscriptstyle n-1}_{j=1} \,\det(I - \Omega^{\scriptscriptstyle *}_{\scriptscriptstyle j}\Omega^{\phantom *}_{\scriptscriptstyle j})^{\scriptscriptstyle\, n-j}
\end{equation}
This will be useful in obtaining the Hessian metric (\ref{eq:btnmetric}) below. \hfill $\blacksquare$
\end{subequations}
 \\[0.1cm]
\indent \textbf{-- Hessian metric\,:} consider now the Hessian metric on $\mathcal{T}^{\scriptscriptstyle N}_{\scriptscriptstyle n}$, introduced in~\cite{jeuris1}.  This is based on the following identification,
\begin{equation} \label{eq:bttnssdecomposition}
 \mathcal{T}^{\scriptscriptstyle N}_{\scriptscriptstyle n} \,=\, \mathcal{T}_{\scriptscriptstyle N} \,\times\, \mathbb{D}_{\scriptscriptstyle N} \,\times\, \ldots \,\times \mathbb{D}_{\scriptscriptstyle N} \hspace{0.5cm} (n \mbox{ factors}) 
\end{equation}
where, as already stated, $\mathbb{D}_{\scriptscriptstyle N}$ denotes the set of $N \times N$ symmetric complex matrices $\Omega$ which verify $I - \Omega^{\scriptscriptstyle *}\Omega \succ 0$. Precisely~\cite{jeuris1}, $\mathcal{T}^{\scriptscriptstyle N}_{\scriptscriptstyle n}$ is diffeomorphic to the product space on the right hand side, and the required diffeomorphism is given by $T \mapsto ( T_{\scriptscriptstyle 0}, \Omega_{\scriptscriptstyle 1},\ldots,\Omega_{\scriptscriptstyle n-1})$, in the notation of  (\ref{eq:mopuc}). In the following, $T_{\scriptscriptstyle 0}$ is denoted $P$, in order to avoid the unnecessary subscript.

Using (\ref{eq:bttnssdecomposition}), the Hessian metric on $\mathcal{T}_{\scriptscriptstyle n}$ can be found directly~\cite{jeuris1}. Note first the entropy function $H(T) = -\log \det(T)$ follows from (\ref{eq:mopuc3}),
\begin{subequations}
$$
 \mbox{entropy function\,:} \hspace{0.5cm}  H(T) \,=\,  -n\log \det(P) - \, \sum^{\scriptscriptstyle n-1}_{j=1} (n-j) \, \log \det(I - \Omega^{\scriptscriptstyle *}_{\scriptscriptstyle j}\Omega^{\phantom *}_{\scriptscriptstyle j})
$$
Then~\cite{jeuris1} (Page $10$), applying the matrix differentiation formulae (\ref{eq:absil}) yields the required Hessian metric
\begin{eqnarray} \label{eq:btnmetric}
  ds^{\scriptscriptstyle 2}_{\scriptscriptstyle T}(dP,d\Omega) \,=\, n\, ds^{\scriptscriptstyle 2}_{\scriptscriptstyle P}(dP) \,+\, 
   \sum^{\scriptscriptstyle n-1}_{j=1} (n-j)\, ds^{\scriptscriptstyle 2}_{\scriptscriptstyle\Omega_ j}(d\Omega_{\scriptscriptstyle j})\hspace{2.3cm}\, \\[0.1cm]
   \nonumber \mbox{where } \;\; ds^{\scriptscriptstyle 2}_{\scriptscriptstyle\Omega}(d\Omega) \,=\,  \mathrm{tr}\left[(I - 
   \Omega\Omega^{\scriptscriptstyle *})^{\scriptscriptstyle -1\,} d\Omega \, (I - \Omega^{\scriptscriptstyle *}\Omega)^{\scriptscriptstyle -1\,} d\Omega^{\scriptscriptstyle *}\, \right]  \;\;\; \Omega \in \mathbb{D}_{\scriptscriptstyle N} \;\;\,\,
\end{eqnarray}
Here, $ds^{\scriptscriptstyle 2}_{\scriptscriptstyle P}$ is the Riemannian metric (\ref{eq:tnmetric}) of the space $\mathcal{T}_{\scriptscriptstyle N}$, evaluated at $P \in \mathcal{T}_{\scriptscriptstyle N}$. On the other hand, $ds^{\scriptscriptstyle 2}_{\scriptscriptstyle\Omega}$ is recognised as the Siegel metric on $\mathbb{D}_{\scriptscriptstyle N}$\!~\cite{siegel}.  Accordingly, the Riemannian distance and Riemannian volume element on $\mathcal{T}_{\scriptscriptstyle n}$ have the following expressions,
\begin{eqnarray} \label{eq:btndistance}
  d^{\scriptscriptstyle \,2}(T^{\scriptscriptstyle (1)},T^{\scriptscriptstyle (2)}) \,=\, n\, d^{\scriptscriptstyle 2}_{\scriptscriptstyle \mathcal{T}_{\scriptscriptstyle N}}(P^{\scriptscriptstyle (1)},P^{\scriptscriptstyle (2)}) \,+\, 
   \sum^{\scriptscriptstyle n-1}_{j=1} \,(n-j)\, d^{\scriptscriptstyle\, 2}_{\scriptscriptstyle \mathbb{D}_{\scriptscriptstyle N}}(\Omega^{\scriptscriptstyle (1)}_{\scriptscriptstyle j},\Omega^{\scriptscriptstyle (2)}_{\scriptscriptstyle j}) \hspace{5.4cm} \\[0.1cm]
   \nonumber \mbox{where } \;\; d^{\scriptscriptstyle\, 2}_{\scriptscriptstyle \mathbb{D}_{\scriptscriptstyle N}}(\Phi,\Psi) \,=\, \,\,\mathrm{tr}\; \left[\mathrm{atanh}^{\scriptscriptstyle 2} \left(R^{\scriptscriptstyle 1/2}\right)\right] \hspace{0.5cm} R = (\Psi - \Phi)(I - \Phi^{\scriptscriptstyle *}\Psi)^{\scriptscriptstyle -1}(\Psi^{\scriptscriptstyle *} - \Phi^{\scriptscriptstyle *})(I - \Phi\Psi^{\scriptscriptstyle *})^{\scriptscriptstyle -1} \;\;\; \Phi,\Psi \in \mathbb{D}_{\scriptscriptstyle N}
\end{eqnarray}
Here, $d^{\scriptscriptstyle 2}_{\scriptscriptstyle \mathcal{T}_{\scriptscriptstyle N}}(P^{\scriptscriptstyle (1)},P^{\scriptscriptstyle (2)})$ is the Riemannian distance (\ref{eq:tndistance}) between $P^{\scriptscriptstyle (1)},P^{\scriptscriptstyle (2)} \in \mathcal{T}_{\scriptscriptstyle N}$.
\begin{eqnarray} \label{eq:btnvolume}
  dv(T) \,=\, {\scriptstyle \sqrt{n}\prod^{n-1}_{j=1} (n-j)} \, \times \, dv_{\scriptscriptstyle \mathcal{T}_{\scriptscriptstyle N}}(P) \, \prod^{\scriptscriptstyle n-1}_{j=1} \, dv(\Omega_{\scriptscriptstyle j})  \hspace{2.1cm}\; \\[0.1cm]
   \nonumber \mbox{where } \;\; dv(\Omega) \,=\,  \frac{\prod_{\scriptscriptstyle i\leq j} \mathrm{Re}\,  d \Omega_{ij}\, \mathrm{Im} \,d \Omega_{ij}}{\strut \det\left(I - \Omega \Omega^{\scriptscriptstyle *}\right)^{\scriptscriptstyle N+1}}  \;\;\; \Omega \in \mathbb{D}_{\scriptscriptstyle N} \hspace{2cm}\;
\end{eqnarray}
Here, $dv_{\scriptscriptstyle \mathcal{T}_{\scriptscriptstyle N}}(P)$ is the Riemannian volume element (\ref{eq:tnvolume}) computed at $P \in \mathcal{T}_{\scriptscriptstyle N}$. In expressions (\ref{eq:btndistance}) and (\ref{eq:btnvolume}),  $d^{\scriptscriptstyle\, 2}_{\scriptscriptstyle \mathbb{D}_{\scriptscriptstyle N}}(\Phi,\Psi)$ and $dv(\Omega)$ denote the Riemannian distance and Riemannian volume element on $\mathbb{D}_{\scriptscriptstyle N}$, corresponding to the Siegel metric $ds^{\scriptscriptstyle 2}_{\scriptscriptstyle\Omega}$ appearing in (\ref{eq:btnmetric}). These may be found in~\cite{jeuris1} (Page $10$) and~\cite{hua} (Chapter IV, Page $84$), respectively. \\[0.1cm]
\indent \textbf{-- Symmetric space\,:} to see that the Riemannian metric (\ref{eq:btnmetric}) makes $\mathcal{T}^{\scriptscriptstyle N}_{\scriptscriptstyle n}$ into a Riemannian symmetric space of non-positive curvature, it is enough to notice that decomposition (\ref{eq:bttnssdecomposition}) identifies $\mathcal{T}^{\scriptscriptstyle N}_{\scriptscriptstyle n}$ with a product of Riemannian symmetric spaces of non-positive curvature. The first factor in this decomposition is the space $\mathcal{T}_{\scriptscriptstyle N}$ of Toeplitz covariance matrices. This was seen to be a Riemannian symmetric space of non-positive curvature in the previous Paragraph \ref{subsec:tn}. Each of the remaining factors is a copy of $\mathbb{D}_{\scriptscriptstyle N}$, equipped with a scaled version of the Siegel metric $ds^{\scriptscriptstyle 2}_{\scriptscriptstyle \Omega}$, as in (\ref{eq:btnmetric}). With this metric, $\mathbb{D}_{\scriptscriptstyle N}$ is a Riemannian symmetric space of non-compact type, known as the Siegel disc~\cite{helgason} (Exercice D.1., Page $527$). 

To apply Facts \ref{fact:1}--\ref{fact:5} to $\mathcal{T}^{\scriptscriptstyle N}_{\scriptscriptstyle n}$, consider the group $G \, = \, G_{\scriptscriptstyle \mathcal{T}_{\scriptscriptstyle N}} \, \times \, Sp(n,\mathbb{R}) \, \times \, \ldots \, \times Sp(n,\mathbb{R})$, where $G_{\scriptscriptstyle \mathcal{T}_{\scriptscriptstyle N}} = \mathbb{R}_{+} \, \times \, SU(1,1) \, \times \, \ldots \, \times \, SU(1,1)$, as in (\ref{eq:gsieg})--(\ref{eq:tngroup}) of Paragraph \ref{subsec:tn}, and where $Sp(n,\mathbb{R})$ is the group of $2N \times 2N$ complex matrices $U$ such that~\cite{siegel} (Page $10$),
\begin{equation} \label{eq:gsieg1}
  U \, = \, \left(\begin{array}{cc} A & B \\ C & D \end{array}\right): \hspace{0.35cm}
  U^{\scriptscriptstyle T} \left(\begin{array}{cc} 0 & -I \\ I & 0 \end{array}\right) U \;= \left(\begin{array}{cc} 0 & -I \\ I & 0 \end{array}\right) \hspace{0.25cm}\, , \hspace{0.25cm}  
  U^{\scriptscriptstyle H} \left(\begin{array}{cc} -I & 0 \\ 0 & I \end{array}\right) U \;= \left(\begin{array}{cc} -I & 0 \\ 0 & I \end{array}\right)
\end{equation}
with $A,B,C,D$ complex $N \times N$ matrices. This group $G$ acts transitively and isometrically on $\mathcal{T}^{\scriptscriptstyle N}_{\scriptscriptstyle n}$ in the following way,
\begin{eqnarray} \label{eq:btngroup}
  g \cdot T \,=\,  (S\cdot P\,,U_{\scriptscriptstyle 1}\cdot \Omega_{\scriptscriptstyle 1}, \ldots, U_{\scriptscriptstyle n-1}\cdot \Omega_{\scriptscriptstyle n-1})
  \hspace{0.25cm},\hspace{0.25cm} g = (S,U_{\scriptscriptstyle 1},\ldots, U_{\scriptscriptstyle n-1})
  \\[0.1cm]
   \nonumber \mbox{where } \;\; U \cdot \Omega \,=\, (A\,\Omega + B )(C \,\Omega + D)^{\scriptscriptstyle -1} \hspace{0.2cm} U \in Sp(n,\mathbb{R}) \, , \, \Omega \in \mathbb{D}_{\scriptscriptstyle N} \hspace{0.5cm}
\end{eqnarray}
where $S\cdot  P$ is given by (\ref{eq:tngroup}) for $S \in G_{\scriptscriptstyle \mathcal{T}_{\scriptscriptstyle N}}$ and $P \in  \mathcal{T}_{\scriptscriptstyle N}$. Then, Facts \ref{fact:1} and \ref{fact:2} state that definitions (\ref{eq:btnmetric})--(\ref{eq:btnvolume}) verify the general identities (\ref{eq:invmetric}) and (\ref{eq:invintegral}). Furthermore, Fact \ref{fact:3} states the existence and uniqueness of Riemannian barycentres in $\mathcal{T}^{\scriptscriptstyle N}_{\scriptscriptstyle n}$. 

The application of Facts \ref{fact:4} and \ref{fact:5} to $\mathcal{T}^{\scriptscriptstyle N}_{\scriptscriptstyle n}$ leads to the two following formulae, which are proved in Example \ref{ex:spn} of Appendix \ref{app:roots},
\begin{eqnarray} \label{eq:btnintegralp}
  \int_{\scriptscriptstyle \mathcal{T}^{\scriptscriptstyle N}_{\scriptscriptstyle n}} \, f(T)\, dv(T) \, = \, \mathrm{C} \times \int_{\scriptscriptstyle \mathcal{T}_{\scriptscriptstyle N}} \,
\int_{\scriptscriptstyle \mathbb{D}_{\scriptscriptstyle N}} \ldots
\int_{\scriptscriptstyle \mathbb{D}_{\scriptscriptstyle N}}
\, f(P,\Omega_{\scriptscriptstyle 1},\ldots,\Omega_{\scriptscriptstyle n-1})  \, dv(\Omega_{\scriptscriptstyle 1}) \ldots dv(\Omega_{\scriptscriptstyle n-1}) \,  dv_{\scriptscriptstyle \mathcal{T}_{\scriptscriptstyle N}}(P)  \\[0.1cm]
\nonumber \int_{\scriptscriptstyle \mathbb{D}_{\scriptscriptstyle N}} \, f(\Omega)\, dv(\Omega) \, = \, \mathrm{C} \times \int_{\scriptscriptstyle SU(N)} \int_{\scriptscriptstyle \mathbb{R}^{\scriptscriptstyle N}} \, f(r,\Theta) \, \prod_{\scriptscriptstyle i < j}\sinh(|r_{\scriptscriptstyle i}-r_{\scriptscriptstyle j}|)\prod_{\scriptscriptstyle i \leq j}\sinh(|r_{\scriptscriptstyle i}+r_{\scriptscriptstyle j}|) \, dr \, d\Theta \hspace{0.7cm} 
\end{eqnarray}
\vspace{-0.05cm}
\begin{eqnarray} \label{eq:btnkilling}
 d^{\scriptscriptstyle \,2}(\mathcal{I},T) \, = \, n\, d^{\scriptscriptstyle 2}_{\scriptscriptstyle \mathcal{T}_{\scriptscriptstyle N}}(I,P) \,+\, 
   \sum^{\scriptscriptstyle n-1}_{j=1} \,(n-j)\, d^{\scriptscriptstyle\, 2}_{\scriptscriptstyle \mathbb{D}_{\scriptscriptstyle N}}(I,\Omega_{\scriptscriptstyle j}) \hspace{0.3cm} (\mathcal{I} \, = \, nN \times nN \mbox{ identity matrix})\, \\
\nonumber   d^{\scriptscriptstyle\, 2}_{\scriptscriptstyle \mathbb{D}_{\scriptscriptstyle N}}(I,\Omega) \, = \, 
   \Vert r \Vert^{\scriptscriptstyle 2} \, = \, \sum^{\scriptscriptstyle n}_{\scriptscriptstyle j=1} \, r^{\scriptscriptstyle 2}_{\scriptscriptstyle j} \hspace{7cm}\;\,
\end{eqnarray} 
For $\mathcal{T}^{\scriptscriptstyle N}_{\scriptscriptstyle n}$, formulae (\ref{eq:btnintegralp}) and (\ref{eq:btnkilling}) play the same role as formulae (\ref{eq:integralp}) and (\ref{eq:killing}), for a general symmetric space $\mathcal{M}$. Here, $dv_{\scriptscriptstyle \mathcal{T}}(P)$ and $d^{\scriptscriptstyle 2}_{\scriptscriptstyle \mathcal{T}_{\scriptscriptstyle N}}(I,P)$ are given by (\ref{eq:tnintegralp}) and (\ref{eq:tnkilling}). Moreover,  $\Omega \in \mathbb{D}_{\scriptscriptstyle N}$ is written in the form,
\begin{equation} \label{eq:takagi}
\Omega = \Theta \, \mathrm{tanh}(r) \, \Theta^{\scriptscriptstyle T} \hspace{0.5cm} \mathrm{tanh}(r) = \mathrm{diag}\left(\mathrm{tanh}(r_{\scriptscriptstyle 1}),\ldots, \mathrm{tanh}(r_{\scriptscriptstyle n})\right)
\end{equation}
where $\Theta \in SU(N)$, the group of unitary $N \times N$ matrices which have unit determinant, and $r \in \mathbb{R}^{\scriptscriptstyle N}$ (this is called the Takagi factorisation in~\cite{horn}). In (\ref{eq:btnintegralp}), $dr$ denotes the Lebesgue measure on $\mathbb{R}^{\scriptscriptstyle N}$ and $d\Theta$ denotes the normalised Haar measure on $SU(N)$.  \\[0.1cm]
 \indent \textbf{-- Gaussian distributions\,:} using formulae (\ref{eq:btnintegralp}) and (\ref{eq:btnkilling}), it is possible to construct the functions $\log Z : \mathbb{R}_+ \rightarrow \mathbb{R}$ and $\Phi : \mathbb{R}_+ \rightarrow \mathbb{R}_+$ for Gaussian distributions on $\mathcal{T}^{\scriptscriptstyle N}_{\scriptscriptstyle n}$. Computation of $\log Z$, based on these two formulae, leads to a similar structure to the one found in (\ref{eq:tnz}) of \ref{subsec:tn}, for the case of Gaussian distributions on $\mathcal{T}_{\scriptscriptstyle n\,}$. 
 
To see this, recall according to (\ref{eq:zgeneral}) in the proof of Proposition \ref{prop:z},
$$
Z(\sigma) \, = \, \int_{\scriptscriptstyle \mathcal{T}^{\scriptscriptstyle N}_{\scriptscriptstyle n}} \, \exp \left[ - \frac{d^{\scriptscriptstyle \,2}(\mathcal{I},T)}{2\sigma^{\scriptscriptstyle 2}}\right]\, dv(T) 
$$
Upon replacing (\ref{eq:btnintegralp}) and (\ref{eq:btnkilling}), this becomes,
$$
 Z(\sigma) \,=\, \mathrm{C} \times \int_{\scriptscriptstyle \mathcal{T}_{\scriptscriptstyle N}} \,
\exp\left[ -\frac{n \,d^{\scriptscriptstyle 2}_{\scriptscriptstyle \mathcal{T}_{\scriptscriptstyle N}}(I,P)}{2\sigma^{\scriptscriptstyle 2}}\right] \, dv_{\scriptscriptstyle \mathcal{T}_{\scriptscriptstyle N}}(P) \, \times \, \prod^{\scriptscriptstyle n-1}_{\scriptscriptstyle j=1} \,\,
 \int_{\scriptscriptstyle \mathbb{D}_{\scriptscriptstyle N}} \, \exp\left[ -\frac{(n-j) \,d^{\scriptscriptstyle\, 2}_{\scriptscriptstyle \mathbb{D}_{\scriptscriptstyle N}}(I,\Omega_{\scriptscriptstyle j})}{2\sigma^{\scriptscriptstyle 2}}\right] 
 \, dv(\Omega_{\scriptscriptstyle j})
$$
This formula states that $Z(\sigma)$ is the product of $n$ integrals. The first of these integrals, which extends over $\mathcal{T}_{\scriptscriptstyle N}\,$, is of the form $Z_{\scriptscriptstyle \mathcal{T}_{\scriptscriptstyle N}}\left(\sigma\middle/\sqrt{ n}\,\right)$, where $Z_{\scriptscriptstyle \mathcal{T}_{\scriptscriptstyle N}}(\sigma)$ is the function given by (\ref{eq:tnz0}) of \ref{subsec:tn}. On the other hand, each one of the remaining $n-1$ integrals, which extend over $\mathbb{D}_{\scriptscriptstyle N}$, is of the form $Z_{\scriptscriptstyle \mathbb{D}_{\scriptscriptstyle N}}\left( \sigma \middle/\sqrt{\scriptstyle n-j}\,\right)$, where $Z_{\scriptscriptstyle \mathbb{D}_{\scriptscriptstyle N}}(\sigma)$ is given by
\begin{equation} \label{eq:zdn}
Z_{\scriptscriptstyle \mathbb{D}_{\scriptscriptstyle N}}(\sigma) \, = \, 
 \int_{\scriptscriptstyle \mathbb{D}_{\scriptscriptstyle N}} \, \exp\left[ -\frac{\,d^{\scriptscriptstyle\, 2}_{\scriptscriptstyle \mathbb{D}_{\scriptscriptstyle N}}(I,\Omega)}{2\sigma^{\scriptscriptstyle 2}}\right] 
 \, dv(\Omega) 
\end{equation}
Using these definitions of $Z_{\scriptscriptstyle \mathcal{T}_{\scriptscriptstyle N}}(\sigma)$ and $Z_{\scriptscriptstyle \mathbb{D}_{\scriptscriptstyle N}}(\sigma)$, it follows that,
\begin{equation} \label{eq:btnz0}
 Z(\sigma) \,=\, \mathrm{C} \times Z_{\scriptscriptstyle \mathcal{T}_{\scriptscriptstyle N}}\left(\sigma\middle/\sqrt{\scriptstyle n}\,\right) \, \times \prod^{\scriptscriptstyle n-1}_{\scriptscriptstyle j=1} \, Z_{\scriptscriptstyle \mathbb{D}_{\scriptscriptstyle N}}\left( \sigma \middle/\sqrt{\scriptstyle n-j}\,\right)
\end{equation}
Then, by taking logarithms, the following expression of $\log Z$ is found,
\begin{equation} \label{eq:btnz}
 \log Z(\sigma) \,=\, \mathrm{C} \, + \, \log Z_{\scriptscriptstyle \mathcal{T}_{\scriptscriptstyle N}}\left(\sigma\middle/\sqrt{\scriptstyle n}\,\right) \, + \,  \sum^{\scriptscriptstyle n-1}_{\scriptscriptstyle j=1} \, \log Z_{\scriptscriptstyle \mathbb{D}_{\scriptscriptstyle N}}
 \left( \sigma\middle/\sqrt{\scriptstyle n-j}\,\right)
\end{equation}
Now, this expression indeed has a similar structure to the one in (\ref{eq:tnz}) of \ref{subsec:tn}. Formally, it is obtained from (\ref{eq:tnz}) after two substitutions\,: $\log(\sigma) \, \rightsquigarrow \, \log Z_{\scriptscriptstyle \mathcal{T}_{\scriptscriptstyle N}}(\sigma)$ and $\log Z_{\scriptscriptstyle \mathbb{D}}
 (\sigma) \, \rightsquigarrow \log Z_{\scriptscriptstyle \mathbb{D}_{\scriptscriptstyle N}}(\sigma)$. Referring to decompositions (\ref{eq:tnssdecomposition}) and (\ref{eq:bttnssdecomposition}), it is seen that, on a geometric level, these correspond to replacing the Euclidean space $\mathbb{R}_+$ by the space $\mathcal{T}_{\scriptscriptstyle N}$, and replacing the Poincar\'e disc $\mathbb{D}$ by the Siegel disc $\mathbb{D}_{\scriptscriptstyle N}$. 

The function $\log Z$ lends itself to numerical evaluation from expression (\ref{eq:btnz}). In this expression, the function $\log Z_{\scriptscriptstyle \mathcal{T}_{\scriptscriptstyle N}}$ can be computed analytically, using (\ref{eq:tnz0}) and (\ref{eq:tnz}) of \ref{subsec:tn}. On the other hand, the function
$\log Z_{\scriptscriptstyle \mathbb{D}_{\scriptscriptstyle N}}$ can be evaluated from (\ref{eq:zdn}) using (\ref{eq:btnintegralp}) and (\ref{eq:btnkilling}). Precisely, by replacing the second line of each one of these formulae in (\ref{eq:zdn}), it follows that,
$$
Z_{\scriptscriptstyle \mathbb{D}_{\scriptscriptstyle N}}(\sigma) \, = \,
\mathrm{C} \times \int_{\scriptscriptstyle SU(N)} \int_{\scriptscriptstyle \mathbb{R}^{\scriptscriptstyle N}} \, 
\exp\left[ -\frac{\,\Vert r \Vert^{\scriptscriptstyle 2}}{2\sigma^{\scriptscriptstyle 2}}\right] 
\, \prod_{\scriptscriptstyle i < j}\sinh(|r_{\scriptscriptstyle i}-r_{\scriptscriptstyle j}|)\prod_{\scriptscriptstyle i \leq j}\sinh(|r_{\scriptscriptstyle i}+r_{\scriptscriptstyle j}|) \, dr \, d\Theta 
$$
and since the function under the integral does not depend on $\Theta$, 
\begin{equation} \label{eq:btnzintegral}
Z_{\scriptscriptstyle \mathbb{D}_{\scriptscriptstyle N}}(\sigma) \, = \,
\mathrm{C} \times \int_{\scriptscriptstyle \mathbb{R}^{\scriptscriptstyle N}} \, 
\exp\left[ -\frac{\,\Vert r \Vert^{\scriptscriptstyle 2}}{2\sigma^{\scriptscriptstyle 2}}\right] 
\, \prod_{\scriptscriptstyle i < j}\sinh(|r_{\scriptscriptstyle i}-r_{\scriptscriptstyle j}|)\prod_{\scriptscriptstyle i \leq j}\sinh(|r_{\scriptscriptstyle i}+r_{\scriptscriptstyle j}|) \, dr 
\end{equation}
\begin{figure} 
\centering
\begin{subfigure}{.5\textwidth}
  \centering
  \includegraphics[width=1\linewidth]{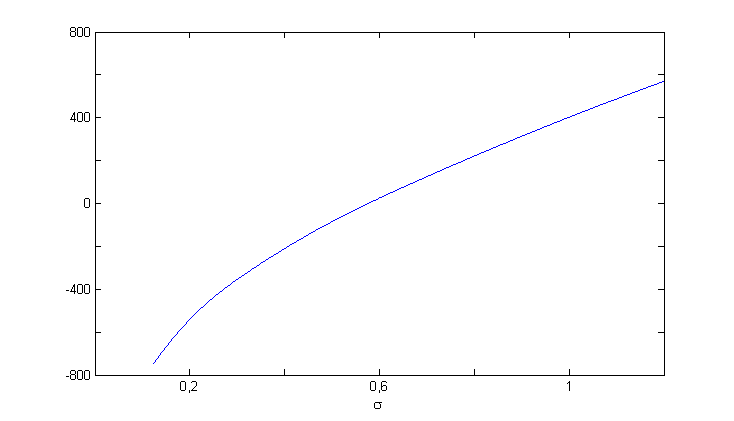}
  \caption{Graph of $\log Z_{\scriptscriptstyle \mathbb{D}_{\scriptscriptstyle N}} : \mathbb{R}_+ \rightarrow \mathbb{R}$}
  \label{fig:hn1}
\end{subfigure}%
\begin{subfigure}{.5\textwidth}
  \centering
  \includegraphics[width=1\linewidth]{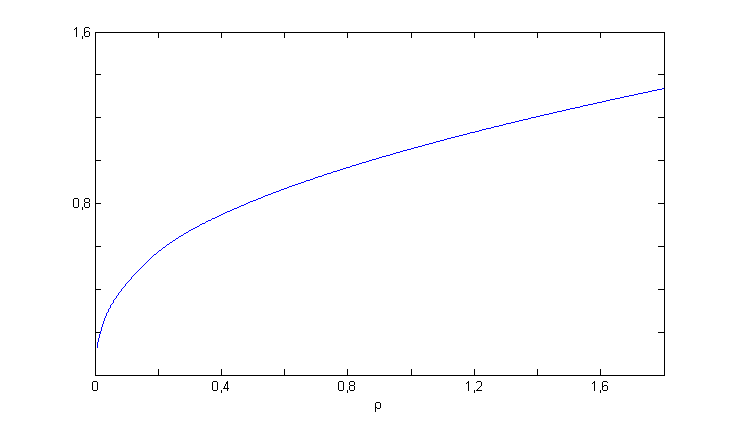}
  \caption{Graph of $\Phi _{\scriptscriptstyle \mathbb{D}_{\scriptscriptstyle N}}:\mathbb{R}_+ \rightarrow \mathbb{R}_+$}
  \label{fig:hn2}
\end{subfigure}
\caption{$\log Z_{\scriptscriptstyle \mathbb{D}_{\scriptscriptstyle N}}$ and $\Phi_{\scriptscriptstyle \mathbb{D}_{\scriptscriptstyle N}}$ for the Siegel disc $\mathbb{D}_{\scriptscriptstyle N}$ with $N = 20$}
\label{fig:btn}
\end{figure}
Using this result, the function $\log Z_{\scriptscriptstyle \mathbb{D}_{\scriptscriptstyle N}}$ can be evaluated by means of the Monte Carlo technique of~\cite{paolo}. Then, complete evaluation of $\log Z$ is achieved by substitution in (\ref{eq:btnz}). 

In addition to being useful in the evaluation of $\log Z$, the function $Z_{\scriptscriptstyle \mathbb{D}_{\scriptscriptstyle N}}$ is interesting in its own right, since it gives the normalising factor for Gaussian distributions on the Siegel disc $\mathbb{D}_{\scriptscriptstyle N}$. Figure \ref{fig:btn} presents the graph of the function $\log Z_{\scriptscriptstyle \mathbb{D}_{\scriptscriptstyle N}}$ and of the corresponding function $\Phi_{\scriptscriptstyle \mathbb{D}_{\scriptscriptstyle N}}$, (defined as in (\ref{eq:phi}) from Proposition \ref{prop:mle} of \ref{subsec:inference}), for $N = 20$.
\end{subequations}
\\[0.15cm]
\indent $\blacktriangleright$ \textbf{Remark 7} : the Siegel disc $\mathbb{D}_{\scriptscriptstyle N}$ is isometric to the Siegel half-plane $\mathbb{H}_{\scriptscriptstyle N}$, the set of $N \times N$ symmetric complex matrices $Z$ such that $\mathrm{Im}\, Z \succ 0$~\cite{siegel}\cite{terras2}. The Siegel half-plane $\mathbb{H}_{\scriptscriptstyle N}$, (in addition to its usefulness in analytical mechanics and number theory), provides a geometric setting which leads to deep results with regard to linear filtering and linear control theory~\cite{zelikin,bougerol1,bougerol2}.  In particular, both the discrete and the continuous forms of the Riccati equation are equivalent to invariant dynamical systems on $\mathbb{H}_{\scriptscriptstyle N}$, which have remarkable geometric properties. This provides a geometric means of studying the asymptotic behavior of the Riccati equation, and therefore of linear filters and control systems.\hfill $\blacksquare$
\section{Density estimation with finite mixtures of Gaussian densities} \label{sec:mixture}
The present section introduces finite mixtures of Gaussian densities, as a tool for dealing with the problem of density estimation in a Riemannian symmetric space of non-positive curvature, $\mathcal{M}$. A finite mixture of Gaussian densities on $\mathcal{M}$ is a probability density of the form
\begin{equation} \label{eq:mixture}
p(x) = \sum^{\scriptscriptstyle K}_{\scriptscriptstyle \kappa \,= 1} \omega_{\scriptscriptstyle \kappa} \times p(x | \,\bar{x}_{\scriptscriptstyle \kappa},\sigma_{\scriptscriptstyle \kappa})  
\end{equation}  
where the weights $\omega_{\scriptscriptstyle 1},\ldots,\omega_{\scriptscriptstyle K} > 0$ are such that $\omega_{\scriptscriptstyle 1} + \ldots + \omega_{\scriptscriptstyle K} = 1$, and where each density $p(x | \,\bar{x}_{\scriptscriptstyle \kappa},\sigma_{\scriptscriptstyle \kappa})$ is a Gaussian density given by (\ref{eq:gaussianpdf}) of Section \ref{sec:gaussian}. Definition (\ref{eq:mixture}) is directly based on the standard definition of a finite mixture of probability densities~\cite{mixtures}\cite{mixture}. 
  
The problem of density estimation in $\mathcal{M}$ is the following~\cite{vapnik}\cite{tibshirani}. Assume a population of data points $x_{\scriptscriptstyle 1},\ldots, x_{\scriptscriptstyle N} \in \mathcal{M}$ is under study. Assume further this population was generated from some unknown probability density $q(x)$. The problem is to approximate this unknown density using the data $x_{\scriptscriptstyle 1},\ldots, x_{\scriptscriptstyle N}$. 

Concretely, this problem arises in the following manner~\cite{rosu,eusipco,gretsi}~\cite{chevallier1,chevallier2,chevallier3}. As in Section \ref{sec:matrices}, the space $\mathcal{M}$ is a space of structured covariance matrices. For example, $\mathcal{M} \,=\, \mathcal{H}_{\scriptscriptstyle n}$, \,$\mathcal{M} \,=\, \mathcal{T}_{\scriptscriptstyle n}$, or\, $\mathcal{M} \,=\, \mathcal{T}^{\scriptscriptstyle N}_{\scriptscriptstyle n}$, (recall these notations from Table 1. in the introduction).  A database of signals or images is under study, and each object in this database is represented by a covariance descriptor. This descriptor is taken to be a structured covariance matrix $x \in \mathcal{M}$. The database of signals or images is thus mapped into a population $x_{\scriptscriptstyle 1},\ldots, x_{\scriptscriptstyle N} \in \mathcal{M}$. In a real-world setting, it is hopeless to attempt any simplifying hypothesis with regard to the probability density $q(x)$ of this population. Rather, one attempts to approximate $q(x)$ within a sufficiently rich class of candidate densities $p(x)$.

The use of mixtures of Gaussian densities as a tool for solving this problem is founded on the idea that the probability density $q(x)$, at least if sufficiently regular, can be approximated to any precision by a density $p(x)$ of the form (\ref{eq:mixture}). Accordingly, consider the mathematical problem of finding a density $p(x)$ of the form (\ref{eq:mixture}) which best approximates the unknown density $q(x)$, in the sens of Kullback-Leibler divergence. Recall this is given by~\cite{cover}, 
\begin{equation} \label{eq:approx1}
\hat{p} = \mathrm{argmin}_{\scriptscriptstyle p} \;D(q\, \Vert \, p) \hspace{0.5cm}   D(q\, \Vert \, p) \ =\, \int_{\scriptscriptstyle \mathcal{M}}\, q(x)\, \log\left[ \frac{q(x)}{p(x)}\right] \, dv(x)
\end{equation}
where $dv$ is the Riemannian volume element of $\mathcal{M}$. Here, the minimum is over all densities $p(x)$ of the form (\ref{eq:mixture}). Following the general idea of ``empirical risk minimisation''~\cite{vapnik} (Page $32$), the above integral with respect to the unknown density $q(x)$ is replaced by an empirical average. The minimisation problem (\ref{eq:approx1}) is then replaced by
\begin{equation} \label{eq:approx2}
\hat{p}_{\scriptscriptstyle N} = \mathrm{argmin}_{\scriptscriptstyle p} \;\frac{1}{N} \, \sum^{\scriptscriptstyle N}_{\scriptscriptstyle n=1} \, 
\log\left[ \frac{q(x_{\scriptscriptstyle n})}{p(x_{\scriptscriptstyle n})}\right] \, = \,
\mathrm{argmin}_{\scriptscriptstyle p} \, - \frac{1}{N} \, \sum^{\scriptscriptstyle N}_{\scriptscriptstyle n=1} \, 
\log\, p(x_{\scriptscriptstyle n}) 
\end{equation}
This new minimisation problem is equivalent to the problem of maximum likelihood estimation of the parameters $\lbrace (\omega_{\scriptscriptstyle \kappa},\bar{x}_{\scriptscriptstyle \kappa},\sigma_{\scriptscriptstyle \kappa}) \, ; \, \kappa = 1, \ldots, K \rbrace$ of the finite mixture of Gaussian densities $p(x)$, under the assumption that the data $x_{\scriptscriptstyle 1},\ldots, x_{\scriptscriptstyle N}$ are independently generated from $p(x)$. 

In~\cite{said}, an original EM (expectation maximisation) algorithm was proposed for solving this problem of maximum likelihood estimation. This algorithm was derived in the special case where $\mathcal{M} = \mathcal{P}_{\scriptscriptstyle n\,}$, the space of $n \times n$ covariance matrices. It was successfully applied to real data, mainly from remote sensing, with the aim of carrying out density estimation as well as classification~\cite{rosu,eusipco,gretsi}. This algorithm generalises directly from $\mathcal{M} = \mathcal{P}_{\scriptscriptstyle n}$ to any $\mathcal{M}$ which is a Riemannian symmetric space of non-positive curvature, and its derivation remains formally the same as in~\cite{said} (Paragraph IV.A., Page $13$). 

In the present general setting, it can be described as follows. The algorithm assumes that $K$, the number of mixture components in (\ref{eq:mixture}), has already  been selected. It iteratively updates an approximation $\hat{\theta} = \lbrace (\hat{\omega}_{\scriptscriptstyle \kappa},\hat{x}_{\scriptscriptstyle \kappa},\hat{\sigma}_{\scriptscriptstyle \kappa})\rbrace$ of the maximum likelihood estimates of the parameters $\lbrace (\omega_{\scriptscriptstyle \kappa},\bar{x}_{\scriptscriptstyle \kappa},\sigma_{\scriptscriptstyle \kappa})\rbrace$. The required update rules (\ref{eq:omegaupdate})--(\ref{eq:sigmaupdate}), as given below, are repeated as long as this leads to a sensible increase in the joint likelihood function of $(\hat{\omega}_{\scriptscriptstyle \kappa},\hat{x}_{\scriptscriptstyle \kappa},\hat{\sigma}_{\scriptscriptstyle \kappa})$. These update rules involve the two quantities $\pi_{\scriptscriptstyle \kappa}(x_{\scriptscriptstyle n},\hat{\theta})$ and $N_{\scriptscriptstyle \kappa}(\hat{\theta})$, defined as follows,
\begin{equation} \label{eq:emdescription1}
   \pi_{\scriptscriptstyle \kappa}(x_{\scriptscriptstyle n},\hat{\theta}) \,\propto\, \hat{\omega}_{\scriptscriptstyle \kappa} \times p(x_{\scriptscriptstyle n}|\, \hat{x}_{\scriptscriptstyle \kappa},\hat{\sigma}_{\scriptscriptstyle \kappa}) \hspace{0.5cm} N_{\scriptscriptstyle \kappa}(\hat{\theta}) = \sum^{\scriptscriptstyle N}_{n=1} \pi_{\scriptscriptstyle \kappa}(x_{\scriptscriptstyle n},\hat{\theta})
\end{equation}
where the constant of proportionality, corresponding to ``$\propto$" in the first expression, is chosen so ${\scriptstyle \sum_{\kappa}}\,\pi_{\scriptscriptstyle \kappa}(x_{\scriptscriptstyle n},\hat{\theta})  \, = \, 1$. Each iteration of the EM algorithm consists in an application of the updatte rules,
 \vspace{0.3cm}
\hrule
\vspace{0.3cm}
\noindent \underline{EM update rules for mixture estimation\,:}\\[0.1cm]
\begin{subequations}
\begin{equation} \label{eq:omegaupdate}
  \hat{\omega}^{\mathrm{new}}_{\scriptscriptstyle \kappa} = \left. N_{\scriptscriptstyle \kappa}(\hat{\theta}) \middle/ N  \right.
\end{equation}
\begin{equation} \label{eq:xbarupdate}
   \hat{x}^{\mathrm{new}}_{\scriptscriptstyle \kappa} = \mathrm{argmin}_{\scriptscriptstyle x \in \mathcal{M}} \, E_{\scriptscriptstyle \kappa}(x) \hspace{0.5cm}
   E_{\scriptscriptstyle \kappa}(x) \ = \, N^{\scriptscriptstyle -1}_\kappa(\hat{\theta}) \, \sum^{\scriptscriptstyle N}_{\scriptscriptstyle n=1} \,  \pi_{\scriptscriptstyle \kappa}(x_{\scriptscriptstyle n},\hat{\theta}) \, d^{\,2}(x,x_{\scriptscriptstyle n})
\end{equation}
\begin{equation} \label{eq:sigmaupdate}
   \hat{\sigma}^{\mathrm{new}}_{\scriptscriptstyle \kappa} = \Phi  (\, E^{\phantom{s}}_{\scriptscriptstyle \kappa}(\hat{x}^{\mathrm{new}}_{\scriptscriptstyle \kappa}))  \hspace{0.5cm} \Phi \mbox{ given by (\ref{eq:phi}) of Proposition \ref{prop:mle}} 
\end{equation}
\end{subequations}
 \vspace{0.3cm}
\hrule
\vspace{0.3cm}

Realisation of update rules (\ref{eq:omegaupdate}) and (\ref{eq:sigmaupdate}) is straightforward. On the other hand, realisation of update rule (\ref{eq:xbarupdate}) requires minimisation of the function $E_{\scriptscriptstyle \kappa}(x)$ over $x \in \mathcal{M}$. It is clear this function is a variance function of the form (\ref{eq:variance}), corresponding to the probability distribution ${\scriptstyle \pi \, \propto \,\, \sum^N_{n=1}\, \pi_\kappa(x_n,\hat{\theta}) \times \delta_{x_n}}$ where $\delta_{\scriptscriptstyle x}$ denotes the Dirac distribution concentrated at $x \in \mathcal{M}$. Therefore, existence and uniqueness of $\hat{x}^{\mathrm{new}}_{\scriptscriptstyle \kappa}$ follow from Fact \ref{fact:3}, and computation of $\hat{x}^{\mathrm{new}}_{\scriptscriptstyle \kappa}$ can be carried out using existing routines for the computation of Riemannian barycentres\!~\cite{lenglet}\cite{ferreira}\cite{aistats}\cite{bonnabel}. \\[0.1cm]
\indent The EM algorithm just described requires prior selection of the number $K$ of mixture components. When this EM algorithm was applied in~\cite{rosu,eusipco,gretsi}, (for the special case $\mathcal{M} = \mathcal{P}_{\scriptscriptstyle n}$)\,, $K$ was selected by maximising a BIC criterion. In general, this criterion is given following~\cite{leroux},
\begin{equation} \label{eq:bic}
\mbox{BIC criterion\,:} \hspace{0.3cm} K^{\scriptscriptstyle *} \, = \, \mathrm{argmax}_{\scriptscriptstyle K}\, BIC(K) \hspace{0.5cm}
BIC(K) \, = \, \ell(K) \,-\, \frac{1}{2} \, DF \times \ln(N)
\end{equation}
Here, $\ell(K)$ is the maximum value of the joint likelihood function of $(\omega_{\scriptscriptstyle \kappa},\bar{x}_{\scriptscriptstyle \kappa},\sigma_{\scriptscriptstyle \kappa})$ for given $K$, and $DF$ is the number of parameters to be estimated, or number of degrees of freedome. This is $DF = K\times [\, 2 + \mathrm{dim}(\mathcal{M})\,] - 1$. For the three spaces of structured covariance matrices considered in Section \ref{sec:matrices}, the dimension of $\mathcal{M}$ is given by, 
\begin{equation} \label{eq:dimension}
\mathrm{dim}(\mathcal{H}_{\scriptscriptstyle n}) \,=\, n^{\scriptscriptstyle 2} \hspace{0.3cm} \mathrm{dim}(\mathcal{T}_{\scriptscriptstyle n}) \,=\, 2n - 1 \hspace{0.3cm} \mathrm{dim}(\mathcal{T}^{\scriptscriptstyle N}_{\scriptscriptstyle n}) \,=\, 2N-1 + (n-1)\times(N^{\scriptscriptstyle 2}+N) 
\end{equation}
The EM algorithm given by update rules (\ref{eq:omegaupdate})--(\ref{eq:sigmaupdate}), and the BIC criterion (\ref{eq:bic}), together provide a complete programme for solving the problem of density estimation in $\mathcal{M}$, using finite mixtures of Gaussian densities. The output of this programme is an approximation of the finite mixture of Gaussian densities $\hat{p}_{\scriptscriptstyle N}$ defined by (\ref{eq:approx2}). Then, with a slight abuse of notation, it is possible to write,
\begin{equation} \label{eq:estimatee}
\hat{p}_{\scriptscriptstyle N}(x) = \sum^{\scriptscriptstyle K}_{\scriptscriptstyle \kappa \,= 1} \hat{\omega}_{\scriptscriptstyle \kappa} \times p(x | \,\hat{x}_{\scriptscriptstyle \kappa},\hat{\sigma}_{\scriptscriptstyle \kappa})  
\end{equation}
In~\cite{said}\cite{rosu,eusipco,gretsi}, an additional problem of classification was addressed, whose aim is to use (\ref{eq:estimatee}) as an explicative model of the data $x_{\scriptscriptstyle 1},\ldots, x_{\scriptscriptstyle N}$, revealing the fact that these data arise from $K$ disjoint classes, where the $\kappa\mathrm{th}$ class is generated from a single Gaussian distribution $G(\hat{x}_{\scriptscriptstyle \kappa},\hat{\sigma}_{\scriptscriptstyle \kappa})$.
In~\cite{said} (Paragraph IV.B., Page $15$), this problem was solved in the special case $\mathcal{M} = \mathcal{P}_{\scriptscriptstyle n\,}$. In general, the proposed solution remains valid as long as $\mathcal{M}$ is a Riemannian symmetric space of non-positive curvature. It is based on a Bayes optimal classification rule. This is a function $\kappa^{\scriptscriptstyle *} : \mathcal{M} \rightarrow \lbrace 1,\ldots, K\rbrace\,$, obtained by minimising the \textit{a posteriori} probability of classification error~\cite{tibshirani} (Section 2.4., Page $21$),
\begin{equation} \label{eq:proberror}
  \kappa^{\scriptscriptstyle *}(x) \, = \, \mathrm{argmin}_{\scriptscriptstyle \kappa = 1,\ldots,K}\, \left[ \,1 - P(\kappa|x)\,\right]
\end{equation}
where $P(\kappa|x)$ is found from Bayes rule, $P(\kappa|x) \propto p(x|\kappa) \times P(\kappa)$, after replacing $p(x|\kappa) = p(x | \,\hat{x}_{\scriptscriptstyle \kappa},\hat{\sigma}_{\scriptscriptstyle \kappa})$ and $P(\kappa) = \hat{\omega}_{\scriptscriptstyle \kappa\,}$, where $p(x | \,\hat{x}_{\scriptscriptstyle \kappa},\hat{\sigma}_{\scriptscriptstyle \kappa})$ is given by (\ref{eq:gaussianpdf}). Carrying out the minimisation (\ref{eq:proberror}) readily gives the expression,
\begin{equation} \label{eq:classrule}
   \kappa^{\scriptscriptstyle *}(x) \, = \, \mathrm{argmin}_{\scriptscriptstyle \kappa = 1,\ldots,K}\, \left [ - \log \, \hat{\omega}_{\scriptscriptstyle \kappa} + \log \, Z(\hat{\sigma}_{\scriptscriptstyle \kappa}) + \frac{d^{\,2}(x,\hat{x}_{\scriptscriptstyle \kappa})}{\strut 2\hat{\sigma}^2_{\scriptscriptstyle \kappa}}  \right ]
\end{equation}
This is a generalisation of the so-called minimum-distance-to-mean classification rule, which corresponds to all $\hat{\omega}_{\scriptscriptstyle \kappa}$ and all $\hat{\sigma}_{\scriptscriptstyle \kappa}$ being constant, (in the sense that they do not depend on $\kappa$). In~\cite{congedo}, the minimum-distance-to-mean classification rule was shown to give very good performance in the area brain-computer interface decoding. Ongoing work, started in~\cite{paolo}, examines the application of  (\ref{eq:classrule}) in this same area. 
\subsection*{\textbf{Comparison to kernel density estimation}} 
In recent literature~\cite{chevallier1,chevallier2,chevallier3}, a fully non-parametric method for density estimation in $\mathcal{M}$ was developed. This is based on the method of kernel density estimation on Riemannian manifolds, first outlined in~\cite{pelletier}. This non-parametric method is alternative to the method of the present section, based on finite mixtures of Gaussian distributions. Currently, no practical comparison of these two methods is available, as they have not yet been tested side-by-side on real or simulated data. However, it is interesting to compare them from a theoretical point of view. 

\begin{subequations}
The main idea of kernel density estimation on Riemannian manifolds is the following~\cite{pelletier}\,: classical performance bounds for kernel density estimation on a Euclidean space can be immediately transposed to any Riemannian manifold $\mathcal{M}$, provided the kernel function $K:\mathbb{R} \rightarrow \mathbb{R}_+$ is adjusted to the Riemannian volume element of $\mathcal{M}$.  

In the present context, where $\mathcal{M}$ is a Riemannian symmetric space of non-positive curvature, consider the case of a Gaussian kernel $K(x) = (2\pi)^{-\scriptscriptstyle 1/2}\exp(-x^{\scriptscriptstyle 2}/2)$. When this Gaussian kernel is adjusted to the Riemannian volume element of $\mathcal{M}$, the result is a family of Riemannian Gaussian kernels $K(x|\,\bar{x},\sigma)$ on $\mathcal{M}$, which provides the following method.
 \vspace{0.3cm}
\hrule
\vspace{0.3cm}
\noindent \underline{Kernel density estimation on $\mathcal{M}$\,:}\\[0.1cm]
\begin{equation} \label{eq:kernel}
  K(x|\,\bar{x},\sigma) \, = \,  (2\pi\sigma^{\scriptscriptstyle 2})^{\scriptscriptstyle -\mathrm{dim}(\mathcal{M})/2} \, \times \, J(x,\bar{x}) \, \times \, \exp \left[-\,\frac{d^{\,\scriptscriptstyle 2}(x,\bar{x})}{2\sigma^{\,\scriptscriptstyle 2}}\right]
\end{equation}
where the ``adjusting factor'' $J(x,\bar{x})$ is expressed in~\cite{chevallier1} (Page $11$),
The kernel density estimate $\hat{q}_{\scriptscriptstyle N}(x)$ of the unknown density $q(x)$ of the data $x_{\scriptscriptstyle 1},\ldots, x_{\scriptscriptstyle N}$ is then defined by~\cite{chevallier1} (Formula (7), Page $8$)
\begin{equation} \label{eq:kernelestimate}
  \hat{q}_{\scriptscriptstyle N}(x) \, = \, \frac{1}{N} \, \sum^{\scriptscriptstyle N}_{\scriptscriptstyle n=1} \, K(x|\,x_{\scriptscriptstyle n},\sigma)
\end{equation}
\end{subequations}
and the integrated mean square error between $\hat{q}_{\scriptscriptstyle N}(x)$ and $q(x)$ can be made to decrease as a negative power of $N$, by a suitable choice of $\sigma = \sigma(N)$, under the condition that $q(x)$ is twice differentiable~\cite{pelletier} (Theorem 3.1., Page $302$).
 \vspace{0.3cm}
\hrule
\vspace{0.3cm}
The EM algorithm (\ref{eq:omegaupdate})--(\ref{eq:sigmaupdate}) and the kernel density estimation method (\ref{eq:kernel})--(\ref{eq:kernelestimate}) provide two approaches to the problem of density estimation in $\mathcal{M}$, which differ in the following ways,\\[0.1cm]
A--\textit{Cost in memory\,:} for the EM algorithm, evaluation of the output $\hat{p}_{\scriptscriptstyle N}$ from (\ref{eq:estimatee}) requires bounded memory, essentially independent of $N$. For kernel density estimation, evaluation of $\hat{q}_{\scriptscriptstyle N}$ from (\ref{eq:kernelestimate}) requires storing $N$ data points, each of dimension $\mathrm{dim}(\mathcal{M})$, which can be found in (\ref{eq:dimension}). \\[0.1cm]
B--\textit{Computational cost\,:} each iteration of the EM algorithm involves computing $K$ Riemannian barycentres of $N$ data points. For kernel density estimation, computing $\hat{q}_{\scriptscriptstyle N}(x)$ at a single point $x \in \mathcal{M}$ requires computing $N$ factors $J(x,x_{\scriptscriptstyle n})$ and $N$ Riemannian distances $d(x,x_{\scriptscriptstyle n})$. Concretely, Riemannian distances on $\mathcal{M}$ are found from highly non-linear formulae, such as (\ref{eq:hndistance}), (\ref{eq:tndistance}) and (\ref{eq:btndistance}).  \\[0.1cm]
C--\textit{Rate of convergence\,:} EM algorithms are known to be slow, and to get trapped in local minima or saddle points. However, stochastic versions of EM, such as the SEM algorithm~\cite{sem}, effectively remedy these problems. For kernel density estimation, (\ref{eq:kernelestimate}) is an exact formula which can be computed directly. However, it should be kept in mind, the precision obtained using this formula suffers from the curse of dimensionality.  \\[0.1cm]
D--\textit{Free parameters\,:} running the EM algorithm requires fixing the mixture order $K$. The BIC criterion (\ref{eq:bic}) proposed to this effect may be difficult to evaluate. Kernel density estimation requires choosing the kernel width $\sigma = \sigma(N)$. This also seems to be a difficult task, especially when the dimension of $\mathcal{M}$ is high. \\[0.1cm]
\indent In conclusion, the EM algorithm (\ref{eq:omegaupdate})--(\ref{eq:sigmaupdate}) provides an approach to the problem of density estimation in $\mathcal{M}$, which can be expected to offer a suitable rate of convergence and which is not greedy in terms of memory. The main computational requirement of this algorithm is the ability to find Riemannian barycentres, a task for which there exists an increasing number of high-performance routines~\cite{ferreira}\cite{aistats}\cite{bonnabel}\cite{jeurissurvey}\cite{congedo1}\cite{congedo2}. The fact that the EM algorithm reduces the problem of probability density estimation in $\mathcal{M}$ to one of repeated computation of Riemannian barycentres is due to the unique connection which exists between Gaussian distributions in $\mathcal{M}$ and the concept of Riemannian barycentre. This connection is stated in the defining Property (\ref{eq:intro1}) of these Gaussian distributions, and was a guiding motivation for the present paper. 
\\[0.15cm]
\indent $\blacktriangleright$ \textbf{Remark 8} : a Riemannian Gaussian kernel $K(x|\,\bar{x},\sigma)$ is a probability density with respect to the Riemannian volume element of $\mathcal{M}$. That is, it is positive and its integral with respect to this volume element is equal to one. Precisely, $K(x|\,\bar{x},\sigma)$ is the probability density of a random point $x$ in $\mathcal{M}$, given by $x = \mathrm{Exp}_{\scriptscriptstyle \bar{x}}(v)$ where $v$ is a tangent vector to $\mathcal{M}$ at $\bar{x}$ with distribution $v \sim \mathcal{N}(0,\sigma)$, and where $\mathrm{Exp}_{\scriptscriptstyle \bar{x}}$ is the Riemannian exponential mapping~\cite{berger}. This probability density $K(x|\,\bar{x},\sigma)$ is not a Gaussian density and therefore is not compatible with Property (\ref{eq:intro1}). \hfill $\blacksquare$ \\[0.1cm]
\section*{Acknowledgement}
This research was in part funded by the NSF grant \texttt{IIS-1525431} to B.C. Vemuri.

\bibliographystyle{IEEEtran}    
\bibliography{saidetalGAUSS}

\vfill
\pagebreak

\appendices

\section{The Cartan decomposition, Killing form, and roots} \label{app:roots}
This appendix provides additional details on the mathematical facts, concerning Riemannian symmetric spaces, which were introduced in Section \ref{sec:geo}. It also justifies the application of these facts to spaces of structured covariance matrices, which was made in Section \ref{sec:matrices}. 

In Section \ref{sec:geo}, Facts \ref{fact:4} and \ref{fact:5} introduced roots in formula (\ref{eq:integralp}) and the Killing form in formula (\ref{eq:killing}), in order to express i) integrals of functions over $\mathcal{M}$, and ii) distances between pairs of points in $\mathcal{M}$, where $\mathcal{M}$ is a Riemannian symmetric space of non-compact type. In Section \ref{sec:matrices}, formulae (\ref{eq:integralp}) and (\ref{eq:killing}) were specified to spaces of structured covariance matrices, according to the following scheme, 
\vspace{0.1cm}

\begin{center}
\begin{tabular}{cll} 
\underline{Section \ref{sec:geo}} & \multicolumn{2}{c}{\hspace{1cm} \underline{\hspace{0.1cm} Section \ref{sec:matrices}\hspace{0.1cm}}}   \\[0.27cm]
 \multirow{3}{*}{\parbox{7mm}{$\phantom{\,}$Facts \\ \ref{fact:4} \& \ref{fact:5}}} 
                               & \ref{subsec:2t2} : Complex covariance matrices                & \hspace{0.2cm} formulae (\ref{eq:hnintegralp}) and (\ref{eq:hnkilling}) \\[0.1cm]
                                & \ref{subsec:tn} : Toeplitz covariance matrices                   & \hspace{0.2cm} formulae (\ref{eq:tnintegralp}) and (\ref{eq:tnkilling})\\[0.1cm]
                                & \ref{subsec:btn} : Block-Toeplitz covariance matrices                   & \hspace{0.2cm} formulae (\ref{eq:btnintegralp}) and (\ref{eq:btnkilling}) \\[0.1cm]
\end{tabular}
\vspace{0.1cm}

\end{center}
Here, the definitions of roots and of the Killing form are briefly recalled, and the specific formulae of Section \ref{sec:matrices} are accordingly proved. The starting point will be the Cartan decomposition. \\[0.1cm]
\indent \textbf{{\sc The Cartan decomposition and Killing form }\,:} A Riemannian symmetric space of non-compact type $\mathcal{M}$ may always be expressed in the form $\mathcal{M} \, = \left. G \middle/\, H\right.$, where $G$ is a semisimple Lie group of non-compact type and $H$ a compact Lie subgroup of $G$~\cite{helgason} (Chapter VI). Let $\mathfrak{g}$ and $\mathfrak{h}$ denote the Lie algebras of $G$ and $H$, respectively. The local geometric properties of $\mathcal{M}$ are equivalent to algebraic properties of the pair $(\mathfrak{g},\mathfrak{h})$. This equivalence is expressed in terms of the Cartan decomposition, the Killing form, and roots. The Cartan decomposition is the following~\cite{helgason} (Chapter III, Page $156$),
\begin{equation} \label{eq:cartan}
 \mbox{Cartan decomposition\,:} \hspace{0.5cm}  \mathfrak{g} \, = \, \mathfrak{h} \, + \, \mathfrak{p}
\end{equation}
where $\mathfrak{p}$ is a subspace of $\mathfrak{g}\,$, which is complementary to $\mathfrak{h}\,$, and such that\,: $[\mathfrak{h},\mathfrak{p}] \subset \mathfrak{p}$ and $[\mathfrak{p},\mathfrak{p}] \subset \mathfrak{h}$. For the geometric meaning of these conditions, see~\cite{berger} (Chapter 4, Page $186$). In particular, for any point $o \in \mathcal{M}$ there exists a Cartan decomposition (\ref{eq:cartan}) such that $T_{\scriptscriptstyle o\,} \mathcal{M}$ is identified with $\mathfrak{p}$, and the isotropy algebra of $T_{\scriptscriptstyle o\,} \mathcal{M}$ is identified with $\mathfrak{h}$. Here, $T_{\scriptscriptstyle o\,} \mathcal{M}$ denotes the tangent space to $\mathcal{M}$ at $o$. The Killing form $B : \mathfrak{g} \times \mathfrak{g} \rightarrow \mathbb{R}$ is an $\mathrm{Ad}$-invariant symmetric bilinear form on $\mathfrak{g}$, given by
\begin{equation} \label{eq:cartankilling}
\mbox{Killing form\,:}\hspace{0.5cm} B(v,w)\, = \, \mathrm{tr}\left(\mathrm{ad}(v) \, \mathrm{ad}(w) \right) \hspace{0.75cm} v,w \in \mathfrak{g}
\end{equation}
where $\mathrm{ad}(v) : \mathfrak{g} \rightarrow \mathfrak{g}$ is the linear mapping $\mathrm{ad}(v)u = [v,u]$. The Killing form is negative definite on $\mathfrak{h}$ and positive definite on $\mathfrak{p}$~\cite{helgason} (Proposition 7.4., Page $158$). 

The Cartan decomposition and Killing form provide an algebraic expression for the local Riemannian geometry of $\mathcal{M}$. Precisely, when $T_{\scriptscriptstyle o\,}\mathcal{M}$ is identified with $\mathfrak{p}$, the Riemannian metric and the Riemannian curvature tensor of $\mathcal{M}$ can be written as follows~\cite{helgason} (Theorem 4.2., Page $180$),
\begin{subequations} 
\begin{eqnarray}
\label{eq:cartanmetric} \mbox{Riemanniain metric\,:}\hspace{0.5cm}&  \Vert v \Vert^{\scriptscriptstyle 2}_{\scriptscriptstyle o} = B(v,v)& \hspace{0.75cm} v \in T_{\scriptscriptstyle o\,}\mathcal{M} \cong \mathfrak{p}
\\[0.1cm]
 \label{eq:curvature} \mbox{Riemanniain curvature\,:}\hspace{0.5cm}&   R_{\scriptscriptstyle o}(v,w)u \, = \, -\, [\,[v,w]\,,\,u\,]& \hspace{0.75cm} v,w,u \in T_{\scriptscriptstyle o\,}\mathcal{M} \cong \mathfrak{p}
\end{eqnarray}
\end{subequations}
 where $\Vert v \Vert_{\scriptscriptstyle o}$ denotes the Riemannian length of $v \in T_{\scriptscriptstyle o\,}\mathcal{M}$ and $R_{\scriptscriptstyle o}$ the Riemannian curvature tensor at $o \in \mathcal{M}$. \\[0.1cm]
\indent \textbf{{\sc Roots }\,:} It is clear that application of expressions (\ref{eq:cartanmetric}) and (\ref{eq:curvature}) requires an efficient means of computing the linear mapping $\mathrm{ad}(v)$ for $v \in \mathfrak{p}$. This is nicely done using the roots of the Lie algebra $\mathfrak{g}$. A root $\lambda$ of $\mathfrak{g}$ is a linear mapping $\lambda : \mathfrak{a} \rightarrow \mathbb{R}$, where $\mathfrak{a}$ is a maximal Abelian subspace of   $\mathfrak{p}$ (Abelian means commutative, so $[a_{\scriptscriptstyle 1},a_{\scriptscriptstyle 2}] = 0$ for $a_{\scriptscriptstyle 1},a_{\scriptscriptstyle 2} \in \mathfrak{a}$). Precisely, $\lambda : \mathfrak{a} \rightarrow \mathbb{R}$ is called a root of $\mathfrak{g}$ if there exists a non-trivial subspace $\mathfrak{g}_{\scriptscriptstyle \lambda} \subset \mathfrak{g}$ such that $\mathrm{ad}(a)w = \lambda(a)w$ for all $w \in \mathfrak{g}_{\scriptscriptstyle \lambda}$, see~\cite{helgason} (Chapter III). 

First, roots provide an efficient means of computing $\mathrm{ad}(a)$ when $a \in \mathfrak{a}$. Note that the subspace $\mathfrak{g}_{\scriptscriptstyle \lambda}$ is called the root space corresponding to $\lambda$. Its dimension is denoted  $m_{\scriptscriptstyle \lambda}$, (as in (\ref{eq:D}) of Section \ref{sec:geo}). The various subspaces $\mathfrak{g}_{\scriptscriptstyle \lambda}$ are in direct sum, so that the linear mappings $\mathrm{ad}(a)$ for $a \in \mathfrak{a}$ are jointly diagonalisable, and equal to $\lambda(a)$ times the identity mapping on each subspace $\mathfrak{g}_{\scriptscriptstyle \lambda}$. Second, any $v \in \mathfrak{p}$ is of the form $v \,=\, \mathrm{Ad}(h)\,a$ for some $a \in \mathfrak{a}$ and $h \in H$. Then, $\mathrm{ad}(v) = \mathrm{A}d(h) \circ \mathrm{ad}(a) \circ \mathrm{Ad}(h^{\scriptscriptstyle -1\,})$, so that the linear mappings $\mathrm{ad}(v)$ and $\mathrm{ad}(a)$ are the same up to similarity. 

Note finally that the system of roots $\lambda$ of $\mathfrak{g}$ possesses a distinct symmetry, which is described by a discrete group of permutations, called the Weyl group. Here, it is only mentioned that if $\lambda$ is a root, then so is $-\lambda$. Moreover, the set of non-zero roots can be split into two disjoint subsets, called positive roots and negative roots, in such a way that the positive roots span the dual space of $\mathfrak{a}$~\cite{helgason} (Chapter VI). To indicate that a root $\lambda$ is positive, one writes $\lambda > 0$, (again, as in (\ref{eq:D}) of Section \ref{sec:geo}).
\begin{example} \label{ex:sln} here, it is shown that formulae (\ref{eq:hnintegralp}) and (\ref{eq:hnkilling}) follow from Facts \ref{fact:4} and \ref{fact:5}, applied to the case where $\mathcal{M} = \mathcal{H}_{\scriptscriptstyle n}$ as in \ref{subsec:2t2}. Precisely, (\ref{eq:hnintegralp}) and (\ref{eq:hnkilling}) can be obtained from Facts \ref{fact:4} and \ref{fact:5}, by substitution from the following table \cite{helgason}\cite{terras2}.
\vspace{0.2cm}

\begin{center}
\begin{tabular}{lll} 
\multicolumn{3}{c}{\underline{Table 2.\,: Data for the symmetric space $\mathcal{H}_{\scriptscriptstyle n}$}}\\[0.4cm]
  Group action\,:  & $G = \lbrace n \times n \mbox{ complex invertible matrices} \rbrace$ &  action given by (\ref{eq:congruence})  \\[0.14cm]
   Quotient space\,: & $\mathcal{H}_{\scriptscriptstyle n} \, = \, G/H$  & $H =  \lbrace n \times n \mbox{ unitary matrices} \rbrace$ \\[0.14cm]                         
 \multirow{2}{*}{\parbox{15mm}{Lie algebra of $G$\,:}} & \\
 & $\mathfrak{g} = \lbrace n \times n \mbox{ complex matrices} \rbrace$ &  \\[0.14cm] 
\multirow{2}{*}{\parbox{15mm}{Lie algebra of $H$\,:}}                          &  \\
& $\mathfrak{h} =  \lbrace n \times n \mbox{ skew-Hermitian matrices}\rbrace$ &  \\[0.14cm]
\multirow{2}{*}{\parbox{18mm}{Cartan decomposition (\ref{eq:cartan})\,:}}  & \\
& $\mathfrak{p} = \lbrace n \times n \mbox{ Hermitian matrices} \rbrace$ & $T_{\scriptscriptstyle I\,}\mathcal{H}_{\scriptscriptstyle n}\, \cong \,\mathfrak{p} $, ($I = n \times n$ identity matrix)\\[0.14cm]
\multirow{2}{*}{\parbox{15mm}{Abelian subspace $\mathfrak{a}$\,:}} & \\
& $\mathfrak{a} = \lbrace n \times n \mbox{ real diagonal matrices}\rbrace$ &
$r \in \mathfrak{a}$ is written $r = \mathrm{diag}(r_{\scriptscriptstyle 1},\ldots,r_{\scriptscriptstyle n})$   \\[0.14cm]
Positive roots \,: &  $\lambda_{\scriptscriptstyle ij}(r) = r_{\scriptscriptstyle i} - r_{\scriptscriptstyle j}$ \,;\, $i < j$ & with $ m_{\scriptscriptstyle \lambda_{\scriptscriptstyle ij}}= 2$\\[0.14cm]
Killing form \,: & $B(r,r) \, = {\scriptstyle \frac{1}{2n}\, \sum_{\scriptscriptstyle ij}\, (r_{\scriptscriptstyle i} - r_{\scriptscriptstyle j})^2}$ &  for $r \in \mathfrak{a}$
\end{tabular}
\vspace{0.25cm}
\end{center}
The data in this table can be substituded into (\ref{eq:integralp}) and (\ref{eq:killing}), to obtain (\ref{eq:hnintegralp}) and (\ref{eq:hnkilling}). While (\ref{eq:hnintegralp}) follows immediately, additional care is required for (\ref{eq:hnkilling}). This is because $\mathcal{H}_{\scriptscriptstyle n}$ is not a Riemannian symmetric space of non-compact type, but only a Riemannian symmetric space of non-positive curvature~\cite{terras2} (Page $262$). Therefore, (\ref{eq:hnkilling}) does not follow directly from (\ref{eq:killing}), but rather from (\ref{eq:killing}) and (\ref{eq:proddistance}), after noticing $\mathcal{H}_{\scriptscriptstyle n} = \mathbb{R} \, \times \, S\mathcal{H}_{\scriptscriptstyle n\,}$, where $S\mathcal{H}_{\scriptscriptstyle n}$ is the subset of $Y \in \mathcal{H}_{\scriptscriptstyle n}$ such that $\det(Y) = 1$. This implies that $r \in \mathfrak{a}$ should be considered as a couple $r = (\mathrm{tr}(r),\bar{r})$ where $\bar{r} = r - \frac{1}{n}\mathrm{tr}(r)\,I \in \mathfrak{a}$. Then, using the notation of  (\ref{eq:cartanmetric}), it follows from (\ref{eq:proddistance}) that $\Vert r \Vert^{\scriptscriptstyle 2}_{\scriptscriptstyle I} \, = \, \frac{1}{n}\mathrm{tr}^{\scriptscriptstyle 2}(r)+ \, B(\bar{r},\bar{r})$, which is the same as the right hand side of (\ref{eq:hnkilling}). \hfill $\blacksquare$ 
\end{example}
\vspace{0.2cm}

\begin{example} \label{ex:sp1} the space $\mathcal{T}_{\scriptscriptstyle n}$ of Toeplitz covariance matrices, treated in \ref{subsec:tn}, is a special case of the space $\mathcal{T}^{\scriptscriptstyle N}_{\scriptscriptstyle n}$ of block-Toeplitz covariance matrices, treated in \ref{subsec:btn}. Precisely, this is the special case obtained when $N = 1$. Therefore, in order to  show that formulae (\ref{eq:tnintegralp}) and (\ref{eq:tnkilling}) follow from Facts \ref{fact:4} and \ref{fact:5}, applied to the case where $\mathcal{M} = \mathcal{T}_{\scriptscriptstyle n\,}$, it is enough to replace $N = 1$ in the argument of the following example. \hfill $\blacksquare$ 
\end{example}
\vspace{0.2cm}

\begin{example} \label{ex:spn} here, it is shown that formulae (\ref{eq:btnintegralp}) and (\ref{eq:btnkilling}) follow from Facts \ref{fact:4} and \ref{fact:5}, applied to the case where $\mathcal{M} \, = \, \mathcal{T}^{\scriptscriptstyle N}_{\scriptscriptstyle n}$ as in \ref{subsec:btn}. Recall that $\mathcal{T}^{\scriptscriptstyle N}_{\scriptscriptstyle n}$ admits the decomposition (\ref{eq:bttnssdecomposition}). Through (\ref{eq:btndistance}) and (\ref{eq:btnvolume}), this decomposition implies the first line in each of the two formulae (\ref{eq:btnintegralp}) and (\ref{eq:btnkilling}). The second line in each of these two formulae can be found from Facts \ref{fact:4} and \ref{fact:5}, by substitution from the following table~\cite{helgason}\cite{terras2}. In this table, $U$, $J$, $C$, $w$, and $v$, denote $2N \times 2N$ complex matrices, where
$$
U \, = \, \left(\begin{array}{cc} A & B \\ C & D \end{array}\right) \hspace{0.25cm} , \hspace{0.25cm}
J \, = \, \left(\begin{array}{cc} 0 & -I \\ I & 0 \end{array}\right)   \hspace{0.25cm} , \hspace{0.25cm}
C \, = \, \left(\begin{array}{cc} -I & 0 \\ 0 & I \end{array}\right) \hspace{0.25cm}  \hspace{0.25cm} (I = N \times N \mbox{ identity matrix})
$$
\vspace{0.2cm}

\begin{center}
\begin{tabular}{lll} 
\multicolumn{3}{c}{\underline{Table 3.\,: Data for the Siegel disc $\mathbb{D}_{\scriptscriptstyle N}$}}\\[0.4cm]
  Group action\,:  & $G = \lbrace U\,|\,U^{\scriptscriptstyle T} J \,U = J \, , \, U^{\scriptscriptstyle H} C \, U = C \rbrace$ &  action given by (\ref{eq:btngroup})  \\[0.2cm]
   Quotient space\,: & $\mathbb{D}_{\scriptscriptstyle N} \, = \, G/H$  & $H =  \left\lbrace U\,|\, A = D^{\scriptscriptstyle H\,}\, , \, A \mbox{ unitary} \right\rbrace$ \\[0.14cm]                         
 \multirow{2}{*}{\parbox{15mm}{Lie algebra of $G$\,:}} & \\
 & \multicolumn{2}{l}{$\mathfrak{g} = \lbrace w\,|\, w^{\scriptscriptstyle T}J + J\,w = 0 \, , \, 
 w^{\scriptscriptstyle H}C + C\,w = 0\rbrace$}  \\[0.4cm] 
\multirow{3}{*}{\parbox{15mm}{Lie algebra of $H$\,:}}                          &  \\
& $\mathfrak{h} =  
\left\lbrace w = 
\left(\begin{array}{cc} \omega & 0 \\[0.1cm] 0 & \omega^{\scriptscriptstyle H} \end{array}\right) \, \middle| \, \omega \mbox{ skew-Hermitian} \right\rbrace$ &  \\[0.4cm]
\multirow{3}{*}{\parbox{18mm}{Cartan decomposition (\ref{eq:cartan})\,:}}  & \\
& $\mathfrak{p} = \left \lbrace v  = \left(\begin{array}{cc} 0 & s \\[0.1cm] s^{\scriptscriptstyle H} & 0 \end{array}\right) \, \middle| \, 
s \mbox{ complex symmetric} \right\rbrace$ & $T_{\scriptscriptstyle 0\,}\mathbb{D}_{\scriptscriptstyle N}\, \cong \,\mathfrak{p} $, ($0= N \times N$ zero matrix)\\[0.4cm]
\multirow{3}{*}{\parbox{15mm}{Abelian subspace $\mathfrak{a}$\,:}} &  \\
& $\mathfrak{a} =  
\left\lbrace
\left(\begin{array}{cc} 0 & r \\[0.1cm] r & 0 \end{array}\right) \, \middle| \, r \mbox{ real diagonal} \right\rbrace$ &
$r = \mathrm{diag}(r_{\scriptscriptstyle 1},\ldots,r_{\scriptscriptstyle N})$   \\[0.5cm]
Positive roots \,: &  $\lambda_{\scriptscriptstyle ij}(r) = r_{\scriptscriptstyle i} - r_{\scriptscriptstyle j}$ \,;\, $i < j$ & with $ m_{\scriptscriptstyle \lambda_{\scriptscriptstyle ij}}= 1$ \,\phantom{ ,} \, (identically zero when $N = 1$)\\[0.1cm]
&  $\nu_{\scriptscriptstyle ij}(r) = r_{\scriptscriptstyle i} + r_{\scriptscriptstyle j}$ \,;\, $i \leq j$ & with $ m_{\scriptscriptstyle \nu_{\scriptscriptstyle ij}}= 1$\,\phantom{ , }\,\;\, (positive for any value of $N$)\\[0.4cm] 
Killing form \,: & $B(r,r) \, = {\scriptstyle\, \sum_{\scriptscriptstyle i}\, r_{\scriptscriptstyle i}^2}$ &  for $r \in \mathfrak{a}$
\end{tabular}
\vspace{0.3cm}
\end{center}

Indeed, the data in this table can be substituded into formulae (\ref{eq:integralp}) and (\ref{eq:killing}), as these directly apply to the Siegel disc $\mathbb{D}_{\scriptscriptstyle N}$, which is a Riemannian symmetric space of non-compact type, in order to obtain the second line in each of (\ref{eq:btnintegralp}) and (\ref{eq:btnkilling}). \hfill $\blacksquare$
\end{example}

\section{Sampling Gaussian structured covariance matrices} \label{app:sample}
Algorithm 1. of \ref{subsec:gaussiandef} gave a generally applicable method for sampling from a Gaussian distribution on a Riemannian symmetric space of non-positive curvature $\mathcal{M}$. Here, Algorithms 2., 3., and 4. are specific versions of Algorithm 1., which deal with the three spaces of structured covariance matrices considered in Section \ref{sec:matrices}. Precisely, Algorithms 2., 3., and 4. sample Gaussian structured covariance matrices, by implementing the instructions of Algorithm 1., labeled \texttt{l1} to \texttt{l12}. Accordingly, in the description below, each instruction in Algorithms 2., 3., and 4. will be given the same label as the corresponding instruction in Algorithm 1. 
\vspace{0.2cm}

\addtocounter{example}{-3}
\begin{example}
here, Algorithm 2. samples from a Gaussian distribution on the space $\mathcal{M} = \mathcal{H}_{\scriptscriptstyle n\,}$, which was considered in \ref{subsec:2t2}. Precisely, the output of Algorithm 2. is a sample $Y$ from a Gaussian distribution $G(\bar{Y},\sigma)$ with parameters $\bar{Y} \in \mathcal{H}_{\scriptscriptstyle n}$ and $\sigma > 0$. In reference to Algorithm 1., only instructions \texttt{l6} to \texttt{l9} are required. In these instructions, the role of the compact group $H$ is played by the group $U(n)$ of $n \times n$ unitary matrices, and  the role of (\ref{eq:polar}) and (\ref{eq:apdf}) is played by the following (\ref{eq:spectral}) and (\ref{eq:hnapdf}), respectively, 
\begin{subequations} \label{eq:samplehn}
\begin{equation} \label{eq:spectral}
\mbox{spectral decomposition\,:} \hspace{0.45cm} Y(r,U) \, = \, U \, e^{\scriptscriptstyle r} \, U^{\scriptscriptstyle H} \hspace{3.5cm}
\end{equation}
\begin{equation} \label{eq:hnapdf}
p(r) \, \propto \, \exp \left[-\,\frac{\Vert r \Vert^{\scriptscriptstyle 2}}{2\sigma^{\,\scriptscriptstyle 2}}\right]\, \prod_{i < j} \sinh^{\scriptscriptstyle 2}(|r_{\scriptscriptstyle i} - r_{\scriptscriptstyle j}|/2) \,\, \hspace{0.3cm} \mbox{for } r = (r_{\scriptscriptstyle 1},\ldots,r_{\scriptscriptstyle n}) \in \mathbb{R}^{\scriptscriptstyle n}
\end{equation}
In instruction \texttt{l9} below, the notation $\bar{Y}^{\scriptscriptstyle 1/2}\cdot\, Y$ is that of (\ref{eq:congruence}), so that $\bar{Y}^{\scriptscriptstyle 1/2}\cdot\, I \,=\, \bar{Y}$, where $I$ is the $n \times n$ identity matrix, (compare to (ii) of Proposition \ref{prop:sample}).
\end{subequations}

 \vspace{0.3cm}
\hrule
\vspace{0.3cm}

\noindent \texttt{l6}.\; sample $U$ from a uniform distribution on $U(n)$ \\[0.1cm]
\noindent \texttt{l6-a}.\; \; generate $Z$ with $n \times n\,$ i.i.d. complex normal entries   \\[0.1cm]
\noindent \texttt{l6-b}.\; \; perform QR decomposition\,: $Z = UT$ \;\;\% $T$ upper triangular with unit diagonal   \\[0.1cm]
\noindent \texttt{\phantom{l6-b}}\,\; \; \phantom{perform QR decomposition\,: $Z = UT$} \;\;\% see~\cite{eaton} (Proposition 7.3., Page $234$) \\[0.1cm]
\texttt{l7}.\; sample $r$ from the multivariate density (\ref{eq:hnapdf}) \hspace{6cm} {\small \sc Algorithm 2.\,: sampling from a } \\[0.1cm]
\texttt{l8}.\; put $Y = Y(r,U)$ as in (\ref{eq:spectral})    \ \hspace{8.1cm}\,\, {\small \sc Gaussian distribution on $\mathcal{H}_{\scriptscriptstyle n}$} \\[0.1cm]
\texttt{l9}.\; reset $Y$ to $\bar{Y}^{\scriptscriptstyle 1/2}\cdot\, Y$ \hspace{2.3cm} \,\;\;\,\% $\bar{Y}^{\scriptscriptstyle 1/2} = \,$ Hermitian square root of $\bar{Y}$
\vspace{0.25cm}
\hrule
\vspace{0.3cm}
\end{example}
\vspace{0.2cm}

\begin{example} here, Algorithm 3. samples from a Gaussian distribution on the space $\mathcal{M} = \mathcal{T}_{\scriptscriptstyle n}\,$, which was considered in \ref{subsec:tn}. The output of this algorithm is a sample $T$ from a Gaussian distribution $G(\bar{T},\sigma)$ with parameters $\bar{T} \in \mathcal{T}_{\scriptscriptstyle n}$ and $\sigma > 0$. This output sample $T$ is obtained via decomposition (\ref{eq:tnssdecomposition}), in the form
\begin{subequations} \label{eq:sampletn}
\begin{equation} \label{eq:coordinates}
  T \,= \, \Psi(\,r\,,\,\alpha_{\scriptscriptstyle 1}\,, \ldots,\,\alpha_{\scriptscriptstyle n-1}\,) \hspace{0.3cm} r \in \mathbb{R}_{+}\,\mbox{ and }\, 
  \alpha_{\scriptscriptstyle 1}\,, \ldots,\,\alpha_{\scriptscriptstyle n-1} \, \in \mathbb{D}
\end{equation}
where the mapping $\Psi(\,r_{\scriptscriptstyle 0}\,,\,\alpha_{\scriptscriptstyle 1}\,, \ldots,\,\alpha_{\scriptscriptstyle n-1}\,) \, = \, (r_{\scriptscriptstyle 0},r_{\scriptscriptstyle 1}, \ldots,r_{\scriptscriptstyle n-1})$ yields the matrix entries of $T$, given by $T_{\scriptscriptstyle ij} = T^{\scriptscriptstyle *}_{\scriptscriptstyle ji} =  r_{\scriptscriptstyle i-j\,}$. Concretely, $\Psi$ is computed by means of Verblunsky's formula, which maps $\alpha_{\scriptscriptstyle 1}\,, \ldots,\,\alpha_{\scriptscriptstyle j}$ to $r_{\scriptscriptstyle 1}, \ldots,r_{\scriptscriptstyle j}$ , for $j = 1,\ldots, n-1$~\cite{simon} (Theorem 1.5.5., Page $60$).

For Algorithm 3., let $\bar{T} = \Psi(\bar{r},\bar{\alpha}_{\scriptscriptstyle 1},\ldots,\bar{\alpha}_{\scriptscriptstyle n-1})$. Instruction \texttt{l3} of this algorithm samples $r$, given $\bar{r}$ and $\sigma$. Then, instructions \texttt{l6} to \texttt{l9} sample each one of the $\alpha_{\scriptscriptstyle j}$, given $\bar{\alpha}_{\scriptscriptstyle j}$ and $\sigma$. Here, in reference to instructions \texttt{l6} to \texttt{l9} of Algorithm 1., the role of the compact group $H$ is played by the circle group $SU(1)$, and the role of (\ref{eq:polar}) and (\ref{eq:apdf}) is played by the following (\ref{eq:dpolar}) and (\ref{eq:tnapdf}).
\begin{equation} \label{eq:dpolar}
  \mbox{polar coordinates on } \mathbb{D}\mbox{\,:} \hspace{0.5cm} \alpha(\rho,\theta) \, = \,
  \mathrm{tanh}(\rho)\,e^{\scriptscriptstyle \mathrm{i}\theta}
\end{equation}
\begin{equation} \label{eq:tnapdf}
 p(\rho|\,\sigma) \, \propto \, 
 \exp\left[ -\frac{\,\rho^{\scriptscriptstyle 2}}{2\sigma^{\scriptscriptstyle 2}}\right] 
 \, \sinh(|\rho|) 
\end{equation}
In instruction \texttt{l9}, the notation $u\cdot \alpha$ is that of (\ref{eq:tngroup}), and $u_{\scriptscriptstyle j} = u(\alpha_{\scriptscriptstyle j})$, where
\begin{equation} \label{eq:uj}
  u(\alpha) \, = \, \left( \begin{array}{cc} e^{\scriptscriptstyle \mathrm{i}\frac{\theta}{2}} \mathrm{cosh}(\rho)  & 
  e^{\scriptscriptstyle \mathrm{i}\frac{\theta}{2}} \mathrm{sinh}(\rho) \\[0.1cm] 
  e^{\scriptscriptstyle -\mathrm{i}\frac{\theta}{2}} \mathrm{sinh}(\rho) & 
  e^{\scriptscriptstyle -\mathrm{i}\frac{\theta}{2}} \mathrm{cosh}(\rho) \\[0.1cm] 
  \end{array} \right) \hspace{0.5cm} \mbox{for } \alpha = \alpha(\rho,\theta)
\end{equation}
so that $u(\alpha)\cdot 0 = \alpha$. 
\end{subequations}
\vspace{0.3cm}
\hrule
\vspace{0.3cm}

\noindent \texttt{l3}.\; sample $r$ from a log-normal distribution \\ 
\phantom{xxx}\, with median $\bar{r}$ and scale $\sigma/\sqrt{\scriptstyle n}$ 
\hspace{4.3cm} \;\,\% $r = e^{\scriptscriptstyle t}$ where  $t \sim  \mathcal{N}\left(\log \bar{r},\sigma/{\sqrt{\scriptstyle n}\,}\right)$  \\[0.1cm] 
\texttt{\phantom{l4}}.\; \texttt{for} $j = 1,\ldots, n-1$ \texttt{do}\\[0.1cm]
\texttt{l6}.\; \; sample $\theta_{\scriptscriptstyle j}$ from a uniform distribution on $[0,2\pi]$ 
\hspace{2.1cm}\;\;\,\% so $e^{\scriptscriptstyle \mathrm{i}\theta_{\scriptscriptstyle j}} \sim \,$ uniform distribution on $SU(1)$ \\[0.1cm]
\texttt{l7}.\; \; sample $\rho_{\scriptscriptstyle j}$ from the density $p\left(\,\rho_{\scriptscriptstyle j}\,|\,{ \sigma/\sqrt{\scriptstyle n-j}}\,\right)$ given by  (\ref{eq:tnapdf}) \\[0.1cm]
\texttt{l8}.\; \; put $\alpha_{\scriptscriptstyle j} = \alpha(\rho,\theta)$ as in (\ref{eq:dpolar})  \\[0.1cm]
\texttt{l9}.\; \; reset $\alpha_{\scriptscriptstyle j}$ to  $u_{\scriptscriptstyle j}\cdot \alpha_{\scriptscriptstyle j}$ where $u_{\scriptscriptstyle j}$ is given by (\ref{eq:uj})
\hspace{5.5cm} {\small \sc Algorithm 3.\,: sampling from a }
\\[0.1cm]
\texttt{l11}.\! \texttt{end for} \hspace{10cm}\,\,\,\,\,\;\; {\small \sc Gaussian distribution on $\mathcal{T}_{\scriptscriptstyle n}$} \\[0.1cm]
\texttt{l12}.\! put $T \,= \, \Psi(\,r\,,\,\alpha_{\scriptscriptstyle 1}\,, \ldots,\,\alpha_{\scriptscriptstyle n-1}\,)$ as in (\ref{eq:coordinates}) 
\vspace{0.25cm}
\hrule
\vspace{0.3cm}
\end{example}
\vspace{0.2cm}

\begin{example}
here, Algorithm 4. samples from a Gaussian distribution on the space $\mathcal{M} \, = \, \mathcal{T}^{\scriptscriptstyle N}_{\scriptscriptstyle n\,}$, which was considered in \ref{subsec:btn}. The algorithm obtains a sample $T$ from a Gaussian distribution $G(\bar{T},\sigma)$, where $\bar{T} \in \mathcal{T}^{\scriptscriptstyle N}_{\scriptscriptstyle n}$ and $\sigma > 0$, via decomposition (\ref{eq:bttnssdecomposition}). This decomposition is expressed in the form,
\begin{subequations} \label{eq:btnsample}
\begin{equation} \label{eq:btcoordinates}
  T \,= \, \Psi(\,P\,,\,\Omega_{\scriptscriptstyle 1}\,, \ldots,\,\Omega_{\scriptscriptstyle n-1}\,) \hspace{0.3cm} P \in \mathcal{T}_{\scriptscriptstyle N}\,\mbox{ and }\, 
  \Omega_{\scriptscriptstyle 1}\,, \ldots,\,\Omega_{\scriptscriptstyle n-1} \, \in \mathbb{D}_{\scriptscriptstyle N}
\end{equation}
where the mapping $\Psi(P,\Omega_{\scriptscriptstyle 1},\ldots,\Omega_{\scriptscriptstyle n-1}) = (T_{\scriptscriptstyle 0},T_{\scriptscriptstyle 1},\ldots,T_{\scriptscriptstyle n-1})$ yields the $N \times N$ blocks of $T$. Concretely, $\Psi$ is computed as in~\cite{jeuris1} (formula (3.4), Page $6$), which recursively gives the blocks $T_{\scriptscriptstyle 1},\ldots,T_{\scriptscriptstyle j}$ in terms of $\Omega_{\scriptscriptstyle 1},\ldots,\Omega_{\scriptscriptstyle n-1\,}$, for $j = 1,\ldots, n-1$.

For Algorithm 4., let $\bar{T} = \Psi(\bar{P},\bar{\Omega}_{\scriptscriptstyle 1},\ldots,\bar{\Omega}_{\scriptscriptstyle n-1})$. Instruction \texttt{l3} of this algorithm samples $P$, given $\bar{P}$ and $\sigma$. Then, instructions \texttt{l6} to \texttt{l9} sample each one of the $\Omega_{\scriptscriptstyle j}$, given $\bar{\Omega}_{\scriptscriptstyle j}$ and $\sigma$. Here, in reference to instructions \texttt{l6} to \texttt{l9} of Algorithm 1., the role of the compact group $H$ is played by the group $SU(N)$ of unitary $N \times N$ matrices with unit determinant, and the roles of (\ref{eq:polar}) and (\ref{eq:apdf}) are respectively played by (\ref{eq:takagi}) and by the following (\ref{eq:btnapdf}),
\begin{equation} \label{eq:btnapdf}
 p(r|\,\sigma) \, \propto \, 
\exp\left[ -\frac{\,\Vert r \Vert^{\scriptscriptstyle 2}}{2\sigma^{\scriptscriptstyle 2}}\right] 
\, \prod_{\scriptscriptstyle i < j}\sinh(|r_{\scriptscriptstyle i}-r_{\scriptscriptstyle j}|)\prod_{\scriptscriptstyle i \leq j}\sinh(|r_{\scriptscriptstyle i}+r_{\scriptscriptstyle j}|) \hspace{0.3cm} \mbox{for } r = (r_{\scriptscriptstyle 1},\ldots,r_{\scriptscriptstyle N}) \in \mathbb{R}^{\scriptscriptstyle N}
\end{equation}
In instruction \texttt{l9}, the notation $U\cdot \Omega$ is that of (\ref{eq:btngroup}), and $U_{\scriptscriptstyle j} = U(\Omega_{\scriptscriptstyle j})$, where
\begin{equation} \label{eq:btnuj}
  U(\Omega) \, = \, \left(\begin{array}{cc} \Theta\cosh(r) & \Theta\sinh(r) \\[0.1cm] \Theta^{\scriptscriptstyle *}\sinh(r) & \Theta^{\scriptscriptstyle *}\cosh(r) \end{array}\right) \hspace{0.5cm} \mbox{for } \Omega=  \mbox{ right-hand side of (\ref{eq:takagi}) }
\end{equation}
so that $U(\Omega)\cdot 0 = \Omega$, where $0$ denotes the $N \times N$ zero matrix. 
\end{subequations}
\vspace{0.3cm}
\hrule
\vspace{0.3cm}

\noindent \texttt{l3}.\; sample $P$ from a Gaussian distribution $G\left(\bar{P},\sigma/\sqrt{\scriptstyle n}\,\right)$ on $\mathcal{T}_{\scriptscriptstyle N}$ 
\hspace{0.8cm} \;\,\% use Algorithm 3.  \\[0.1cm] 
\texttt{\phantom{l4}}.\; \texttt{for} $j = 1,\ldots, n-1$ \texttt{do}\\[0.1cm]
\texttt{l6}.\; \; sample $\Theta_{\scriptscriptstyle j}$ from a uniform distribution on $SU(N)$ \hspace{2cm} \% compare to \texttt{l6} in Algorithm 2. \\[0.1cm]
\noindent \texttt{l6-a}.\; \; generate $Z$ with $n \times n\,$ i.i.d. complex normal entries   \\[0.1cm]
\noindent \texttt{l6-b}.\; \; perform QR decomposition\,: $Z = \Theta T$  \\[0.1cm]
\noindent \texttt{l6-c}.\; \; reset $\Theta$ to $\left(\det{\Theta}\right)^{\scriptscriptstyle 1/N}\,\Theta$ \\[0.1cm]
\texttt{l7}.\; \; sample $r_{\scriptscriptstyle j}$ from the multivariate density $p(\,r_{\scriptscriptstyle j}\,|\,{ \sigma/\sqrt{\scriptstyle n-j}}\,)$ given by  (\ref{eq:btnapdf}) \\[0.11cm]
\texttt{l8}.\; \; put $\Omega_{\scriptscriptstyle j} = \Omega(r_{\scriptscriptstyle j},\Theta_{\scriptscriptstyle j})$  
\hspace{5.4cm}\;\;\,\% $\Omega(r_{\scriptscriptstyle j},\Theta_{\scriptscriptstyle j}) = \,$ right-hand side of (\ref{eq:takagi}) \\[0.14cm]
\texttt{l9}.\; \; reset $\Omega_{\scriptscriptstyle j}$ to  $U_{\scriptscriptstyle j}\cdot \Omega_{\scriptscriptstyle j}$ where $U_{\scriptscriptstyle j}$ is given by (\ref{eq:btnuj})
\hspace{5.5cm} {\small \sc Algorithm 4.\,: sampling from a }
\\[0.1cm]
\texttt{l11}.\! \texttt{end for} \hspace{10cm}\,\,\,\,\,\;\; {\small \sc Gaussian distribution on $\mathcal{T}^{\scriptscriptstyle N}_{\scriptscriptstyle n}$} \\[0.1cm]
\texttt{l12}.\! put $T \,= \, \Psi(\,P\,,\,\Omega_{\scriptscriptstyle 1}\,, \ldots,\,\Omega_{\scriptscriptstyle n-1}\,)$ as in (\ref{eq:coordinates}) 
\vspace{0.25cm}
\hrule
\vspace{0.3cm}

\end{example}

\end{document}